\input amstex

\documentstyle{amsppt}
\magnification=\magstep1

\hsize=5.2in
\vsize 6.9in
\topmatter
\title W$^*$--SUPERRIGIDITY FOR BERNOULLI ACTIONS\\
 OF PROPERTY (T) GROUPS \endtitle
\rightheadtext{W$^*$--superrigidity for Bernoulli actions}
\vskip 0.05in

\author { ADRIAN IOANA}\endauthor\footnote{The author was supported by a Clay Research Fellowship}
\address Mathematics Department,  UCLA, Los Angeles, CA 91125 and IMAR, Bucharest ROMANIA
\endaddress
\email adiioana\@math.ucla.edu \endemail

\thanks 2010 {\it Mathematics Subject Classification.} Primary 46L36; Secondary 28D15, 37A20.
\endthanks

\abstract 
 We consider group measure space II$_1$ factors $M=L^{\infty}(X)\rtimes\Gamma$ arising from Bernoulli actions 
 of ICC property (T) groups $\Gamma$ (more generally, of groups $\Gamma$ containing an infinite normal subgroup with the relative property (T)) and prove a rigidity result for $*$--homomorphisms $\theta:M\rightarrow M\overline{\otimes}M$.
We deduce that the action $\Gamma\curvearrowright X$ is W$^*$--superrigid, i.e. if $\Lambda\curvearrowright Y$ is {\bf any} free, ergodic, measure preserving action such that the  factors $M=L^{\infty}(X)\rtimes\Gamma$ and $L^{\infty}(Y)\rtimes\Lambda$ are isomorphic,  then the  actions $\Gamma\curvearrowright X$ and $\Lambda\curvearrowright Y$ must be conjugate.
 Moreover, we show that if $p\in M\setminus\{1\}$ is a projection, then $pMp$  does not admit a group measure space decomposition nor a group von Neumann algebra decomposition (the latter under the additional assumption that $\Gamma$ is torsion free).

We also prove a rigidity result for $*$--homomorphisms $\theta:M\rightarrow M$, this time for $\Gamma$ in a larger class of groups than above, now  including products of non--amenable groups. For certain groups $\Gamma$, e.g. $\Gamma=\Bbb F_2\times\Bbb F_2$, we deduce that $M$ does not embed into  $pMp$, for any projection $p\in M\setminus\{1\}$, and obtain a 
description of the endomorphism semigroup of $M$.

\endabstract
\noindent
\endtopmatter
\document

%\vskip -.2in \centerline{ CONTENTS}

\head {Contents}\endhead
\settabs\+\indent\indent & 10. \ & 1.1 \quad & More on rigidity
\cr \tenpoint{\+& 0. & Introduction \cr \+& 1. & Preliminaries  \cr \+& 2. & Rigid subalgebras of II$_1$ factors coming from generalized Bernoulli actions \cr \+& 3. & A spectral gap argument \cr \+& 4. &
Lower bound on height \cr \+& 5. & A conjugacy result for subalgebras \cr \+& 6. &
A dichotomy result for subalgebras \cr \+& 7. & Popa's conjugacy criterion for actions \cr \+& 8. & Strong rigidity for embeddings of II$_1$ factors \cr \+& 9. & W$^*$--superrigidity \cr \+& 10. & Further applications \cr }

\newpage

\head \S 0. {Introduction}.\endhead
\vskip 0.2in
\noindent
{\bf 0.1}\hskip 0.05in The {\it group measure space construction} of Murray and von Neumann  associates to every  probability measure preserving (p.m.p.) action $\Gamma\curvearrowright X$ of a countable group $\Gamma$ on a probability space $(X,\mu)$, a finite von Neumann algebra $L^{\infty}(X)\rtimes\Gamma$ ([MvN36]). If the action is essentially free and ergodic then this algebra is a II$_1$ factor which contains $L^{\infty}(X)$ as a  {\it Cartan subalgebra}. A central problem in the theory of von Neumann algebras is understading how much of the group action is ``remembered'' by its von Neumann algebra.
The main goal of this paper is to prove that, for a large, natural family of group actions (Bernoulli actions of property (T) groups), their associated II$_1$ factor completely remembers the group and the action. 
We start by giving some motivation for this result.

By a celebrated result of A. Connes, if $\Gamma$ is infinite amenable, then the II$_1$ factor of any free ergodic p.m.p. action   $\Gamma\curvearrowright X$  is isomorphic to the hyperfinite II$_1$ factor ([C76], see also [OW80] and [CFW81]). 
In contrast, the study of group measure space algebras  arising from actions of non--amenable groups  has led to a deep rigidity theory (see the survey [Po07b] and the introduction of [PV09]).

In particular, S. Popa's seminal work [Po06ab] shows that, for actions belonging to a large class, 
isomorphism of their group measure space algebras implies conjugacy of the actions.  More precisely, assume that $\Gamma$ is an ICC (infinite conjugacy class) group which has an infinite, normal subgroup with  the relative property (T) and let $\Gamma\curvearrowright X$ be a free ergodic  action. Additionally, suppose that $\Lambda\curvearrowright Y=Y_0^{\Lambda}$ is a  Bernoulli action.  Popa proves that if the actions $\Gamma\curvearrowright X$ and $\Lambda\curvearrowright Y$  are W$^*$--{\it equivalent}, i.e. if they produce isomorphic von Neumann algebras (also known as W$^*$--algebras), then the actions are {\it conjugate} ([Po06b]). This means that $\Gamma\cong\Lambda$ and there exists an isomorphism of probability spaces $\theta:X\rightarrow Y$ such that $\theta\hskip 0.01in\Gamma\hskip 0.01in\theta^{-1}=\Lambda$.
	
The natural question underlying and motivating this result, explicitly formulated in the introduction of [Po06b], is whether the same is true if one imposes all the conditions on only one of the actions. In other words, assume that $\Gamma$ is an ICC group having an infinite, normal subgroup with  the relative property (T) and suppose that $\Gamma\curvearrowright X=X_0^{\Gamma}$ is a Bernoulli action. Is it true that any free ergodic p.m.p. action $\Lambda\curvearrowright Y$ which is W$^*$--equivalent to $\Gamma\curvearrowright X$,  must be conjugate to it? 

Further evidence that the answer to this question is true was provided by Popa. In [Po07a,Po06b], he proves that this is the case if the actions are assumed {\it orbit equivalent} (OE). This amounts to  the existence of an isomorphism of probability spaces $\theta:X\rightarrow Y$ satisfying $\theta(\Gamma x)=\Lambda (\theta x)$, for almost every $x\in X$. 
Recall that orbit equivalence is stronger than W$^*$--equivalence of actions.
More precisely, as first noticed in [Si55] (see also [FM77]), two actions $\Gamma\curvearrowright X$ and $\Lambda\curvearrowright Y$ are orbit equivalent if and only if there exists an isomorphism of their group measure space algebras, $L^{\infty}(X)\rtimes\Gamma\cong L^{\infty}(Y)\rtimes\Lambda$, which identifies the  Cartan subalgebras $L^{\infty}(X)$ and $L^{\infty}(Y)$. 

\vskip 0.05in
\noindent
{\bf 0.2} \hskip 0.05in
The main result of this paper answers affirmatively the above question. Before stating it, we first review a few concepts and then introduce some terminology.

 An inclusion $(\Gamma_0\subset\Gamma)$ of countable groups has the {\it relative property (T) of Kazhdan--Margulis} if any unitary representation of $\Gamma$ which has almost invariant vectors must have a non--zero $\Gamma_0$--invariant vector ([Ka67],[Ma82]). When $\Gamma_0=\Gamma$, this condition is equivalent to the {\it property (T)} of the group $\Gamma$. Examples of
relative property (T) inclusions of groups are given by $(\Bbb Z^2\subset\Bbb Z^2\rtimes\Gamma)$, for any non--amenable
subgroup $\Gamma$ of SL$_2(\Bbb Z)$ ([Bu91]), and by $(\Gamma_0\subset\Gamma_0\times\Gamma_1)$, for a property (T) group $\Gamma_0$ 
(e.g. $\Gamma_0=$ SL$_n(\Bbb Z)$, $n\geq 3$) and an arbitrary countable group $\Gamma_1$ ([Ka67]).

For a probability space $(X_0,\mu_0)$, the {\it Bernoulli action} $\Gamma\curvearrowright (X_0^{\Gamma},\mu_0^{\Gamma})$ is given by $\gamma\cdot (x_g)_g=(x_{\gamma^{-1}g})_g$, for all $(x_g)_{g\in\Gamma}\in X_0^{\Gamma}$ and $\gamma\in\Gamma$.

\vskip 0.05in
 In this paper, by a W$^*$-- or OE--{\it rigidity result} we mean a result proving that two actions $\Gamma\curvearrowright X,$ $\Lambda\curvearrowright Y$, which are W$^*$-- or OE--equivalent, must be conjugate. If this happens when only one of these actions is in a fixed class, while the other can be {\it any} free ergodic action of {\it any} countable group, we have a {\it superrigidity} result.

\proclaim {Theorem A (W$^*$--superrigidity)} Let $\Gamma$ be a countable ICC  group which admits an infinite normal subgroup $\Gamma_0$ such that the inclusion $(\Gamma_0\subset\Gamma)$ has the relative property (T). Let $(X_0,\mu_0)$ be a non--trivial probability space and let $\Gamma\curvearrowright(X,\mu)=(X_0^{\Gamma},\mu_0^{\Gamma})$ be the Bernoulli action. Denote $M=L^{\infty}(X)\rtimes\Gamma$ and let $p\in M$ be a projection. 
Let $\Lambda\curvearrowright (Y,\nu)$ be a free ergodic p.m.p. action of a countable group $\Lambda$. Denote $N=L^{\infty}(Y)\rtimes\Lambda$.
\vskip 0.05in
\noindent
If $N\cong pMp$, then $p=1$, $\Gamma\cong \Lambda$ and the actions $\Gamma\curvearrowright X$, $\Lambda\curvearrowright Y$ are conjugate. 

\endproclaim
Moreover,  any $*$--isomorphism $\theta:N\rightarrow M$ comes from a conjugacy of the actions $\Gamma\curvearrowright X$, $\Lambda\curvearrowright Y$ and a character of $\Gamma$ (see Theorem 9.1). Note that Theorem A provides a new, large class of II$_1$ factors which are not group measure space factors.
\vskip 0.05in
Theorem A is proven in the framework of Popa's {\it deformation/rigidity theory} by playing against each other the ``rigidity" of $\Gamma$ (manifested here in the form of relative property (T)) and the ``deformation properties" of Bernoulli actions $\Gamma\curvearrowright X$ (see 1.5). 
Before discussing its method of proof in more detail,  let us put into context.	 

\vskip 0.05in
\noindent
{\bf 0.3} \hskip 0.05in
Despite remarkable progress in both the areas of  W$^*$-- and OE--rigidity during the last decade (see [Po07b], [Fu09], [PV09]), W$^*$--superrigidity results remained elusive until very recently ([Pe10], [PV09]). Prior to these results, several large classes of OE--superrigid actions were found ([Fu99], [Po07a], [Po08], [Ki06], [Io08], [Ki09]). In addition, for other classes of actions $\Gamma\curvearrowright X$, it was shown that $L^{\infty}(X)$ is the unique group measure space Cartan subalgebra of  $L^{\infty}(X)\rtimes\Gamma$, up to unitary conjugacy ([OP07],[Pe10]). However, since no examples of actions having both of these properties were known,  W$^*$--superrigidity still could not be concluded.

The situation changed starting with the work of J. Peterson who was able to show the existence of virtually W$^*$--superrigid actions ([Pe10]). He considered profinite actions of certain groups having infinite subgroups with the relative property (T) and used, among other things, techniques and results from [Pe09], [OP07], [Io08].

Shortly after, S. Popa and S. Vaes discovered the first concrete families of W$^*$--superrigid actions ([PV09]). Furthermore, they succeeded to prove sweeping W$^*$--superrigidity results. 
First, they  showed that for groups $\Gamma$ in a certain class $\Cal G$ of amalgamated free product groups, the II$_1$ factor of any free ergodic action $\Gamma\curvearrowright X$ has a unique group measure space Cartan subalgebra.
Then, by applying the OE--superrigidity results from  [Po07a, Po08] and   [Ki09] they respectively deduced that the following actions are W$^*$--superrigid: a) Bernoulli actions, generalized Bernoulli actions, Gaussian actions  of groups $\Gamma\in\Cal G$   and b) any free mixing action of $\Gamma=$ PSL$_n(\Bbb Z)*_{T_n}\text{PSL}_n(\Bbb Z)\in\Cal G$, where $T_n\subset$ PSL$_n(\Bbb Z)$ is the subgroup of upper triangular matrices and $n\geq 3$. 

The class $\Cal G$ consists of groups of the form $\Gamma=\Gamma_1*_{\Sigma}\Gamma_2$ subject to several conditions: 1) $\Sigma$  is  amenable, 2) $\Sigma$ is  ``almost malnormal" in $\Gamma_1$, normal in $\Gamma_2$ and 3) $\Gamma_1$ is ``rigid" (e.g. has property (T) or is a product of non--amenable groups). As is easy to see, 2) precludes groups $\Gamma\in\Cal G$ from having an infinite normal subgroup with the relative property (T). Therefore, Theorem A and the results of [PV09] cover disjoint classes of groups. In particular, Theorem A is the first W$^*$--superrigidity result that applies to property (T) groups.

\vskip 0.05in
\noindent{\bf 0.4}  \hskip 0.05in
In this paper, we introduce a general strategy for analyzing group measure space decompositions of II$_1$ factors and  implement it successfully to
derive Theorem A.
To outline our  approach, let  $M$ be the II$_1$ factor associated to some fixed, ``known" action $\Gamma\curvearrowright X$. 
The first part consists of ``classifying" all unital  $*$--homomorphisms $\theta:M\rightarrow M\overline{\otimes}M$ (or {\it embeddings} of $M$ into $M\overline{\otimes}M$) in terms of the decomposition $M=L^{\infty}(X)\rtimes\Gamma$ (i.e. one would like to have a rigidity result concerning the structure of such $\theta$ (see Thm. C)).

Now, suppose that $M$ also arises as the II$_1$ factor of a ``mystery" action $\Lambda\curvearrowright Y$. 
The decomposition of $M$ as $L^{\infty}(Y)\rtimes\Lambda$ induces a unital $*$--homomorphism $\theta:M\rightarrow M\overline{\otimes}M$ given by $\theta(av_{\lambda})=av_{\lambda}\otimes v_{\lambda}$, for all $a\in L^{\infty}(Y)$ and  every $\lambda\in\Lambda$. 
This embedding of $M$ into $M\overline{\otimes}M$ has been introduced and used by Popa and Vaes in the proof of [PV09, Lemma 3.2.].
In the second part, we apply the above classification to $\theta$. 
Thus, we have some information about the form of $\theta$ with respect to both group measure space decompositions of $M$. This relates the two decompositions and, ideally, the relationship will be powerful enough to imply that the two decompositions coincide, up to unitary conjugacy.
If in addition we know that the action $\Gamma\curvearrowright X$ is OE--superrigid then we can deduce that the involved actions must be conjugate.

\vskip 0.05in
\noindent
{\it Remark.} Using this strategy it can be readily seen that if a II$_1$ factor $M$ has property (T), in the sense of Connes--Jones [CJ85], then it  admits only countably many group measure space decompositions (see Section 10).

\vskip 0.05in

\noindent
{\bf 0.5} \hskip 0.05in
 Following the above strategy, we prove that the II$_1$ factor $M=L^{\infty}(X)\rtimes\Gamma$ defined in Theorem A has a unique  group measure space Cartan subalgebra. In combination with Popa's OE superrigidity result [Po07a] this proves Theorem A.

To apply our approach, we first need to describe all $*$--homomorphisms $\theta:M\rightarrow M\overline{\otimes}M$. We start by addressing the seemingly easier question of describing all $*$--homomorphisms $\theta:M\rightarrow M$.  
In this context, we prove a strong rigidity result in the spirit of [Po06b]. 
 
\proclaim {Theorem B} Let $M=L^{\infty}(X)\rtimes\Gamma$ be as in Theorem A.
Let $\theta:M\rightarrow M$ be a (not necessarily unital)  $*$--homomorphism and assume that  $\Gamma$ is torsion free. 
\vskip 0.05in
\noindent
Then one of the following holds true:

\noindent
(1) $\theta(M)\prec_{M}L(\Gamma)$ or

\noindent
(2) $\theta$ is unital and we can find  a character $\eta$ of $\Gamma$, a morphism 
$\delta:\Gamma\rightarrow\Gamma$ and a unitary $u\in M$ such that  
 $u\theta(L^{\infty}(X))u^*\subset L^{\infty}(X)$ and 
 $u\theta(u_{\gamma})u^*=\eta(\gamma)u_{\delta(\gamma)},$ for all $\gamma\in\Gamma$.
\endproclaim

Furthermore, generalizing the main result of [Po06b] (in the case of {\it w--rigid} groups),  we prove a rigidity result for $*$--homomorphisms between two II$_1$ factors, each of which share some, but not all, of the defining properties of $M$  (see Theorem 8.1). 
\vskip 0.01in
Before elaborating on conditions (1) and (2), let us recall the notations used in Theorem B.
Thus, we denote  by $L(\Gamma)\subset M$ the copy of the group von Neumann algebra of $\Gamma$ generated by the canonical unitaries $\{u_{\gamma}\}_{\gamma\in\Gamma}\subset M$.  Also, we use the notation $Q\prec_{M}B$ to indicate that ``a corner of a subalgebra $Q\subset M$ can be embedded into a subalgebra $B\subset M$ inside a factor $M$'',  in the sense of Popa ([Po06a], see 1.3.1). This roughly means that we can conjugate $Q$ into $B$ with a unitary element from $M$. 

Now, condition (1) essentially says that $\theta$ comes from an embedding of $M$ into $L(\Gamma)$. 
To better explain condition (2), we rephrase it differently  using the following example. 

\vskip 0.05in
\noindent
{\it Example.} Quotient actions naturally induce $*$--homomorphisms between  group measure space algebras. 
Recall that  an action $\Gamma\curvearrowright Y$ is a {\it quotient} (or a {\it factor}) of  an action $\Gamma\curvearrowright Z$ if there exists a $\Gamma$--equivariant,  measure preserving map $q:Z\rightarrow Y$. This is equivalent to the existence of a $*$--homomorphism $\theta:L^{\infty}(Y)\rtimes\Gamma\rightarrow L^{\infty}(Z)\rtimes\Gamma$ such that $\theta(L^{\infty}(Y))\subset L^{\infty}(Z)$ and $\theta(u_{\gamma})=u_{\gamma}$, for all $\gamma\in\Gamma$.
\vskip 0.05in

 With this in mind, (2) says that there exists a subgroup $\Gamma'$ of $\Gamma$ isomorphic to $\Gamma$ such that $\theta$ comes from a conjugacy between $\Gamma\curvearrowright X$ and a  quotient of the action $\Gamma'\curvearrowright X$.

\vskip 0.05in

\noindent
{\bf 0.6}\hskip 0.05in Next, returning to our initial question of describing  $*$--homomorphisms $\theta:M\rightarrow M\overline{\otimes}M$, we prove the following classification result.

\proclaim {Theorem C}  Let $M=L^{\infty}(X)\rtimes\Gamma$ be as in Theorem A.
Let $\theta:M\rightarrow M\overline{\otimes}M$ be a (not necessarily unital) $*$--homomorphism and assume that $\Gamma$ is torsion free.
\vskip 0.05in
\noindent
Then one of the following holds true:

\noindent
(1) $\theta(L(\Gamma_0))\prec_{M\overline{\otimes}M}L(\Gamma)\otimes 1$ or $\theta(L(\Gamma_0))\prec_{M\overline{\otimes}M}1\otimes L(\Gamma)$.

\noindent
(2) $\theta(M)\prec_{M\overline{\otimes}M}L(\Gamma)\overline{\otimes}M$ or $\theta(M)\prec_{M\overline{\otimes}M}M\overline{\otimes}L(\Gamma)$.

\noindent
(3)  $\theta$ is unital  and we can find  a character $\eta$ of $\Gamma$, two group morphisms 
$\delta_1,\delta_2:\Gamma\rightarrow\Gamma$ and a unitary $u\in M$ such that  
$u\theta(L^{\infty}(X))u^*\subset L^{\infty}(X)\overline{\otimes}L^{\infty}(X)$ and 
 $u\theta(u_{\gamma})u^*=\eta(\gamma)(u_{\delta_1(\gamma)}\otimes u_{\delta_2(\gamma)})$, for all $\gamma\in\Gamma$.
\endproclaim

This result holds true without the torsion freeness assumption, but in this case its statement becomes considerably more complicated (see Theorem 8.2). Therefore, for simplicity, we assume only here that  $\Gamma$ is torsion free although we will later state and prove Theorem C in full generality. This generalization  will enable us to prove Theorem A for groups with torsion, including the interesting examples SL$_n(\Bbb Z)$ ($n\geq 3$) and  $\Bbb Z^2\rtimes$SL$_2(\Bbb Z)$. 

\vskip 0.05in

\noindent
{\bf 0.7} \hskip 0.05in  In Section 9, we will combine Theorem C and the strategy described in 0.4 to deduce Theorem A. Since the argument needed is quite involved, we will not elaborate more on this here.
Rather, we mention a few words on the proofs of Theorem B and C. Since these proofs are analogous, let us only discuss the simpler case of Theorem B.

 Consider, for simplicity, a unital embedding $\theta$ of $M=L^{\infty}(X)\rtimes\Gamma$ into itself.
The starting point of the proof of Theorem B is Popa's ``absorption" result for relative property (T) subalgebras of $M$ which enables us to assume that $\theta(L(\Gamma))\subset L(\Gamma)$ ([Po06a]). 

The main ingredient of the proof consists of the following dichotomy  for abelian subalgebras of II$_1$ factors coming from Bernoulli actions.

\proclaim {Theorem D} Let  $\Gamma\curvearrowright (X,\mu)$ be a Bernoulli action of a countable group $\Gamma$. Denote $M=L^{\infty}(X)\rtimes\Gamma$ and let $D\subset M$ be a unital abelian von Neumann subalgebra. 

\noindent
Suppose that there exists a sequence of unitaries $\{u_n\}_{n\geq 1}\subset L(\Gamma)$ such that $u_n\rightarrow 0$ (weakly)
 and  $u_n$ normalizes $D$ (i.e. $u_nDu_n^*=D$), for all $n\geq 1$.

\vskip 0.05in
\noindent
 Then either: 

\noindent
(1) $D\prec_{M}L(\Gamma)$ or

\noindent
(2) $D'\cap M$ is of type I and there exists a unitary $u\in M$ such that $uL^{\infty}(X)u^*\subset D'\cap M$.
\endproclaim

The proof of Theorem D is based on several new ideas as well as on an idea of J. Peterson ([Pe10], see Section 3). 
Let us emphasize two key tools that we introduce, both of which are related to a certain notion of ``height".

First, fix an element $x$ in a group von Neumann algebra $L(\Gamma)$ and write $x=\sum_{\gamma\in\Gamma}x_{\gamma}u_{\gamma}$, where $x_{\gamma}\in\Bbb C$. The {\it height} of $x$ is defined as $h(x)=\max_{\gamma\in\Gamma}|x_{\gamma}|$  and measures the distance between $x$ and the group $\Gamma=\{u_{\gamma}\}_{\gamma\in\Gamma}\subset L(\Gamma)$. With this notation, we show that if condition (1) of Theorem D fails, then $\inf_{n}h(u_n)>0$. Roughly, this means that the $u_n$'s are uniformly close to $\Gamma$. This fact is crucial  because it allows us to carry subsequent calculations by analogy with the case when $u_n\in\Gamma$.

More generally, if $x=\sum_{g\in\Gamma}x_gu_g$
is an element of a  crossed product algebra  $M= A\rtimes\Gamma$, then we define its  {\it height over} $A$ as  $h_{A}(x)=\max_{g\in\Gamma} ||x_g||_2$. Our second tool is a new ``intertwining theorem" for subalgebras $D\subset M$. Thus, we prove that if the height function $h_A$ is bounded from below on the unitary group of $D$, i.e. $\inf_{x\in \Cal U(D)}h_A(x)>0,$ then a corner of $D$ can be embedded into $A$ inside $M$ (see Theorem 1.3.2).

Returning to the explanation of Theorem B, fix a sequence $\{\gamma_n\}_{n\geq 1}\subset\Gamma$  with $\gamma_n\rightarrow\infty$. Notice that the unitary elements $u_n=\theta(u_{\gamma_n})$ belong to $L(\Gamma)$ and normalize $D=\theta(L^{\infty}(X))$. By applying Theorem D we derive that either (1) $D\prec_{M}L(\Gamma)$ or (2) $D'\cap M=uL^{\infty}(X)u^*$, for some unitary $u$ (this is true only up to finite index).

 In the first case, since  $\theta(L(\Gamma))\subset L(\Gamma)$, we deduce that $\theta(M)\prec_{M}L(\Gamma)$, thus condition (1) in the conclusion of Theorem B is satisfied. 

In the second case, denote by $\tilde\Gamma$ the group $\{\theta(u_{\gamma})\}_{\gamma\in\Gamma}$. Then $\tilde\Gamma$  normalizes  $\tilde D=D'\cap M$. The crossed product algebra $N=\tilde D\rtimes\tilde\Gamma$, viewed inside $M=L^{\infty}(X)\rtimes \Gamma$, has the property that  its core, $\tilde D$, is a unitary conjugate of $L^{\infty}(X)$ in addition to satisfying $L(\tilde\Gamma)\subset L(\Gamma)$. A generalization of Popa's conjugacy criterion for actions ([Po06b], see Section 7) implies that we may in fact assume that $\tilde\Gamma\subset\Gamma$ (modulo scalars). In other words, condition (2) in the conclusion of Theorem B holds true.

\vskip 0.05in
\noindent {\bf 0.8}\hskip 0.05in
We continue with an application of Theorem B. 

\proclaim {Corollary E}   Let $\Gamma\curvearrowright (X,\mu)$ be a non--trivial Bernoulli action of $\Gamma=\Bbb F_m\times\Bbb F_n$, with $2\leq m,n\leq\infty$. Denote $M=L^{\infty}(X)\rtimes\Gamma$.
\vskip 0.05in		
\noindent
If $\theta:M\rightarrow pMp$ is a unital $*$--homomorphism, for some projection $p\in M$, then $p=1$. Moreover, there exist a character $\eta$ of $\Gamma$, a group morphism $\delta:\Gamma\rightarrow\Gamma$ and a unitary $u\in M$ such that $u\theta(L^{\infty}(X))u^*\subset L^{\infty}(X)$ and $u\theta(u_{\gamma})u^*=\eta(\gamma)u_{\delta(\gamma)}$, for all $\gamma\in\Gamma$.
\endproclaim

The proof of Corollary E  is the combination of two facts: an extension of Theorem B to groups $\Gamma$ that are products of non--amenable groups (see Theorem 8.1) and  a result of N. Ozawa and S. Popa  guaranteeing that, if $\Gamma$ is non--amenable and has the complete metric approximation property (e.g. if $\Gamma=\Bbb F_m\times\Bbb F_n$), then  there is no embedding of $M$ into $L(\Gamma)$ ([OP07]).

Corollary E provides the first examples of II$_1$ factors $M$ which do not embed into any of their corners. (The question of whether such factors exist was posed by N. Ozawa during a seminar at UCLA in 2007.) 
Corollary E also reduces the 
calculation of the endomorphism semigroup of $M$ (i.e. the semigroup, End$(M)$, of unital $*$--homomorphisms $\theta:M\rightarrow M$) to an ergodic--theoretic problem. More precisely, any $\theta\in$ End$(M)$ is determined by the following data: a character $\eta$ of $\Gamma$, a morphism $\delta:\Gamma\rightarrow\Gamma$ and  a measure preserving map $q:X\rightarrow X$ satisfying $q(\delta(\gamma)x)=\gamma q(x)$, for every $\gamma\in\Gamma$ and almost all $x\in X$ (see the example from 0.5).

 \vskip 0.05in
\noindent {\bf 0.9}\hskip 0.05in 
Finally, let us present an application of Theorem C.  

After introducing the group measure space construction in [MvN36],  Murray and von Neumann later found a simpler way of constructing II$_1$ factors. Thus, to every countable group $\Lambda$, one can associate its {\it group von Neumann algebra}, $L(\Lambda)$  ([MvN43]). This algebra is finite in general and it is a II$_1$ factor if and only if $\Lambda$ is ICC. 

The first examples of II$_1$ factors which are {\it not} group von Neumann algebras were discovered by A. Connes  ([C75], see also [Jo80]). Recently, more examples have been exhibited in [IPP08] and [PV08a]. All of these examples were obtained by analyzing anti--automorphisms. More precisely, one  first proves that the II$_1$ factors involved either do not have anti--automorphisms ([C75],[IPP08],[PV08a]) or do not have anti--automorphisms of order 2 ([Jo80]). Then, since any group von Neumann algebra  $L(\Lambda)$ admits an anti--automorphism of order 2 (given by $\Phi(v_{\lambda})=v_{\lambda^{-1}}$, for $\lambda\in\Lambda$), one  deduces lack of group von Neumann algebra decomposition.

As a consequence of Theorem C, we obtain a wide class of new examples through a different method.

\proclaim {Corollary F} Let $=L^{\infty}(X)\rtimes\Gamma$ be as in Theorem A. Assume that $\Gamma$ is torsion free.
\vskip 0.05in
\noindent
 Then, for any projection $p\in M\setminus\{0,1\}$, the II$_1$ factor $pMp$ is not isomorphic to the group von Neumann algebra,  $L(\Lambda)$, of a countable group $\Lambda$.  
\endproclaim

For a more general statement, see Theorem 10.1. 

  Corollary F provides in particular the first examples of II$_1$ factors which  have involutary anti--automorphisms and yet are not group von Neumann algebras (see Remark 10.3). This answers a question posed by V.F.R. Jones in [Jo80, Remark 5.7.].

The proof of Corollary F borrows its main idea from the strategy described in 0.4. If $pMp=L(\Lambda)$,  then $\theta:pMp\rightarrow pMp\overline{\otimes}pMp$ given by $\theta(v_{\lambda})=v_{\lambda}\otimes v_{\lambda}$ is a unital $*$--homomorphism. By taking amplifications, we get a  $*$--homomorphism $\tilde\theta:M\rightarrow M\overline{\otimes}M$ with $\tau(\tilde\theta(1))=\tau(p)$. It is immediate to see that $\tilde\theta$ does not verify conditions (1) and (2) of Theorem C.
Therefore, $\tilde\theta$ is forced to be unital and hence $p=1$. 
\vskip 0.05in
Furthermore, as we explain in Section 10 (IV), an extension of our results can be used to show
 that the II$_1$ factors $pMp$ from Corollary F are in fact not isomorphic to any {\it twisted} group von Neumann algebra $L_{\alpha}(\Lambda)$.

\proclaim {Corollary G} Let $=L^{\infty}(X)\rtimes\Gamma$ be as in Theorem A. Assume that $\Gamma$ is torsion free.
\vskip 0.05in
\noindent
 Then, given a projection $p\in M\setminus\{0,1\}$, the II$_1$ factor $pMp$ is not isomorphic to  $L_{\alpha}(\Lambda)$, for any countable group $\Lambda$ and any 2--cocycle $\alpha\in$ H$^2(\Lambda,\Bbb T)$.  
\endproclaim

\vskip 0.1in
\noindent
{\it Organization of the paper.} Besides the introduction, this paper has ten other sections. In Section 1, we recall a few notions and results regarding von Neumann algebras and establish a crucial intertwining result (Theorem 1.3.2). In Section 2, we prove an ``absorption" result for relative property (T) subalgebras of II$_1$ factors coming from generalized Bernoulli actions. The proof of Theorem D occupies Sections 3--6. Further, in Section 7, we generalize Popa's conjugacy criterion for actions. In Section 8, we combine the results of the previous sections to derive Theorems B and C, while in Section 9 we deduce Theorem A. Our last section is devoted to several applications of Theorems B and C (including Corollaries E, F and G).

\vskip 0.1in
\noindent
{\it Acknowledgments.} I would like to thank Jesse Peterson for explaining [Pe10] and  Sorin Popa for many useful suggestions. I am also grateful to Vaughan Jones for bringing to my attention the question that led to Corollary G, Stefaan Vaes for providing me with an argument which simplifies the initial proof of Theorem 9.1 and the anonymous referee for several suggestions which helped improve the exposition. This paper was written while I was visiting the Department of Mathematics at UCLA.
\vskip 0.2in

\head \S 1. {Preliminaries.}\endhead

\vskip 0.2in
\noindent {\bf 1.1 Finite von Neumann algebras}. We first review a few concepts and constructions involving von Neumann algebras (see e.g. [BrOz08]).  In this paper we work with finite von Neumann algebras $M$ always endowed with a fixed faithful normal tracial state $\tau$. 
We denote by  $L^2(M)$ the Hilbert space obtained by completing  $M$ with respect to the 2--norm
$||x||_2=\tau(x^*x)^{\frac{1}{2}}$. Hereafter, we see every $x\in M$ both as an element of $L^2(M)$ and as a bounded (left multiplication) operator on $L^2(M)$.  The inequality $||axb||_2\leq ||a||\hskip 0.02in||x||_2\hskip 0.02in||b||$, for all $a,b,x\in M$, will be used often.

We denote by $\Cal U(M)$ the group of unitaries of $M$, by Aut$(M)$ the group of automorphisms of $M$ endowed with the pointwise $||.||_2$ topology and by id$_{M}$ the identity automorphism of $M$. For every $u\in\Cal U(M)$, the inner automorphism Ad$(u)$ of $M$ is given by Ad$(u)(x)=uxu^*$. We also denote by $\Cal P(M)$ the set of projections of $M$, by $\Cal Z(M)$ the {\it center} of $M$ and by $(M)_1$  the set of $x\in M$ with $||x||\leq 1$. A Hilbert space $\Cal H$ is called an {\it $M$--bimodule} if it is endowed with commuting left and right Hilbert $M$--module structures. 
 For a subset $X$ of $M$, we denote by $L^2(X)$ its closure  in $L^2(M)$, by $X'\cap M$ its  commutant in $M$ and by $X''$  the von Neumann algebra it generates.

\vskip 0.05in

 Let $B$ be a unital von Neumann subalgebra of $(M,\tau)$.  We denote by $ qN_{M}(B)$ the {\it quasi--normalizer} of $B$ in $M$, i.e. the set of $x\in M$ for which there exist $x_1,..,x_n\in M$ satisfying $xB\subset \sum_{i=1}^nBx_i$, $Bx\subset\sum_{i=1}^nx_iB$ (see [Po06c]). 
Note that $qN_{M}(B)$ contains $\Cal N_{M}(B)$, the {\it normalizer} of $B$ in $M$, i.e. the set of unitaries $u\in M$ such that Ad$(u)(B)=B$. If $q\Cal N_{M}(B)''=M$, we say that $B$ is {\it quasi--regular} in $M$ and if $\Cal N_{M}(B)''=M$, we say that $B$ is {\it regular} in $M$.  If $\Gamma_0\subset\Gamma$ is an {\it almost normal}  inclusion of countable groups, then $L(\Gamma_0)$ is quasi--regular in $L(\Gamma)$.

Recall that if $e_B$ is the orthogonal projection from $L^2(M)$ onto $L^2(B)$, then its restriction to $M$ is equal to $E_B$, the conditional expectation from $M$ onto $B$. Jones' {\it  basic construction} $\langle M,e_B\rangle$ is the von Neumann algebra generated by $M$ and $e_B$ inside $\Bbb B(L^2(M))$. The basic construction $\langle M,e_{B}\rangle$ contains the span of $\{xe_By|x,y\in M\}$ as a dense $*$--subalgebra and is endowed with a faithful normal semi--finite trace $Tr$ given by $Tr(xe_By)=\tau(xy)$. We denote by $L^2(\langle M,e_B\rangle)$ the associated Hilbert space.
\vskip 0.05in

Next, let $\omega$ be a free ultrafilter on $\Bbb N$.
For a finite von Neumann algebra $(M,\tau)$, we denote by $(M^{\omega},\tau_{\omega})$ its {\it ultrapower} algebra, i.e. the finite von Neumann algebra $\ell^{\infty}(\Bbb N,M)/\Cal I$, where the trace $\tau_{\omega}$ is given by $\tau_{\omega}((x_n)_n)=\lim_{n\rightarrow\omega}\tau(x_n)$ and $\Cal I$ is the ideal of $x=(x_n)_n\in \ell^{\infty}(\Bbb N,M)$ such that $\tau_{\omega}(x^*x)=0$.  Note that $M$ embeds into $M^{\omega}$ through the map $x\rightarrow (x_n)_n,$ where $x_n=x$, for all $n$, and that any automorphism $\theta$ of $M$ extends to an automorphism of $M^{\omega}$ by letting $\theta((x_n)_n)=(\theta(x_n))_n$, for all $x=(x_n)\in M^{\omega}$.

\vskip 0.05in

Finally, let $M$ be a II$_1$ factor and  $t>0$. Let $n\geq t$ be an integer and  $p\in\Bbb M_n(\Bbb C)\otimes M$ be a projection of normalized trace  $\frac{t}{n}$. The isomorphism class of the algebra $p(\Bbb M_n(\Bbb C)\otimes M)p$ is independent of the choice of $n$ and $p$, is called the {\it $t$--amplification of} $M$ and is denoted by $M^t$. If $N$ is a II$_1$ factor and $\theta:N\rightarrow M$ is a unital $*$--homomorphism, then for every $t>0$ there exists a natural $*$--homomorphism $\theta^t:N^t\rightarrow M^t$. Moreover, $\theta^t$ in uniquely defined, up to composition with a inner automorphism.

\vskip 0.1in
\noindent {\bf 1.2 The crossed product construction}. Let $\sigma:\Gamma\rightarrow$ Aut$(B)$ be a trace preserving action a countable group $\Gamma$ on a finite von Neumann algebra $(B,\tau)$ ([MvN36]). Set $\Cal H=L^2(B)\overline{\otimes}\ell^2(\Gamma)$ and for every $b\in B$, $\gamma\in\Gamma$ define the operators $L_f,u_{\gamma}\in\Bbb B(\Cal H)$ through the formulas $$L_b(b'\otimes \delta_{\gamma'})=bb'\otimes\delta_{\gamma'},\hskip 0.1in u_{\gamma}(b'\otimes \delta_{\gamma'})=\sigma(\gamma)(b')\otimes \delta_{\gamma\gamma'},\hskip 0.1in \forall b'\in L^2(B),\gamma'\in\Gamma.$$

Since $u_{\gamma}u_{\gamma'}=u_{\gamma\gamma'},L_bL_{b'}=L_{bb'},u_{\gamma}L_bu_{\gamma}^*=L_{\sigma(\gamma)(b)}$, for every $\gamma,\gamma'\in\Gamma$ and $b,b'\in B$, the linear span of $\{L_bu_{\gamma}|b\in B,\gamma\in\Gamma\}$ is a $*-$subalgebra of $\Bbb B(\Cal H)$. The strong operator closure of this algebra, denoted $B\rtimes_{\sigma}\Gamma$, is called the {\it crossed product von Neumann algebra} associated to  $\sigma$. The trace $\tau$ extends to a trace on $B\rtimes_{\sigma}\Gamma$ given by $\tau(bu_{\gamma})=\delta_{\gamma,e}\tau(b)$, making $B\rtimes_{\sigma}\Gamma$ is a finite von Neumann algebra. Every element $x\in B\rtimes_{\sigma}\Gamma$ decomposes as $x=\sum_{\gamma\in\Gamma}b_{\gamma}u_{\gamma}$, where the convergence holds in $||.||_2$ and $b_{\gamma}\in B$, for every $\gamma\in\Gamma$. 
\vskip 0.05in

Two important examples of crossed product algebras arise when $B$ is either trivial or abelian.
 Firstly, if $B=\Bbb C1$ (with $\Gamma$ acting trivially), then the associated crossed product algebra is precisely {\it the group von Neumann algebra} $L(\Gamma)$ of $\Gamma$ ([MvN43]). In general, we have the natural embedding $L(\Gamma)\cong\{u_{\gamma}|\gamma\in\Gamma\}''\subset B\rtimes_{\sigma}\Gamma$.
Secondly, if $B=L^{\infty}(X)$, for some standard probability space $(X,\mu)$, then $\sigma$ comes from a probability measure preserving (p.m.p.) action $\Gamma\curvearrowright^{\sigma} (X,\mu)$. The crossed product algebra $L^{\infty}(X)\rtimes_{\sigma}\Gamma$ is called  the {\it group measure space construction} associated with $\sigma$ ([MvN36]). If $\Gamma$ is infinite and $\sigma$ is (essentially) free and ergodic, then $L^{\infty}(X)\rtimes_{\sigma}\Gamma$ is a II$_1$ factor and $L^{\infty}(X)$ is a {\it Cartan subalgebra}, i.e. it is regular and maximal abelian.
\vskip 0.05in

Two free, ergodic p.m.p. actions $\Gamma\curvearrowright (X,\mu)$ and $\Lambda\curvearrowright (Y,\nu)$ are said to be
\vskip 0.03in
\noindent
$\bullet$ {\it conjugate} if there exist a measure space isomorphism $\theta:X\rightarrow Y$ and a group isomorphism $\delta:\Gamma\rightarrow\Lambda$ such that $\theta(\gamma x)=\delta(\gamma)\theta(x)$, for almost every $x\in X$,

\noindent
$\bullet$ {\it orbit equivalent} if there exists a measure space isomorphism $\theta:X\rightarrow Y$ such that $\theta(\Gamma x)=\Lambda\theta(x)$, for almost every $x\in X$, and

\noindent
$\bullet$ W$^*$--{\it equivalent} (or {\it von Neumann equivalent}) if $L^{\infty}(X)\rtimes\Gamma\cong L^{\infty}(Y)\rtimes\Lambda.$
\vskip 0.05in

If two actions are conjugate, then they are orbit equivalent. In turn, two actions are orbit equivalent if and only if the inclusions $(L^{\infty}(X)\subset L^{\infty}(X)\rtimes\Gamma)$ and $(L^{\infty}(Y)\subset L^{\infty}(Y)\rtimes\Gamma)$ are isomorphic ([Si55], [FM77]). This shows that orbit equivalent actions are W$^*$--equivalent; the converse is false by [CJ82].

A p.m.p. action $\Gamma\curvearrowright (X,\mu)$ is {\it weakly mixing} if for any measurable sets $A_1,..,A_n$ of $X$ and every $\varepsilon>0$ we can find $\gamma\in\Gamma$ such that $|\mu(\gamma A_i\cap A_j)-\mu(A_i)\mu(A_j)|<\varepsilon$, for all $i.j$. It is
 {\it mixing} if for any measurable sets $A_1,A_2\subset X$ and every $\varepsilon>0$ we can find $F\subset\Gamma$ finite such that  $|\mu(\gamma A_i\cap A_j)-\mu(A_i)\mu(A_j)|<\varepsilon$, for all $\gamma\in\Gamma\setminus F$.

\vskip 0.1in 
\noindent {\bf 1.3 Popa's intertwining technique}. In [Po06a], S. Popa  introduced a very powerful technique for proving unitary conjugacy of two subalgebras of a finite von Neumann algebra. Throughout the paper, this technique will play a central role.  Here, we recall Popa's result and  establish a strengthening of a particular case of it that will be crucial in the proof of our main results.

\proclaim{1.3.1 Theorem [Po06a, Theorem 2.1. and Corollary 2.3.]}
  Let $(M,\tau)$ be a finite von Neumann algebra together with two, possibly non--unital, von Neumann subalgebras $B$ and $Q$ (with units $1_B$ and $1_Q$, respectively). 
Then the following are equivalent:
\vskip 0.03in
\noindent
(1) There exist  non--zero projections $q\in Q, p\in B$, a $*$--homomorphism $\psi:qQq\rightarrow pBp$ and a non--zero partial isometry $v\in pMq$ such that $\psi(x)v=vx$, for all $x\in qQq$.

\noindent
(2) There exist $a_1,..,a_n\in 1_BM1_Q$ and $\varepsilon>0$ such that $\sum_{i,j=1}^n||E_{B}(a_iua_j^*)||_2^2\geq \varepsilon$, for all $u\in\Cal U(Q)$.

\noindent
(3) There exist $a_1,..,a_n\in 1_BM1_Q$, $\varepsilon>0$ and a group $\Cal U\subset\Cal U(Q)$ such that $\Cal U''=Q$ and $\sum_{i,j=1}^n||E_{B}(a_iua_j^*)||_2^2\geq \varepsilon$, for all $u\in\Cal U$.
\endproclaim
If one of these conditions holds true,  we say that {\it a corner of $Q$ embeds into $B$ inside $M$} and write $Q\prec_{M}B$. Note that if $B_1,B_2,..$ is a sequence of von Neumann subalgebras of $M$ such that $Q\nprec_{M}B_i$, for all $i\geq 1$, then for every group $\Cal U\subset \Cal U(Q)$ with $\Cal U''=Q$ we can find a sequence $\{u_n\}_{n\geq 1}\subset\Cal U$ such that $\lim_{n\rightarrow\infty}||E_{B_i}(au_nb)||_2=0$, for all $i\geq 1$ and every $a,b\in M$. This statement follows from Theorem 1.3.1 (see [IPP08, proof of Theorem 4.3] or [Va08, Remark 3.3]).
\vskip 0.1in

If $M=B\rtimes\Gamma$, for some action of a countable group $\Gamma$, then condition (2) is equivalent to the following: there exist $F\subset\Gamma$ finite and $\varepsilon>0$ such that for every $u\in U(Q)$, the Fourier coefficients $b_{\gamma}=E_B(uu_{\gamma}^*)$ satisfy $\max_{\gamma\in F}||b_{\gamma}||_2\geq\varepsilon$.   Below, we show that this is equivalent with the apparently weaker condition $\max_{\gamma\in\Gamma}||b_{\gamma}||_2\geq\varepsilon$, for all $u\in\Cal U(Q)$.

\proclaim {1.3.2 Theorem} Let  $\sigma:\Gamma\rightarrow$ Aut$(B)$ be an action of a countable group $\Gamma$ on a finite von Neumann algebra $(B,\tau)$ and set $M=B\rtimes_{\sigma}\Gamma$. Let $Q\subset M$ be a von Neumann subalgebra. If there exists $\varepsilon>0$ such that $\max_{\gamma\in\Gamma}||E_{B}(uu_{\gamma}^*)||_2\geq\varepsilon,$ for all $u\in\Cal U(Q),$ then $Q\prec_{B}M$.\endproclaim

\noindent
{\it Proof.} We begin the proof with the following:
\vskip 0.05in
\noindent {\it Claim 1.}
Let $\alpha$ be the flip automorphism of $Q\overline{\otimes}Q$, i.e. $\alpha(x\otimes y)=y\otimes x$, for all $x,y\in Q$, and let $\Cal Q=(Q\overline{\otimes}Q)^{\alpha}$ be the von Neumann algebra of $\alpha$--fixed points. Then $\Cal U:=\{u\otimes u|u\in\Cal U(Q)\}\subset\Cal U(\Cal Q)$  satisfies $\Cal U''=\Cal Q$.
\vskip 0.05in
\noindent{\it Proof of claim 1.}  Let  $\Cal H$ be  the $||.||_2$--closure of the span of $\Cal U$. Then the claim is equivalent to $L^2(\Cal Q)=\Cal H$.
First notice that the span of $\{x\otimes y+y\otimes x|x,y\in\ Q\}$ is dense in $L^2(\Cal Q)$. Since $x,y\in Q$ are finite linear combinations of hermitian elements of $Q$, the span of $\{h_1\otimes h_2+h_2\otimes h_1|h_1=h_1^*\in Q, h_2=h_2^*\in Q\}$ is also dense in $L^2(\Cal Q)$. Next, remark that for every $h_1,h_2\in Q$ we have that $h_1\otimes h_2+h_2\otimes h_1=(h_1+h_2)\otimes (h_1+h_2)-h_1\otimes h_1-h_2\otimes h_2$.

Altogether, to prove the claim it suffices to show that for every hermitian element $h\in Q$ we have that $h\otimes h\in\Cal H$. Towards proving this, fix $h=h^*\in Q$. Let $P$ be the orthogonal projection from $L^2(\Cal Q)$ onto $\Cal H$ and let $a:\Bbb R\rightarrow \Cal H$ be given by $a(t)=P(e^{ith}\otimes e^{ith})$.  By the definition of $\Cal H$ we get that $a(t)=0$, for all $t\in \Bbb R$. On the other hand, we have that $e^{ith}\otimes e^{ith}=\sum_{m,n\geq 0}\frac{1}{m!n!}i^{m+n}t^{m+n}(h^m\otimes h^n)$, where the sum is absolutely convergent  in the uniform norm and thus in $||.||_2$. Therefore we have that $0=a(t)=\sum_{m,n\geq 0}\frac{1}{m!n!}i^{m+n}t^{m+n}P(h^m\otimes h^n)$, for all $t\in\Bbb R$, and by analyticity all the coefficients of $a$ must be equal to 0. In particular, we derive that $h\otimes 1+1\otimes h\in\Cal H$ and $\frac{1}{2}(h^2\otimes 1)+h\otimes h+\frac{1}{2}(1\otimes h^2)
\in\Cal H$, for all hermitians $h\in Q$. Thus $h\otimes h\in\Cal H$, which concludes the proof of the claim.\hfill$\square$
\vskip 0.05in
Further, we can assume, that $Q$ is diffuse, otherwise the conclusion of the lemma is trivial. Then we can find $v\in \Cal U(Q\overline{\otimes}Q)$ such that $\alpha(v)=-v$. Indeed, $Q$ contains a copy of $L^{\infty}([0,1])$, so we can just take $v\in L^{\infty}([0,1]\times [0,1])\subset Q\overline{\otimes}Q$ given by $v(x,y)=1$, if $x\geq y$ and  $v(x,y)=-1$, if $x<y$.  If we let $\Cal V=\Cal U(\Cal Q)\cup\Cal U(\Cal Q)v$, then $\Cal V$ is a subgroup of $\Cal U(Q\overline{\otimes}Q)$ and $\Cal V''=Q\overline{\otimes}Q$. 

Now, write $M\overline{\otimes}M=(B\overline{\otimes} B)\rtimes(\Gamma\times\Gamma)$, let $\Delta(\Gamma)=\{(\gamma,\gamma)|\gamma\in\Gamma\}\subset\Gamma\times\Gamma$ and set $N=(B\overline{\otimes}B)\rtimes\Delta(\Gamma)$.
\vskip 0.05in
\noindent{\it Claim 2.} $Q\overline{\otimes}Q\prec_{M\overline{\otimes}M}N$.
\vskip 0.05in
\noindent{\it Proof of claim 2.} We start by noticing that for every $u\in\Cal U(Q)$ $$||E_{N}(u\otimes u)||_2^2=||E_{N}(\sum_{(\gamma_1,\gamma_2)\in\Gamma\times\Gamma}(E_{B}(uu_{\gamma_1}^*)\otimes E_{B}(uu_{\gamma_2}^*))u_{(\gamma_1,\gamma_2)})||_2^2=$$ $$||\sum_{\gamma\in\Gamma}(E_{B}(uu_{\gamma}^*)\otimes E_{B}(uu_{\gamma}^*))u_{(\gamma,\gamma)})||_2^2=\sum_{\gamma\in\Gamma}||E_B(uu_{\gamma}^*)||_2^4\geq\varepsilon^4$$

By combining claim 1 with Popa's result (Theorem  1.3.1) we get that $\Cal Q\prec_{M\overline{\otimes}M}N$. Since $\Cal V=\Cal U(\Cal Q)\cup\Cal U(\Cal Q)v$ is a group and satisfies $\Cal V''=Q\overline{\otimes}Q$, it easily follows that $Q\overline{\otimes}Q\prec_{M\overline{\otimes}M}N$, as claimed.\hfill$\square$
\vskip 0.05in

We are now ready to show that $Q\prec_{M}B$. By claim 2, we can find $a_1,..,a_n\in M\overline{\otimes}M$ and $c>0$ such that $\sum_{i,j=1}^n||E_N(a_iua_j^*)||_2^2\geq c$, for all $u\in\Cal U(Q\overline{\otimes}Q)$.  By using $||.||_2$ approximations (see 1.2) and the fact that $E_N$ is $N$--bimodular, we may assume that $a_i=u_{(e,\gamma_i)}$, for some $\gamma_i\in\Gamma$.  Then, for every $u\in\Cal U(Q)$, $u\otimes 1\in\Cal U(Q\overline{\otimes}Q)$,  hence  $$c\leq\sum_{i,j}^n||E_{N}(a_i(u\otimes 1)a_j^*||_2^2=\sum_{i,j=1}^{n}||E_{N}(u\otimes u_{\gamma_i\gamma_j^{-1}})||_2^2=\sum_{i,j=1}^{n}||E_B(uu_{\gamma_i\gamma_j^{-1}}^*)||_2^2,$$
which proves the lemma. \hfill$\blacksquare$

\vskip 0.05in
We end this subsection with the following lemma, whose proof we leave as exercise.

\proclaim {1.3.3 Lemma} Let $(M,\tau)$ be a finite von Neumann algebra and let $Q,B\subset M$ be two von Neumann subalgebras. Let $p\in \Cal Z(Q)$ be a central projection and  set $p'=\bigvee_{u\in\Cal N_{M}(Q)}upu^*\in\Cal Z(Q)$. If $Qp'\prec_{M}B$ then $Qp\prec_{M}B$.

\endproclaim
%
%\noindent
%{\it Proof.}
%If  $Qp'\prec_{M}B$, then we can find non--zero projections $q\in Qp'$, $r\in B$, a $*$--homomorphism $\psi:q(Qp')q\rightarrow rBr$ and a partial isometry $0\not=v\in rMq$ such that $\psi(x)v=vx$, for all $x\in q(Qp')q$.  Let $q_0\in Q$ be a projection such that $q=q_0p'$.
%Since $0\not=vq=v(q_0p')$ we can find $u\in\Cal N_{M}(Q)$ such that $v'=v(q_0(upu^*))\not=0$. Now, we have that $upu^*, p'\in \Cal Z(Q)$ and thus $q_0Qq_0(upu^*)\subset q_0Qq_0p'=q(Qp')q$. The restriction $\phi$ of $\psi$ to $q_0Qq_0(upu^*)$ satisfies $\phi(x)v'=v'x$, for all $x\in q_0Qq_0(upu^*)$. Since $0\not= v'\in rM(q_0(upu^*))$, we get that a corner of $Q(upu^*)=u(Qp)u^*$ embeds into $B$, proving the conclusion.\hfill$\blacksquare$ 

\vskip 0.1in
\noindent {\bf 1.4 Rigid inclusions of von Neumann algebras}. We next recall S. Popa's notion of rigidity for inclusions of finite von Neumann algebras. Let $(M,\tau)$ be a finite von Neumann algebra together with a von Neumann subalgebra $B$. 
The inclusion $(B\subset M)$ is {\it rigid} (or, has the {\it relative property (T)}) if  any sequence of subunital ($\Phi_n(1)\leq 1$), subtracial ($\Phi_n\circ\tau\leq\tau$),  normal completely positive maps $\Phi_n:M\rightarrow M$ which satisfy $\lim_{n\rightarrow\infty}||\Phi_n(x)-x||_2=0$, for all $x\in M$, converges to the identity uniformly on the unit ball of $B$, i.e. $\lim_{n\rightarrow\infty}\sup_{x\in B,\hskip 0.02in ||x||\leq 1}||\Phi_n(x)-x||_2=0$ ([Po06c]). In the case when $B=M$, we say that $M$ has property (T) ([CJ85]).

For two countable groups $\Gamma_0\subset\Gamma$, the inclusion $(L(\Gamma_0)\subset\ L(\Gamma))$ is rigid if and only if the inclusion $(\Gamma_0\subset\Gamma)$ has the {\it relative property (T)} ([Po06c, Proposition 5.1.]). By the classical results of Kazhdan and Margulis, the inclusions (SL$_n(\Bbb Z)\subset$ SL$_n(\Bbb Z))$ ($n\geq 3)$ and $(\Bbb Z^2\subset\Bbb Z^2\rtimes$ SL$_2(\Bbb Z))$  have the relative property (T) ([Ka67],[Ma82]). Recall that a group $\Gamma$ has property (T) of Kazhdan if and only if the inclusion $(\Gamma\subset\Gamma)$ has the relative property (T). For examples of property (T) groups, see the extensive monograph [BHV08].

\vskip 0.1in
\noindent {\bf 1.5 Weakly malleable deformations of  Bernoulli actions.} S. Popa discovered that Bernoulli actions $\Gamma\curvearrowright (X_0,\mu_0)^{\Gamma}$  have a remarkable deformation property, and called it {\it malleability} ([Po06a]). By pairing it with  property (T) of the group $\Gamma$, he proved striking rigidity results concerning the associated II$_1$ factor ([Po06ab]).
Since then, malleable deformations have been found in several other contexts and are now a central tool in Popa's deformation/rigidity theory (see [Po07b, Section 6]). In [Io07], in order to extend some of Popa's results [Po06ab] to II$_1$ factors coming from Bernoulli actions with arbitrary base, we introduced a new class of malleable deformations.
\vskip 0.05in
To recall their construction from [Io07, Section 2], let $(B,\tau)$ be a finite von Neumann algebra 
and let $\Gamma\curvearrowright I$ be an action of a countable group on a countable set. For every set $J\subset I$, we denote ${B}^{J}=\overline{\bigotimes}_{i\in J}(B)_i$. Whenever $J\subset J'$,  we view $B^{J}$ as a subalgebra of $B^{J'}$, in the natural way. 
Let $\sigma:\Gamma\rightarrow$ Aut$(B^I)$ be the {\it generalized Bernoulli action} defined by $\sigma(\gamma)(\otimes_{i\in I}x_i)=\otimes_{i\in I}x_{\gamma^{-1}\cdot i},$ for every $x=\otimes_{i\in I}x_i\in B^I$ and  $\gamma\in \Gamma$.
Let $M=B^{I}\rtimes_{\sigma}\Gamma$ be the associated crossed product von Neumann algebra.
 
Next,  we augment $M$ to a von Neumann algebra $\tilde M$ and define a $1$--parameter group of automorphisms $\{\theta_t\}_{t\in\Bbb R}$ of $\tilde M$ such that $\theta_t\rightarrow$ id$_{\tilde M}$,  in the  pointwise $||.||_2$ topology, as $t\rightarrow 0$. Towards this, we define  the free product von Neumann algebra $\tilde B=B*L(\Bbb Z)$.
 Let $\tilde\sigma:\Gamma\rightarrow \text{Aut}({\tilde B}^I)$ be the generalized Bernoulli action. It is clear that  $B^I\subset {\tilde B}^I$ and that $\tilde\sigma$ extends $\sigma$, hence we have the inclusion $M\subset \tilde M:={\tilde B}^I\rtimes_{\tilde\sigma}\Gamma.$

Now, let $u\in L(\Bbb Z)$ be a Haar unitary such that $L(\Bbb Z)=\{u^n|n\in\Bbb Z\}''$ and let $h\in L(\Bbb Z)$ be a hermitian element such that $u=e^{ih}$. For every $t\in\Bbb R$, we define the unitary element  $u_t=e^{ith}\in L(\Bbb Z)$ and consider the automorphism $\theta_t=\otimes_{i\in I}\text{Ad}(u_t)_i$ of ${\tilde B}^I.$
Since  $\theta_t$ commutes with $\tilde\sigma$,  it extends to an automorphism of $\tilde M$ through the formula  $\theta_t(x)=\sum_{\gamma\in\Gamma}\theta_t(x_{\gamma})u_{\gamma},$ for all $x=\sum_{\gamma\in\Gamma}x_{\gamma}u_{\gamma}\in \tilde M$. Since  $\lim_{t\rightarrow 0}||u_t-1||_2=0$, we get that $\theta_t\rightarrow$ id$_{\tilde M}$, , as $t\rightarrow 0$.
We end this section by noticing that $\{\theta_t\}_{t\in\Bbb R}$ admits a certain $\beta-$symmetry (see e.g. [Po06a, section 1.4]).

\proclaim {1.5.1 Lemma} There exists an automorphism $\beta$ of $\tilde M$ such that $\beta^2=\text{id}_{\tilde M}$, $\beta_{|M}=\text{id}_{M}$ and $\beta\theta_t\beta=\theta_{-t}$, for all $t\in\Bbb R$. 
\endproclaim
\noindent
{\it Proof.} Since $u$ generates $L(\Bbb Z)$, it follows that $M$ and $\{\otimes_{i\in I}u^{n_i}|n_{i}\in\Bbb Z, |\{i\in I|n_i\not=0\}|<\infty\}$ generate $\tilde M$ as a von Neumann algebra.  Using this observation it is immediate to see that $\beta$  defined by $\beta_{|M}=\text{id}_{M}$ and $\beta(\otimes_{i\in I}u^{n_i})=\otimes_{i\in I}u^{-n_i}$, for every $n_i\in\Bbb Z$ such that $\{i\in I|n_i\not=0\}$ is finite, extends to an automorphism of $\tilde M$ with the desired properties.
\hfill$\blacksquare$

\vskip 0.2in
\head \S 2. {Rigid subalgebras of II$_1$ factors arising from\\ generalized Bernoulli actions}\endhead

\vskip 0.2in
Let $M=B^{\Gamma}\rtimes\Gamma$ be the II$_1$ factor associated with a Bernoulli action. S. Popa showed, as part of his {\it deformation/rigidity theory}, that if $B$ is abelian then any  rigid subalgebra $Q$ of $M$ whose normalizer generates a factor can be conjugated inside $L(\Gamma)$ by a unitary element ([Po06a, Theorem 4.1.]). The proof exploits the tension between the rigidity of the inclusion $Q\subset M$ and the malleability property of Bernoulli actions.

In [Io07, Theorems 3.6., 3.7.], by using the weakly malleable deformations associated with $M$  we showed that any rigid subalgebra $Q$ of $M$ has a corner which embeds into $L(\Gamma)$, regardless of any regularity property of $Q$. In this section, by relying on ideas and techniques from [Io07] and by following a deformation/rigidity strategy we extend the last result to the class of generalized Bernoulli actions.  

\proclaim {2.1 Theorem} Let $\Gamma$ a countable group acting on a countable set $I$. For $i\in I$, denote by $\Gamma_i$ its stabilizer in $\Gamma$. Let $(B,\tau)$ be an abelian von Neumann algebra. Let $\sigma:\Gamma\rightarrow$ Aut$(B^I)$ be the generalized Bernoulli action  and denote $M=B^I\rtimes_{\sigma}\Gamma$. Suppose that $Q$ is von Neumann subalgebra of $M$ such that the inclusion $(Q\subset M)$ is rigid. 
\vskip 0.03in
\noindent
Then $Q\prec_{M} L(\Gamma).$ 
\vskip 0.03in
\noindent
Moreover, if $\Gamma$ is ICC and $Q\nprec_{M}L(\Gamma_i)$, for every $i\in I$, then there exists a unitary $u\in M$ such that $uQu^*\subset L(\Gamma)$. Furthermore, in this case, we have that $uPu^*\subset L(\Gamma)$, where $P$ is the von Neumann algebra generated by the quasi--normalizer of $Q$ in $M$.
\endproclaim		
\vskip 0.05in
This theorem generalizes three other results in the literature. The first part and the moreover part of 2.1 extend respectively [Va08, Lemma 5.2.] and [PV08, Theorem 6.5.] in whose hypotheses it is additionally assumed that $P\nprec_{M}B^{I}\rtimes_{\sigma}\Gamma_i$, for all $i\in I$.  Theorem 2.1 also generalizes (a particular case of) [Io09, Proposition 3.3] which shows that $M$ admits no rigid subalgebra $Q$ such that $Q\subset B^I$. 

\vskip 0.1in

\noindent
{\it Proof.}  Let $\tilde M$ and $\{\theta_t\}_{t\in\Bbb R}$ be defined as in Section 1.5.
 Since the inclusion $Q\subset M$ is rigid, the inclusion $Q\subset\tilde M$ is rigid (by  [Po06c, Proposition 4.6.]). Since $\theta_t\rightarrow$ id$_{\tilde M}$, we can find $t>0$ such that $||\theta_t(u)-u||_2\leq \frac{1}{2}$, for every $ u\in \Cal U(Q).$ 
Let $v$ be unique element of minimal $2-$norm in the $||.||_2$--closed convex hull of the set $\{\theta_t(u)u^*|u\in\Cal U(Q)\}$. Then $||v-1||_2\leq \frac{1}{2}$, hence $v\not=0$. Also, by construction, $v\in \tilde M$, $||v||\leq 1$ and $\theta_t(u)v=vu$, for all $u\in\Cal U(Q)$.
\vskip 0.05in

Next, let $\varepsilon>0$ such that $\delta:=\|v\|_2-(1+\frac{2}{\sqrt{1-|\tau(u_t)|^4}})\varepsilon>0$. Then we can find finite sets $F\subset I$, $K\subset \Gamma$ and $w$ in the linear span of $\{\tilde B^F u_{\gamma}|\gamma\in K\}$ such that $||w-v||_2\leq\varepsilon$ and $||w||\leq 1$. Suppose for simplicity that $K=\{\gamma^{-1}|\gamma\in K\}$. Also, we  use the notation $KF=\{k\cdot i|k\in K,i\in F\}$. Denote by $T$  be the orthogonal projection from $L^2(\tilde M)$ onto the closed linear span of $\{(\tilde B^{F}\otimes B^{\gamma KF\setminus F})u_{\gamma}|\gamma\in\Gamma\}$.

\vskip 0.1in
\noindent
{\it Claim 1.}  For all $u\in\Cal U(Q)$, we have \hskip 0.05in $||T(wu)||_2\geq\delta.$
\vskip 0.1in

\noindent
{\it Proof of claim 1.} Fix $u\in\Cal U(Q)$. Then  $||\theta_t(u)w-wu||_2=||\theta_t(u)(w-v)-(w-v)u||_2\leq 2\varepsilon$. Decompose $wu=\sum_{\gamma\in\Gamma}a_{\gamma}u_{\gamma}$ and $\theta_t(u)w=\sum_{\gamma\in\Gamma}b_{\gamma}u_{\gamma}$.
Since $u\in M$,  $w\in\sum_{\gamma\in K}\tilde B^F u_{\gamma}$ and $K=K^{-1}$, an easy calculation shows that  $a_{\gamma}\in \tilde B^{F}\overline{\otimes}B^{I\setminus F}$ and $b_{\gamma}\in\tilde B^{\gamma KF}\overline{\otimes}(u_tBu_t^*)^{I\setminus \gamma KF}$, for all $\gamma\in\Gamma$. 

Denote by $S_{\gamma}$ and $R_{\gamma}$ the orthogonal projections from $L^2(\tilde B^I)$ onto  $L^2(\tilde B^{F\cup\gamma KF})$ and onto $L^2(\tilde B^{\gamma KF}\overline{\otimes}(u_tBu_t^*)^{I\setminus \gamma KF})$, respectively. Since $a_{\gamma}\in \tilde B^{F}\overline{\otimes}B^{I\setminus F}$, [Io09, Lemma 3.5] implies that 
$$||R_{\gamma}(a_{\gamma})||_2^2\leq |\tau(u_t)|^4||a_{\gamma}||_2^2 + (1-|\tau(u_t)|^4)||S_{\gamma}(a_{\gamma})||_2^2,\forall\gamma\in\Gamma\tag 2.a$$

Now, since $R_{\gamma}(b_{\gamma})=b_{\gamma}$, for every $\gamma\in\Gamma$, we get that $$\sum_{\gamma\in\Gamma} {||a_{\gamma}-R_{\gamma}(a_{\gamma})||}_2^2\leq \sum_{\gamma\in\Gamma}||a_{\gamma}-b_{\gamma}||_2^2= {||wu-\theta_t(u)w||}_2^2\leq 4\varepsilon^2\tag 2.b$$
On the other hand, (2.a) implies that $$ \sum_{\gamma\in\Gamma}||a_{\gamma}-R_{\gamma}(a_{\gamma})||_2^2\geq (1-|\tau(u_t)|^4)\sum_{\gamma\in\Gamma}||a_{\gamma}-S_{\gamma}(a_{\gamma})||_2^2\tag 2.c$$  Since $a_{\gamma}\in \tilde B^{F}\overline{\otimes}B^{I\setminus F}$ we have $T(wu)=\sum_{\gamma\in\Gamma}S_{\gamma}(a_{\gamma})u_{\gamma}$. By combining (2.b) and (2.c) we get that $\|T(wu)-wu\|_2^2\leq 4\varepsilon^2(1-|\tau(u_t)|^4)^{-1}$, which yields the claim.\hfill$\square$

\vskip 0.1in

Next, let $\{B_n\}_{n\geq 1}$ be an increasing sequence of unital, finite dimensional subalgebras of $B$ such that $\cup_{n\geq 1}B_n$ is dense in $B$. For $n\geq 1$, set  $M_n=B_n^I\rtimes\Gamma$ and denote by $E_n$ the conditional expectation from $M$ onto $M_n$. Since the inclusion $Q\subset M$ is rigid and $E_n\rightarrow\text{id}_{M}$, we can find $n_0\geq 1$ such that $||E_{n_0}(u)-u||_2\leq\frac{\delta}{2}$, for all $u\in\Cal U(Q)$. Together with (2.c) this yields $$||T(wE_{n_0}(u))||_2\geq\frac{\delta}{2},\forall u\in\Cal U(Q)\tag 2.d$$

\vskip 0.05in
\noindent
{\it Claim 2.} We have $Q\prec_{M}L(\Gamma)$.
\vskip 0.05in
\noindent
{\it Proof of claim 2.} Assuming this is false, we get a sequence $\{u_m\}_{m\geq 1}\subset\Cal U(Q)$ such that $||E_{L(\Gamma)}(au_mb)||_2\rightarrow 0,$ for all $a,b\in M.$ We will prove that $$||T(cu_gE_{n_0}(u_m))||_2\rightarrow 0,\forall c\in\tilde B^F, g\in \Gamma\tag 2.e$$
Note that since $w$ belongs to the linear span of $\{\tilde B^F u_{\gamma}|\gamma\in K\}$, (2.e) implies that $||T(wE_{n_0}(u_m))||_2\rightarrow 0$, which contradicts (2.d).

To deduce (2.e), write $u_m=\sum_{\gamma\in\Gamma}x_m^{\gamma}u_{\gamma}$, where $x_m^{\gamma}\in B^{\Gamma}$, and assume $||c||\leq 1$. Since ${E_{n_0}}_{|B^{\Gamma}}=E_{B_{n_0}^{\Gamma}}$ and $\sigma$ commute we have  $cu_gE_{n_0}(u_m)=\sum_{\gamma\in\Gamma}cE_{B_{n_0}^{\Gamma}}(\sigma_g(x_m^{\gamma}))u_{g\gamma}$.
 As $c\in \tilde B^F$, we get that $T(cu_gE_{n_0}(u_m))=\sum_{\gamma\in\Gamma}E_{\tilde B^F\overline{\otimes}B^{g\gamma KF\setminus F}}(cE_{B_{n_0}^{\Gamma}}(\sigma_g(x_m^{\gamma}))u_{g\gamma}=c\sum_{\gamma\in\Gamma}E_{B_{n_0}^{F\cup g\gamma KF}}(\sigma_g(x_m^{\gamma}))u_{g\gamma}.$
Thus, as $||c||\leq 1$, we have $$||T(cu_gE(u_m))||_2^2\leq \sum_{\gamma\in\Gamma}||E_{B_{n_0}^{F\cup g\gamma KF}}(\sigma_g(x_m^{\gamma}))||_2^2=\sum_{\gamma
\in\Gamma}||E_{B_{n_0}^{g^{-1}F\cup \gamma KF}}(x_m^{\gamma})||_2^2\tag 2.f$$
\noindent
Let $\{\xi_i\}_{i=1}^l$ be an orthonormal basis for $B_{n_0}$. For  $J\subset I$ finite and $s=(s_i)_{i\in J}\in\{1,..,l\}^J$, set $\xi_s=\otimes_{i\in J}\xi_{s_i}\in B_{n_0}^J$. Then $\{\xi_s\}_{s\in \{1,..,l\}^J}$ is an orthonormal basis for $B_{n_0}^J$, hence $$||E_{B_{n_0}^{g^{-1}F\cup \gamma KF}}(x_m^{\gamma})||_2^2=\sum_{s\in \{1,..,l\}^{g^{-1}F}}\sum_{t\in \{1,..,l\}^{\gamma KF\setminus g^{-1}F}}|\tau(\xi_s x_m^{\gamma}\xi_t)|^2\leq\tag 2.g$$ $$\sum_{s\in \{1,..,l\}^{g^{-1}F}}\sum_{t\in \{1,..,l\}^{KF}}|\tau(\xi_sx_m^{\gamma}\sigma_{\gamma}(\xi_t))|^2$$
The combination of (2.f) and (2.g) gives that $$||T(cu_gE(u_m))||_2^2\leq\sum_{s\in \{1,..,l\}^{g^{-1}F}}\sum_{t\in \{1,..,l\}^{KF}} ||E_{L(\Gamma)}(\xi_su_m\xi_t)||_2^2\rightarrow 0.$$

\vskip 0.05in
\noindent
For the moreover part of the statement, assume that $\Gamma$ is ICC, i.e. $L(\Gamma)$ is a factor, and that $Q\nprec_{M}L(\Gamma_i)$, for all $i\in I$.

\vskip 0.1in

\noindent
{\it Claim 3.} For every relatively rigid von Neumann subalgebra $Q\subset M$ which verifies $Q\nprec_{M}L(\Gamma_i)$, for all $i\in I$, we can find a unitary $u\in M$ and a non--zero projection $q'\in Q'\cap M$ such that $u(Qq')u^*\subset L(\Gamma)$.

\vskip 0.05in
\noindent
{\it Proof of claim 3.} By the first part of the proof we have that $Q\prec_M L(\Gamma)$. Since $Q\nprec_{M}L(\Gamma_i)$, by [Va08, Remark 3.8.] we can find non--zero projections $q\in Q$, $p\in L(\Gamma)$, a $*$--homomorphism $\psi:qQq\rightarrow pL(\Gamma)p$ and a partial isometry $0\not=v\in pMq$  such that $\psi(x)v=vx$, for all $x\in qQq$, and $\psi(qQq)\nprec_{L(\Gamma)}L(\Gamma_i)$, for all $i\in I$. Thus, by [Va08, Lemma 4.2.$(1)$] we deduce that $\psi(qQq)'\cap pL(\Gamma)p\subset pL(\Gamma)p$. In particular, it follows that $vv^*\in pL(\Gamma)p$ and hence $vQv^*\subset L(\Gamma)$. Since $v^*v\in(qQq)'\cap qMq$ we can find a projection  $q''\in Q'\cap M$ such that $v^*v=qq''$. Let $w\in M$ be a unitary extending $v$, then $w(qQqq'')w^*=vQv^*\subset L(\Gamma)$. Since $L(\Gamma)$ is a factor it follows that we can find a unitary $u\in M$ such that $u((Qq'')z){u}^*\subset L(\Gamma)$, where $z$ is the central support of $qq''$ in $Qq''$. Since the center of $Qq''$ is contained in $Q'\cap M$ we  deduce that $0\not=q'=q''z\in\Cal P(Q'\cap M)$ satisfies the claim. 
\hfill$\square$
\vskip 0.05in

Now, let $q\in Q,q'\in Q'\cap M$ such that $\tau(qq')=\frac{1}{n}$, where $n\geq 1$ is an integer, and consider the natural embedding of  $Q_0:=M_{n\times n}(qQqq')$ into $M$. By [Po06c], the inclusion $Q_0\subset M$ is rigid. Since $Q\nprec_{M}L(\Gamma_i)$, we get that $Q_0\nprec_{M} L(\Gamma_i)$, for all $i\in I$. Thus the conclusion of claim 3 holds true for $Q_0$.
A close inspection of the last part of the proof of [Po06a, Theorem 4.4. {\it (ii)}] reveals that this fact and the factoriality of $L(\Gamma)$ are enough to imply (via a maximality argument) that we can find $u\in\Cal U(M)$ such that $uQu^*\subset L(\Gamma)$. Finally, [Va08, Lemma 4.2.$(1)$] implies  that $uPu^*\subset L(\Gamma)$.\hfill$\blacksquare$

\vskip 0.1in
\noindent
{\bf 2.2 Remark.} Let $N$ be a finite von Neumann algebra. The above proof  can be modified  to show that if a Neumann subalgebra $Q\subset M\overline{\otimes}N$ satisfies $(\theta_t\otimes\text{id}_N)(x)v=vx$, for every $x\in Q$, for some fixed  $0\not=v\in\tilde M\overline{\otimes} N$ and  $t\not=0$, then $Q\prec_{M\overline{\otimes}N}L(\Gamma)\overline{\otimes}N$.  To see this, just  replace throughout the proof  any subalgebra $A\subset M$ with the subalgebra $A\overline{\otimes}N$ of $M\overline{\otimes}N$.

\vskip 0.2in

\head \S 3.{ A spectral gap argument}.\endhead

\vskip 0.2in
In the next sections, we prove a series of results concerning {\it tensor products} of II$_1$ factors arising from Bernoulli actions (Theorems 3.2, 4.1 and 6.1). These results have natural analogs for a {\it single} II$_1$ factor $M=B^{\Gamma}\rtimes\Gamma$ coming from a Bernoulli action with $B$ abelian. For simplicity, we will only describe the latter, although below we state and prove the former. Denote $A=B^{\Gamma}$.

In this section, we prove that if $D\subset M$ is a von Neumann subalgebra such that there exists a sequence $x=(x_n)_n\in D'\cap M^{\omega}$ on which the deformations $\theta_t$ converge uniformly  but  $x\notin A^{\omega}\rtimes\Gamma$, then a corner of $D$ embeds into $L(\Gamma)$.

This statement was inspired by an idea of J. Peterson. In the context when $\Gamma$ is a free product group and the action of $\Gamma$ on $A$ is compact, he analyzed sequences $\{x_n\}_n\subset D'\cap M$ on which certain deformations converge uniformly in order to derive conjugacy results for $D$ inside $M$  (see [Pe10, Theorem 4.1]).
The proof of our result, rather than using arguments from [Pe10], relies on Popa's spectral gap argument ([Po08]).

\vskip 0.05in

Before continuing, we fix some notations that we will use  throughout the paper.
\vskip 0.05in
\noindent
{\bf 3.1 Notations:} 

\noindent
$\bullet$ Let $\Gamma_1,\Gamma_2$ be countable groups and $B_1,B_2$ be abelian von Neumann algebras. 

\noindent
$\bullet$
Consider the Bernoulli actions $\sigma_i:\Gamma_i\rightarrow$ Aut$(B_i^{\Gamma_i})$ and denote $M_i=B_i^{\Gamma}\rtimes_{\sigma_i}\Gamma_i$.

\noindent
$\bullet$ Denote $M=M_1\overline{\otimes}M_2$, $A=B_1^{\Gamma_1}\overline{\otimes}B_2^{\Gamma_2}$ and $\Gamma=\Gamma_1\times\Gamma_2$.

\noindent
$\bullet$ Note that $M=A\rtimes_{\sigma}\Gamma$, where $\sigma(\gamma_1,\gamma_2)=\sigma_1(\gamma_1)\otimes\sigma_2(\gamma_2)$.

\proclaim {3.2 Theorem} Let $\tilde M_i=(B_i*L(\Bbb Z))^{\Gamma}\rtimes\Gamma\supset M_i$ and $\{\theta_t^i\}_{t\in\Bbb R}\subset$ \text{Aut}$(\tilde M_i)$ be defined as in section 1.5. Let $\tilde M=\tilde M_1\overline{\otimes}\tilde M_2$ and for every $t\in\Bbb R$, set $\theta_t=\theta_t^1\otimes\theta_t^2\in$ \text{Aut}$(\tilde M)$. 
\vskip 0.03in
\noindent
Let $D\subset M$ be a von Neumann subalgebra and  $x\in D'\cap M^{\omega}$ satisfying $\lim_{t\rightarrow 0}||\theta_t(x)-x||_2=0$. Then one of the following holds true:
\vskip 0.03in
\noindent
(1) $x\in A^{\omega}\rtimes\Gamma$.

\noindent
(2) $D\prec_{M}M_1\overline{\otimes} L(\Gamma_2)$ or $D\prec_{M}L(\Gamma_1)\overline{\otimes} M_2.$
\endproclaim
\noindent
Here, we view $\theta_t$ as an automorphism of $\tilde M^{\omega}$ given by $\theta_t((x_n)_n)=(\theta_t(x_n))_n$. Also, we denote by $A^{\omega}\rtimes\Gamma$ the von Neumann algebra generated by $A^{\omega}$ and $\{u_{\gamma}\}_{\gamma\in\Gamma}$ inside $M^{\omega}$.

\vskip 0.05in
The proof of Theorem 3.2 is based on the following general technical result.

\proclaim {3.3 Lemma} Let $\sigma:\Gamma\rightarrow$ Aut$(\tilde A)$ be an action on a finite von Neumann algebra $(\tilde A,\tau)$ which leaves invariant a von Neumann subalgebra $A\subset\tilde A$. Denote $M=A\rtimes\Gamma$, $\tilde M=\tilde A\rtimes\Gamma$. Assume that the Hilbert $M$--bimodule $L^2(\tilde M)\ominus L^2(M)$ is isomorphic to $\oplus_{i\in I}L^2(\langle M,e_{A_i}\rangle)$, where $A_i\subset M$ are von Neumann subalgebras such that for every $i\in I$, we have $A_i\subset A\rtimes\Gamma_i$, for some finite subgroup $\Gamma_i\subset\Gamma$.
 Let $\{\theta_t\}_{t\in\Bbb R}$ be a 1-parameter group of automorphisms of $\tilde M$ and $\beta\in\text{Aut}(\tilde M)$  such that $\beta^2=1_{\tilde M}$, $\beta_{|M}=1_{M}$, $\beta\theta_t\beta=\theta_{-t}$, for all $t\in\Bbb R$.

\vskip 0.03in
\noindent
Let $D\subset M$ be a von Neumann subalgebra. Assume that there exist no $t\not=0$ and $0\not= v\in\tilde M$ such that $\theta_t(y)v=vy$, for all $y\in D$. Then  we have  $$||x-E_{A^{\omega}\rtimes\Gamma}(x)||_2\leq 4\sqrt{2} \liminf_{t\rightarrow 0}||\theta_t(x)-x||_2,\hskip 0.05in \forall x\in D'\cap M^{\omega}.$$ 
\endproclaim

\noindent
{\it Proof of Lemma 3.3.} Let  $x=(x_n)_n\in D'\cap M^{\omega}$. After rescaling we may assume that $||x||\leq 1$.  
Let $y=x-E_{A^{\omega}\rtimes\Gamma}(x)\in M^{\omega}$. Since $||y||\leq 2$, we can represent $y=(y_n)_n$, where $y_n\in M$ satisfy $||y_n||\leq 2$, for all $n$. Assuming the conclusion if false, we can find $0<\varepsilon<\frac{||y||_2}{4\sqrt{2}}$ and $t\not=0$ such that $||\theta_t(x)-x||_2=\lim_{n\rightarrow\omega}||\theta_t(x_n)-x_n||_2\leq\varepsilon$.
\vskip 0.05in
 For $u\in M$, we denote $\delta_t(u):=\theta_t(u)-E_M(\theta_t(u))\in L^2(\tilde M)\ominus L^2(M)$.
\vskip 0.1in
\noindent
{\it Claim 1.} For every $u\in\Cal U(D)$ we have that $||[\delta_t(u),y]||_2\leq 8 \varepsilon$.
\vskip 0.05in

\noindent{\it Proof of claim 1.} Let $u\in \Cal U(D)$. Since $[u,x]=0$, by using an idea from [Po08, proof of Theorem 4.1.], we have that  $$||[\theta_t(u),x]||_2=||[u,\theta_{-t}(x)||_2=||[u,(\theta_{-t}(x)-x)]||_2\leq\tag 3.a$$ $$2||\theta_{-t}(x)-x||_2=2||x-\theta_t(x)||_2\leq 2\varepsilon.$$

Next, notice that $E_{{\tilde A}^{\omega}\rtimes\Gamma}(x)=(E_{{\tilde A}^{\omega}\rtimes\Gamma}\circ E_{M^{\omega}})(x)=E_{A^{\omega}\rtimes\Gamma}(x)$. By combining this and the fact that $\theta_t(u)\in\tilde M\subset {\tilde A}^{\omega}\rtimes\Gamma$ with (3.a) we get that $$||[\theta_t(u),y]||_2=||[\theta_t(u),(x-E_{A^{\omega}\rtimes\Gamma}(x))]||_2\leq\tag 3.b$$ $$ ||[\theta_t(u),x]||_2+||\theta_t(u)E_{A^{\omega}\rtimes\Gamma}(x)-E_{A^{\omega}\rtimes\Gamma}(x)\theta_t(u)||_2\leq$$ $$2\varepsilon+||E_{{\tilde A}^{\omega}\rtimes\Gamma}(\theta_t(u)x-x\theta_t(u))||_2\leq 2\varepsilon+||\theta_t(u)x-x\theta_t(u)||_2\leq 4\varepsilon.$$

Finally, by using (3.b) we have that $$||[\delta_t(u),y]||_2\leq 4\varepsilon+||E_M(\theta_t(u))y-yE_M(\theta_t(u))||_2=$$
 $$4\varepsilon+||E_{M^{\omega}}(\theta_t(u)y-y\theta_t(u))||_2\leq 8\varepsilon,$$ and thus the claim is proven.\hfill$\square$ 
\vskip 0.1in

\noindent
{\it Claim 2.} For every $\xi,\eta\in L^2(\tilde M)\ominus L^2(M)$, we have that $\lim_{n\rightarrow\omega}\langle y_n\xi,\eta y_n\rangle=0$.

\vskip 0.05in
\noindent{\it Proof of claim 2.} Recall that $L^2(\tilde M)\ominus L^2(M)=\oplus_{i\in I}L^2(\langle M,e_{A_i}\rangle)$. Since $||y_n||\leq 2,$ for all $n$, it is clearly sufficient to prove the claim for $\xi$ and $\eta$ of the form $\xi=\xi_1 e_{A_i}\xi_2\in L^2(\langle M,e_{A_i}\rangle)$ and $\eta=\eta_1 e_{A_i}\eta_2\in L^2(\langle M,e_{A_i}\rangle)$, for some $i\in I$, $\xi_1,\xi_2,\eta_1,\eta_2\in M$. But, in this case, if  $Tr$ denotes the natural semi--finite trace on $\langle M,e_{A_i}\rangle$, then  $$|\langle y_n\xi,\eta y_n\rangle|=|Tr(y_n^*\eta^* y_n\xi)|=|Tr(y_n^*\eta_2^*e_{A_i}\eta_1^* y_n\xi_1 e_{A_i}\xi_2)|=$$ $$|\tau(E_{A_i}(\xi_2y_n^*\eta_2^*)E_{A_i}(\eta_1^* y_n\xi_1))|\leq ||E_{A_i}(\eta_2 y_n^*\xi_2^*)||_2||E_{A_i}(\eta_1^* y_n\xi_1)||_2$$

This reduces the claim to showing that $\lim_{n\rightarrow\omega}||E_{A_i}(\eta_1^* y_n\xi_1)||_2=0$, for all $i\in I$, $\xi_1,\eta_1\in M$. To this end, let $\Gamma_i\subset\Gamma$ be a finite group such that $A_i\subset A\rtimes\Gamma_i$. Then, for every $\zeta\in M$ we have that $||E_{A_i} (\zeta)||_2^2\leq ||E_{A\rtimes\Gamma_i}(\zeta)||_2^2=\sum_{\gamma\in\Gamma_i}||E_{A}(\zeta u_{\gamma}^*)||_2^2$. Hence, we can further reduce the claim to showing that $\lim_{n\rightarrow\omega}||E_A(\eta_1^* y_n\xi_1)||_2=0$, for all $\xi_1,\eta_1\in M$. This in turn is an immediate consequence of the fact that $y=(y_n)_n\perp A^{\omega}\rtimes\Gamma$.\hfill$\square$

\vskip 0.05in

Next, if $u\in\Cal U(D)$, then by claim 2 we get that $\delta_t(u)y\perp y\delta_t(u)$. Using claim 1 we therefore derive that $$64\varepsilon^2\geq ||[\delta_t(u),y]||_2^2=|| \delta_t(u)y||_2^2+||y\delta_t(u)||_2^2.$$ In particular, we obtain that $$8\varepsilon\geq || \delta_t(u)y||_2,\forall u\in \Cal U(D)\tag 3.c$$

Now, we use the same trick as in the proof of [Po08, Lemma 2.1.]. Since $E_M(\theta_t(u))\in M$, $u\in M$, $y\in M^{\omega}$ we get that $\beta(E_M(\theta_t(u)))=E_M(\theta_t(u))$, $\beta(u)=u$ and $\beta(y)=y$. Thus $$||\theta_{-t}(u)y-E_M(\theta_t(u))y||_2=||\beta(\theta_{-t}(u)y-E_M(\theta_t(u))y)||_2=\tag 3.d$$ $$ ||\beta(\theta_{-t}(u))y-E_M(\theta_t(u))y||_2=||\theta_t(\beta(u))y-E_M(\theta_t(u))y||_2=$$ $$||\delta_t(u)y||_2\leq 8\varepsilon,\forall u\in\Cal U(D)$$

By combining (3.c) and (3.d) we deduce that $$||\theta_t(u)y-\theta_{-t}(u)y||_2\leq 16\varepsilon,\forall u\in\Cal U(D)\tag 3.e$$
Thus we have that $$\Re\tau(\theta_t(u)yy^*\theta_{-t}(u^*))=\frac{1}{2}(||\theta_t(u)y||_2^2+||\theta_{-t}(u)y||_2^2-||\theta_t(u)y-\theta_{-t}(u)y||_2^2)\geq\tag 3.f$$ $$ ||y||_2^2-128\varepsilon^2>0$$

\noindent
If $z=E_{\tilde M}(yy^*)\in\tilde M$, then  $\tau(\theta_t(u)yy^*\theta_{-t}(u^*))=\tau(\theta_t(u)z\theta_{-t}(u^*))$, hence  (3.f) implies that $\Re\tau(\theta_t(u)z\theta_{-t}(u^*))\geq ||y||_2^2-128\varepsilon>0$, for all $u\in\Cal U(D)$. A standard averaging argument provides $0\not=w\in\tilde M$ such that $\theta_t(u)w=w\theta_{-t}(u)$, for all $u\in\Cal U(D)$. Thus, if $v=\theta_t(w)\not=0$, then $\theta_{2t}(u)v=vu$, for all $u\in\Cal U(D)$, which gives a contradiction. \hfill$\blacksquare$

\vskip 0.1in
\noindent
{\it Proof of Theorem 3.2.} Assuming that (2) is false we will demonstrate that (1) holds true. 
Set ${\Cal M}=\tilde M_1\overline{\otimes}M_2$ and observe that $M\subset{\Cal M}\subset\tilde M$. Let $\theta_t':=\theta_t^1\otimes$id$_{M_2}\in$ Aut$({\Cal M})$. Let $\beta\in$ Aut$(\tilde M_1)$ be defined as in section 1.5. By Lemma 1.5.1 $\beta'=\beta\otimes$id$_{M_2}\in$ Aut(${\Cal M})$  satisfies $$\beta'_{|M}=1_{M}, {\beta'}^2=1_{{\Cal M}},\beta'\theta_t'\beta'=\theta_{-t}',\forall t\in\Bbb R\tag 3.g$$ 

\noindent
Denote $A_1=B_1^{\Gamma_1}$ and $A_2=B_2^{\Gamma_2}$. Remark that we can write $M=(A_1\overline {\otimes}M_2)\rtimes\Gamma_1$, where $\Gamma_1$ acts by Bernoulli shifts on $A_1$ and trivially on $M_2$. Similarly, we have $M=(M_1\overline{\otimes}A_2)\rtimes\Gamma_2$.

\vskip 0.1in 
\noindent
{\it Claim 1.} There exists a family $\{C_i\}_{i\in I}$ of von Neumann subalgebras of $M$ and a family $\{\Gamma_i\}_{i\in I}$ of finite subgroups of $\Gamma_1$ such that $C_i\subset (A_1\overline{\otimes}M_2)\rtimes\Gamma_i$, for all $i\in I$, and $$L^2({\Cal M})\ominus L^2(M)\cong\oplus_{i\in I}L^2(\langle M,e_{C_i}\rangle),$$ as Hilbert $M$--bimodules.
\vskip 0.05in

\noindent{\it Proof of claim 1.} The proof of [CI08, Lemma 5] provides a family $\{A_i\}_{i\in I}$ of von Neumann subalgebras of $M_1$ and a family $\{\Gamma_i\}_{i\in I}$ of finite subgroups of $\Gamma_1$ such that $A_i\subset A\rtimes\Gamma_i$, for all $i\in I$, and $L^2(\tilde{ M_1})\ominus L^2(M_1)\cong\oplus_{i\in I}L^2(\langle M_1,e_{A_i}\rangle),$ as Hilbert $M_1$--bimodules. By letting $C_i=A_i\overline{\otimes} M_2\subset (A_1\rtimes\Gamma_i)\overline{\otimes}M_2=(A_1\overline{\otimes}M_2)\rtimes\Gamma_i$ and using the decomposition  $L^2({\Cal M})\ominus L^2(M)=(L^2(\tilde M_1)\ominus L^2(M_1))\overline{\otimes}L^2(M_2)$, the claim follows. \hfill$\square$

\vskip 0.05in

Now, recall that $\theta_t=\theta_t^1\otimes\theta_t^2\in$ Aut$(\tilde M=\tilde M_1\overline{\otimes}\tilde M_2)$ and that $\tilde M\supset{\Cal M}\supset M$.
\vskip 0.1in
 
\noindent
{\it Claim 2.} For every $x\in M$ and all $t\in\Bbb R$, we have that $||\theta_t'(x)-x||_2\leq 2||\theta_{\frac{t}{2}}(x)-x||_2$.
\vskip 0.05in

\noindent{\it Proof of claim 2.} Fix $x\in M$. Let us show first that $||E_{M}(\theta_t'(x))||_2\geq ||E_{M}(\theta_t(x))||_2$. Fix $t\in\Bbb R$ and for $i\in\{1,2\}$, denote by $T_i$ the bounded operator on $L^2(M_i)$ induced by $E_{M_i}\circ\theta_t^i:M_i\rightarrow M_i$.
Then $T_i$ is a contraction, hence $T_i^*T_i\leq 1$, as operators on $L^2(M_i)$. Also, we have that $E_{M}(\theta_t'(x))=(T_1\otimes 1)(x)$ and $E_{M}(\theta_t(x))=(T_1\otimes T_2)(x)$. Our assertion now follows from the following estimate 

\vskip 0.05in
$||(T_1\otimes T_2)(x)||_2^2=\langle (T_1^*T_1\otimes T_2^*T_2)(x),x\rangle\leq\langle (T_1^*T_1\otimes 1)(x),x\rangle=||(T_1\otimes 1)(x)||_2^2.$
\vskip 0.05in
 
Next, since $x\in M$, by using (3.g) and [Po08, Lemma 2.1.] we deduce that $||\theta_t'(x)-x||_2\leq 2||\theta_{\frac{t}{2}}'(x)-E_{M}( \theta_{\frac{t}{2}}'(x))||_2$. In combination with the above assertion we get that $$||\theta_t'(x)-x||_2^2\leq 4||\theta_{\frac{t}{2}}'(x)-E_{M}( \theta_{\frac{t}{2}}'(x))||_2^2=4(||\theta_{\frac{t}{2}}'(x)||_2^2-||E_{M}( \theta_{\frac{t}{2}}'(x))||_2^2)\leq $$ $$4(||x||_2^2-||E_{M}(\theta_{\frac{t}{2}}(x))||_2^2=4||x-E_{M}(\theta_{\frac{t}{2}}(x))||_2^2\leq 4||\theta_{\frac{t}{2}}(x)-x||_2^2, $$
and the claim is proven. \hfill$\square$

\vskip 0.05in
Since $x=(x_n)_n\in D'\cap M^{\omega}$ verifies $\lim_{t\rightarrow 0}||\theta_t(x)-x||_2=0$, by claim 2 we conclude that $\lim_{t\rightarrow 0}||\theta_t'(x)-x||_2=0$. Claim 1 says that we are in position to apply 
Lemma 3.3 and deduce that one of the following happens:  $(x_n)_n\in (A_1\overline{\otimes}M_2)^{\omega}\rtimes\Gamma_1$ or  there exists $t\not=0$ and $0\not=v\in\Cal M$ such that $\theta_t'(y)v=vy,$ for all $y\in D$. 
If the latter is true, since $\theta_t'=\theta_t^1\otimes$id$_{M_2}$, Theorem 2.1 and Remark 2.2 imply that $D\prec_{M}
L(\Gamma_1)\overline{\otimes}M_2$.  This contradicts our assumption that (2) is false so we must have that $$(x_n)_n\in (A_1\overline{\otimes}M_2)^{\omega}\rtimes\Gamma_1\tag 3.h$$

Now, let ${\Cal M}=M_1\overline{\otimes}\tilde M_2$ and $\theta_t'=$id$_{M_1}\otimes\theta_t^2\in$ Aut$({\Cal M})$. By arguing in the same way as above and using the fact that (2) is false, we get that $$(x_n)_n\in (M_1\overline{\otimes}A_2)^{\omega}\rtimes\Gamma_2,\tag 3.i$$

Finally, we show that (3.h) and (3.i) together imply that $(x_n)_n\in A^{\omega}\rtimes\Gamma$. 
For finite sets $F_1\subset\Gamma_1$ and $F_2\subset\Gamma_2$, denote by $Q_{F_1}$ and $R_{F_2}$ the orthogonal projections from $L^2(M)$ onto the closed linear span of $\{A_1u_{\gamma}\otimes M_2|\gamma\in F_1\}$ and of $\{M_1\otimes A_2u_{\gamma}|\gamma\in F_2\}$, respectively.
Then (3.h) and (3.i) can be rewritten as $$\inf_{F_1\subset\Gamma_1 \hskip 0.02in\text{finite}}\lim_{n\rightarrow\omega}||x_n-Q_{F_1}(x_n)||_2=0,\hskip 0.05in \inf_{F_2\subset\Gamma_2\hskip 0.02in \text{finite}}\lim_{n\rightarrow\omega}||x_n-R_{F_2}(x_n)||_2=0.$$

\noindent
Thus, we have that $\inf_{F_1\subset\Gamma_1,F_2\subset\Gamma_2\hskip 0.05in\text{finite}}\lim_{n\rightarrow\omega}||x_n-(R_{F_2}\circ Q_{F_1})(x_n)||_2=0$. Since $R_{F_2}\circ Q_{F_1}$ is the orthogonal projection from $L^2(M)$ onto the closed linear span of $\{(A_1\otimes A_2)u_{\gamma}|\gamma\in F_1\times F_2\}$, it follows that $(x_n)_n\in A^{\omega}\rtimes\Gamma$, therefore (1) holds true.
\hfill$\blacksquare$
\vskip 0.1in

\noindent
{\bf 3.4 Remark}. 
In the sequel, we will only use the following corollary of Theorem 3.2: let $D\subset M$ be an abelian von Neumann subalgebra such that condition (2) is false and suppose that $u_n\in L(\Gamma)$ is a sequence of unitaries which normalize $D$ and satisfy $||E_{L(\Gamma_1)}(au_nb)||_2,||E_{L(\Gamma_2)}(au_nb)||_2\rightarrow 0$, for every $a,b\in L(\Gamma)$.  Then for all $x\in D$ we have that $(u_nxu_n^*)_n\in A^{\omega}\rtimes_{\sigma}\Gamma$. Since $\theta_t(u_n)=u_n$, for all $t$ and $n$, this statement is indeed a consequence of 3.2. 

This corollary can be alternatively proved by adapting Popa's ``clustering coefficients" method (see [Po06b]). 
Following the strategy from [Po06b],  one first shows that the Fourier coefficients of $x_n=u_nxu_n^*$ cluster (i.e. they asymptotically belong to $B_1^{\Gamma_1\setminus F_1}\overline{\otimes}B_2^{\Gamma_2\setminus F_2}$ for any finite sets $F_1\subset\Gamma_1$ and $F_2\subset\Gamma_2$). Second, one uses the fact that the clustering sequence $(x_n)_n\in M^{\omega}$ commutes with elements $y\in D$ which are ``almost perpendicular" onto both $M_1\overline{\otimes}L(\Gamma_2)$ and $L(\Gamma_2)\overline{\otimes}M_2$ to conclude that $(x_n)_n\in A^{\omega}\rtimes_{\sigma}\Gamma$.

\vskip 0.05in
We end this section by showing that condition (2) in Theorem 3.2 can be exploited to get information on the quasi--normalizer of $D$. 
The next lemma is a straightforward generalization of [Po06a, Theorem 3.1.], but we include a  proof for the reader's convenience.

\proclaim {3.5 Proposition [Po06a]} Let $\Gamma\curvearrowright (X,\mu)$ be a p.m.p. mixing action and denote $M=L^{\infty}(X)\rtimes\Gamma$. Let $N$ be a finite von Neumann algebra.

\vskip 0.03in
\noindent
Assume that $D$ is a von Neumann subalgebra of $M\overline{\otimes}N$ such that $D\prec_{M\overline{\otimes}N}L(\Gamma)\overline{\otimes}N$. Denote by $P$ the von Neumann algebra generated by the quasi--normalizer of $D$ in $M\overline{\otimes}N$. If $D\nprec_{M\overline{\otimes}N}1\otimes N$, then $P\prec_{M\overline{\otimes}N}L(\Gamma)\overline{\otimes}N$.
\endproclaim

\noindent
{\it Proof.} We  first claim that if $p\in L(\Gamma)\overline{\otimes}N$ is a projection and $D\subset p(L(\Gamma)\overline{\otimes}N)p$ is a  unital von Neumann subalgebra such that $D\nprec_{ L(\Gamma)\overline{\otimes}N}1\otimes N$, then $q\Cal N_{p(M\overline{\otimes}N)p}(D)\subset p(L(\Gamma)\overline{\otimes}N)p$. By [Po06a, Section 3] (or [IPP08, proof of Theorem 1.1.]) in order to prove the claim it suffices to show that for every $\varepsilon>0$, $\eta_1,..,\eta_k\in (M\overline{\otimes}N)\ominus (L(\Gamma)\overline{\otimes}N)$, we can find $u\in\Cal U(D)$ such that $||E_{L(\Gamma)\overline{\otimes}N}(\eta_iu\eta_j^*)||_2\leq\varepsilon,$ for all $i,j\in\{1,..,k\}$.
Since $E_{L(\Gamma)\overline{\otimes}N}$ is $L(\Gamma)\overline{\otimes}N$--bimodular, by using Kaplansky's density theorem  it is enough to prove the last assertion for $\eta_i\in L^{\infty}(X)\otimes 1$ with $\tau(\eta_i)=0$ and $||\eta_i||\leq 1$. 

Next, decompose $u\in \Cal U(D)$ as $u=\sum_{\gamma}u_{\gamma}\otimes x_{\gamma}$, where $x_{\gamma}\in N$. Denote by $\alpha$ the induced action of $\Gamma$ on $L^{\infty}(X)$. We have that $$||E_{L(\Gamma)\overline{\otimes}N}(\eta_iu\eta_j^*)||_2^2=\sum_{\gamma\in\Gamma}|\tau(\eta_i\alpha_{\gamma}(\eta_j^*))|^2||x_{\gamma}||_2^2\tag 3.j$$

Since $\alpha$ is mixing and $\tau(\eta_i)=0$, for all $i$, we can find $F\subset\Gamma$ finite such that $$|\tau(\eta_i\alpha_{\gamma}(\eta_j^*))|\leq \frac{\varepsilon}{2},\forall i,j\in\{1,..,k\},\forall\gamma\in\Gamma\setminus F\tag 3.k$$ Since $D\nprec_{ L(\Gamma)\overline{\otimes}N}1\otimes N$, by Popa's theorem we can find $u\in\Cal U(D)$ such that $$\sum_{\gamma\in F}||x_{\gamma}||_2^2=\sum_{\gamma\in F}||E_{1\otimes N}(u(u_{\gamma}^*\otimes 1))||_2^2\leq\frac{\varepsilon}{2}\tag 3.l$$

As $|\tau(\eta_i\alpha_{\gamma}(\eta_j^*))|\leq ||\eta_i||\hskip 0.03in||\eta_j||\leq 1$, for all $\gamma\in\Gamma$, it is clear that the combination of (3.j), (3.k) and (3.l) implies that $||E_{L(\Gamma)\overline{\otimes}N}(\eta_iu\eta_j^*)||_2\leq\varepsilon,$ for all $i,j\in\{1,..,k\}$, as claimed.
\vskip 0.05in
Going back to the general case, let $D$ be a von Neumann  algebra such that $D\prec_{M\overline{\otimes}N}L(\Gamma)\overline{\otimes}N$ but $D\nprec_{M\overline{\otimes}N}1\otimes N$. By [Va08, Remark 3.8.] we can find non--zero projections $q\in D$, $p\in L(\Gamma)\overline{\otimes}N$, a $*$--homomorphism $\psi:qDq\rightarrow p(L(\Gamma)\overline{\otimes}N)p$ and a non--zero partial isometry $v\in p(M\overline{\otimes}N)q$ such that $\psi(x)v=vx$, for all $x\in qDq$, and $\psi(qDq)\nprec_{L(\Gamma)\overline{\otimes}N}1\otimes N$. By first part of the proof (applied to $\psi(qDq)$) we deduce that $\psi(qDq)'\cap p(L(\Gamma)\overline{\otimes}N)p\subset p(L(\Gamma)\overline{\otimes}N)p$, hence $vv^*\in p(L(\Gamma)\overline{\otimes}N)p$. This further implies that $vDv^*\subset L(\Gamma)\overline{\otimes}N$. Since $v^*v\in (qDq)'\cap q(M\overline{\otimes}N)q$, we can find $q'\in D'\cap (M\overline{\otimes}N)\subset P$ such that $v^*v=qq'$. Let $u$ be any unitary element extending $v$. Then  $u(qDqq')u^*=vDv^*\subset L(\Gamma)\overline{\otimes}N$.

Finally, since $D$ is quasi--regular in $P$, by [Po06a, Lemma 3.5.] we get that $qDqq'$ is quasi--regular in $qq'Pqq'$. Also, from the hypothesis we  have that no corner of $qDqq'$ embeds into $1\otimes N$. Thus,  by the first part of the proof, we deduce that $u(qq'Pqq')u^*\subset L(\Gamma)\overline{\otimes}N$, hence a corner of $P$ embeds into $L(\Gamma)\overline{\otimes}N$.\hfill$\blacksquare$ 

\vskip 0.2in

\head \S 4. {Lower bound on height.}\endhead
\vskip 0.2in
Let $M=A\rtimes\Gamma$ be a II$_1$ factor associated with a Bernoulli action. 
For every $v\in L(\Gamma)$, we define the {\it height} of $v$ as $h(v)=\max_{\gamma\in\Gamma}|\tau(vu_{\gamma}^*)|$.
Consider a sequence of elements $\{v_n\}_{n\geq 1}\subset (L(\Gamma))_1$. In this section, we provide a set of conditions which guarantee that the heights of the $v_n$'s  are uniformly bounded away from 0.

\proclaim {4.1 Theorem} 
Assume the notations from 3.1, i.e. $M_i=B_i^{\Gamma_i}\rtimes\Gamma_i$, $M=M_1\overline{\otimes}M_2$, $A=B_1^{\Gamma_1}\overline{\otimes}B_2^{\Gamma_2}$, $\Gamma=\Gamma_1\times\Gamma_2$.
Let $\{v_n\}_{n\geq 1}\subset (L(\Gamma))_1$ be a sequence for which there exist $C,c>0$ and $x\in M$ such that $C>\sqrt{2}c$, 
\vskip 0.03in
\noindent
(1) $||E_{A^{\omega}\rtimes\Gamma}((v_nxv_n^*)_n)||_2\geq C$,

\noindent
(2) $||E_{M_1\overline{\otimes}L(\Gamma_2)}(x)||_2\leq c$ and
 
\noindent
(3) $||E_{L(\Gamma_1)\overline{\otimes}M_2}(x)||_2\leq c.$
\vskip 0.03in
\noindent
Then $\lim_{n\rightarrow\omega}h(v_n)>0$.
\endproclaim

\noindent
{\it Proof.} Let $\varepsilon>0$ such that $C-\sqrt{2}c-3\varepsilon>0$ and let $x\in M$ satisfying conditions (1)--(3).  The first condition yields a finite set $F\subset\Gamma$ such that if $P_F$ denotes the orthogonal projection from $L^2(M)$ onto the closed linear span of $\{Au_{\gamma}|\gamma\in F\}$, then $$\lim_{n\rightarrow\omega}||P_F(v_nxv_n^*)||_2\geq C-\varepsilon\tag 4.a$$

Next, let $E_1=E_{M_1\overline{\otimes}L(\Gamma_2)}$ and $E_2=E_{L(\Gamma_1)\overline{\otimes}M_2}$. Since $E_1$ and $E_2$ commute,  $y=(($id$_{M}-E_1)\circ ($id$_{M}-E_2))(x)$ verifies $y\perp M_1\overline{\otimes}L(\Gamma_2)$ and $y\perp L(\Gamma_1)\overline{\otimes}M_2$. Thus, if we write $y=\sum_{\gamma\in\Gamma}a_{\gamma}u_{\gamma},$ where $a_{\gamma}\in A$, then $E_{B_1^{\Gamma_1}\otimes 1}(a_{\gamma})=E_{1\otimes {B_2^{\Gamma_2}}}(a_{\gamma})=0$, for all $\gamma\in\Gamma$. 
Also, conditions (2) and (3) yield $||x-y||_2=||E_1(x)+(E_2(x)-(E_1\circ E_2)(x))||_2\leq \sqrt{2}c.$

Let $K\subset\Gamma$ finite such that $z=P_K(y)$ satisfies $||z-y||_2\leq\varepsilon$. 
Let $S_1\subset\Gamma_1,S_2\subset\Gamma_2$ finite such that if $b_{\gamma}=E_{B_1^{S_1}\overline{\otimes}B_2^{S_2}}(a_{\gamma})$, then $w=\sum_{\gamma\in K}b_{\gamma}u_{\gamma}$ satisfies $||w-z||_2\leq\varepsilon.$
 The triangle inequality implies that $||w-x||_2\leq \sqrt{2}c+||p||_2+2\varepsilon$ and hence by (4.a) we have $$\lim_{n\rightarrow\omega}||P_{F}(v_nwv_n^*)||_2\geq (C-\varepsilon)-(\sqrt{2}c+2\varepsilon)>0\tag 4.b$$

Let $S=S_1S_1^{-1}\times S_2S_2^{-1}\subset\Gamma$, where $S_1S_1^{-1}=\{gh^{-1}|g,h\in S_1\}$.
We claim that $$\tau(\sigma_{\gamma}(b_g)b_h^*)=0,\hskip 0.03in\forall\gamma\in\Gamma\setminus S,\hskip 0.03in\forall g,h\in K\tag 4.c$$

  Indeed, let $\gamma=(\gamma_1,\gamma_2)\in\Gamma\setminus S$ and assume that $\gamma_1\notin S_1S_1^{-1}$ (the case $\gamma_2\notin S_2S_2^{-1}$ being analogous). Let $g,h\in K$. Since $b_h\in B_1^{S_1}\overline{\otimes}B_2^{\Gamma_2}$ and $\gamma_1S_1\cap S_1=\emptyset$, we have  $E_{B_1^{\gamma_1 S_1}\overline{\otimes}B_2^{\Gamma_2}}(b_h)=E_{1\otimes B_2^{\Gamma_2}}(b_h)=E_{1\otimes B_2^{S_2}}(a_h)=0$.  Since $\sigma_{\gamma}(b_g)\in B_1^{\gamma_1S_1}\overline{\otimes}B_2^{\Gamma_2}$, we derive that $\tau(\sigma_{\gamma}(b_g)b_h^*)=\tau(\sigma_{\gamma}(b_g)E_{B_1^{\gamma_1 S_1}\overline{\otimes}B_2^{\Gamma_2}}(b_h)^*)=0,$ which proves (4.c).
\vskip 0.05in

Next, fix $v\in L(\Gamma)$ and let us estimate $||P_F(vwv^*)||_2$. Let $\delta=\max_{g,h\in K}||b_g||_2||b_h||_2$ and write $v=\sum_{\gamma\in\Gamma}c_{\gamma}u_{\gamma}$, where $c_{\gamma}\in\Bbb C$, for all $\gamma\in\Gamma$. Then $$vwv^*=\sum_{\gamma,\gamma'\in\Gamma,g\in K}c_{\gamma}\overline{c_{\gamma'}}\sigma_{\gamma}(b_g) u_{\gamma g{\gamma'}^{-1}}=\sum_{h\in\Gamma} (\sum_{\gamma\in\Gamma,g\in K}c_{\gamma}\overline{c_{h^{-1}\gamma g}}\sigma_{\gamma}(b_g))u_h.$$

Thus, $$||P_F(vxv^*)||_2^2=\sum_{h\in F,\gamma,\gamma'\in\Gamma,g,g'\in K}c_{\gamma}\overline{c_{(h^{-1}\gamma g)}}\overline{c_{\gamma'}\overline{c_{(h^{-1}\gamma' g')}}}\tau(\sigma_{\gamma'}(b_{g'}^*)\sigma_{\gamma}(b_g)),$$ which by (4.c) is further equal to the real part of $$\sum_{\gamma\in\Gamma,s\in S,h\in F,g,g'\in K}c_{\gamma}\overline{c_{(h^{-1}\gamma g)}}\overline{c_{\gamma s^{-1}}}c_{(h^{-1}\gamma s^{-1} g')}\tau((b_{g'}^*)\sigma_{s}(b_g))$$ 
\noindent
Since $\Cal Re (a_1a_2a_3a_4)\leq\frac{1}{4}\sum_{i=1}^4|a_i|^4$, for every $\{a_i\}_{i=1}^{4}\in\Bbb C$, the last term is majorized by $$\frac{\delta}{4}\sum_{\gamma\in\Gamma,s\in S,h\in F,g,g'\in K}(|c_{\gamma}|^4+|c_{(h^{-1}\gamma g)}|^4+|c_{\gamma s^{-1}}|^4+|c_{(h^{-1}\gamma s^{-1} g')}|^4)\leq$$ $$\frac{\delta}{4}(1+|S|)(1+|F||K|)\sum_{\gamma\in\Gamma}|c_{\gamma}|^4\leq \frac{\delta}{4}(1+|S|)(1+|F||K|)\max_{\gamma\in\Gamma}|c_{\gamma}|^2\sum_{\gamma'\in\Gamma}|c_{\gamma'}|^2=$$ $$\frac{\delta}{4}(1+|S|)(1+|F||K|) h(v)^2||v||_2^2.$$
Altogether, we have shown that there exists a constant $\kappa>0$ such that $h(v)||v||_2\geq \kappa||P_F(vwv^*)||_2$, for all $v\in L(\Gamma)$. Since $||v_n||\leq 1$ we get  that $h(v_n)\geq \kappa||P_F(v_nwv_n^*)||_2$, for all $n$, and (4.b) implies the conclusion of the theorem.\hfill$\blacksquare$
\vskip 0.1in
\noindent
{\bf 4.2 Remark}. If $M=A\rtimes_{\sigma}\Gamma$ is a II$_1$ factor coming from a Bernoulli action then the above proof shows the following:
assume that $\{v_n\}_{n\geq 1}\subset (L(\Gamma))_1$ is a sequence for which there exists $x\in M$ such that 
$||E_{A^{\omega}\rtimes\Gamma}((v_nxv_n^*)_n)||_2>||E_{L(\Gamma)}(x)||_2$.
Then $\lim_{n\rightarrow\omega}h(v_n)>0$.

\vskip 0.2in

\head \S 5. {A conjugacy result for subalgebras.}\endhead

\vskip 0.2in
Let $M=B^{\Gamma}\rtimes\Gamma$ be a II$_1$ factor coming from a Bernoulli action. In this section, we give a criterion for proving that a von Neumann subalgebra $C$ of $M$ has a corner which embedds into $B^{\Gamma}$.  More generally, our criterion applies to subalgebras of $M\overline{\otimes}N$, where $N$ is an arbitrary finite von Neumann algebra. 
Before stating it we need to introduce some notation.

\vskip 0.03in
\noindent
{\bf Notations}. We  consider the orthogonal projections onto certain Hilbert subspaces of $L^2(M\overline{\otimes}N)$. For subsets $S$, $F$ and $G$ of $\Gamma$ we denote by
\vskip 0.02in
\noindent
$\bullet$ $P_S$ the orthogonal projection onto the closed linear span of $\{B^{\Gamma}u_{\gamma}\otimes N|\gamma\in S\}$.

\vskip 0.02in
\noindent
$\bullet$  $\Cal H_{F}$ the closed linear span of $\{B^Fu_{\gamma}\otimes N|\gamma\in\Gamma\}$. 
\vskip 0.02in
\noindent$\bullet$
 $Q_F$ the orthogonal projection  onto  $\Cal H_{F}$. 
\vskip 0.02in 
\noindent
$\bullet$ $Q_F^0:=Q_F-Q_{\emptyset}$.
\vskip 0.02in
\noindent
$\bullet$ $Q_F^G:=Q_{F\cup G}-Q_{G}$ the orthogonal projection onto $\Cal H_F^G:=\Cal H_{F\cup G}\ominus\Cal H_G$.
\vskip 0.05in
\noindent
Next, we record some useful boundedness and modularity properties of these projections.
\proclaim {5.1 Lemma} Let $S,F,G\subset\Gamma$ be finite sets. Then  $P_S$ commutes with $Q_F$, $Q_F^G$ and 
 $$||P_S(x)||\leq |S|, \hskip 0.03in||(P_S\circ Q_F)(x)||\leq |S|,\hskip 0.03in ||(P_S\circ Q_F^G)(x)||\leq 2|S|,\forall x\in (M\overline{\otimes}N)_1.$$
\noindent
If  $F$ contains $SG$, then $Q_F(P_S(x)y)=Q_F(P_S(x))y$ for every $x\in M\overline{\otimes}N$ and $y\in \Cal H_G$.
\endproclaim

\noindent
{\it Proof.} The commutativity assertion is trivial. Let $x\in M\overline{\otimes}N$ with $||x||\leq 1$. Write $M\overline{\otimes}N=(B^{\Gamma}\overline{\otimes}N)\rtimes_{\rho}\Gamma$, where $\Gamma$ acts through Bernoulli action on $B^{\Gamma}$ and trivially on $N$. Decompose $x=\sum_{\gamma\in\Gamma}x_{\gamma}u_{\gamma}$, where $x_{\gamma}=E_{B^{\Gamma}\overline{\otimes}N}(xu_{\gamma}^*)$. Then  $||x_{\gamma}||\leq 1$, for all $\gamma\in\Gamma$. Since $P_S(x)=\sum_{\gamma\in S}x_{\gamma}u_{\gamma}$, $(P_S\circ Q_F)(x)=\sum_{\gamma\in S}E_{B^{F}\overline{\otimes}N}(x_{\gamma})u_{\gamma}$ and $Q_F^G=Q_{F\cup G}-Q_{G}$, the three inequalities follow.

For the last assertion, it suffices to show that if $F\supset SG$ and if $x=\sum_{\gamma\in S}x_{\gamma}u_{\gamma}$ and $y=\sum_{g\in \Gamma}y_{g}u_{g}$ with $x_{\gamma}\in B^{\Gamma}\overline{\otimes}N$ and $y_g\in B^G\overline{\otimes}N$, for all $\gamma\in S$ and $g\in\Gamma$, then $Q_F(xy)=Q_F(x)y$. Notice that if $\gamma\in S$, then $\rho_{\gamma}(y_g)\in B^{\gamma G}\overline{\otimes}N\subset B^F\overline{\otimes}N$, for all $g\in \Gamma$. Thus, we have that $$Q_F(xy)=\sum_{\gamma\in S,g\in\Gamma}E_{B^{F}\overline{\otimes}N}(x_{\gamma}\rho_{\gamma}(y_g))u_{\gamma g}=\sum_{\gamma\in S,g\in\Gamma}E_{B^{F}\overline{\otimes}N}(x_{\gamma})\rho_{\gamma}(y_g)u_{\gamma g}=Q_F(x)y.$$ 

\noindent
We are now ready to state and prove the main result of this section.

\proclaim {5.2 Theorem} Let $\Gamma$ be a countable group, $B$ be a finite von Neumann algebra and denote $M=B^{\Gamma}\rtimes\Gamma$. Let $N$ be a finite von Neumann algebra and $p\in M\overline{\otimes}N$ be a projection.
\vskip 0.03in
\noindent
Let  $C\subset p(M\overline{\otimes}N)p$ be a unital von Neumann subalgebra. Suppose that there exist  a sequence $\{x_n\}_{n\geq 1}\subset\Cal U(p(M\overline{\otimes}N)p)$ and finite subsets $S$, $\{F_n\}_{n\geq 1}$ of $\Gamma$ such that
\vskip 0.02in
$\bullet$ $||x_ny-yx_n||_2\rightarrow 0$, for all $y\in C$, 

$\bullet$ $F_n\rightarrow\infty$, as $n\rightarrow\infty$ (i.e. if $\gamma\in\Gamma$ then $\gamma\notin F_n$,  for large enough $n$), 

$\bullet$ $\sup_{n\geq 1}|F_n|<\infty$, and

$\bullet$ $\limsup_{n\rightarrow\infty}(||x_n-Q_{F_n}^0(x_n)||_2+3||x_n-P_S(x_n)||_2) < ||p||_2.$
\vskip 0.02in
\noindent
Then $C\prec_{M\overline{\otimes}N}B^{\Gamma}\overline{\otimes}N$.

\endproclaim

\noindent
{\it Proof.}  By replacing $(x_n)_n$ with a subsequence, we may assume that
 $||x_n-Q_{F_n}^0(x_n)||_2\leq c_1$  and $||x_n-P_S(x_n)||_2\leq c_2$, for all $n$,  for some $c_1,c_2\geq 0$ satisfying $c_1+3c_2<||p||_2$. Let  $\varepsilon>0$ such that $c=||p||_2-(c_1+3c_2+11\varepsilon)>0.$
We begin with the following:
\vskip 0.05in
\noindent {\it Claim}. Let $y\in\Cal U(C)$ and  $G\subset\Gamma$ be a finite subset such that $||y-Q_{G}(y)||_2\leq \frac{\varepsilon}{|S|}$.
Let $K_n\subset\Gamma$ be a sequence of finite subsets such that $||y-P_{K_n}(y)||_2\leq\frac{\varepsilon}{|S|}$ and $K_nF_n$ is disjoint from $G\cup SG$. Then we can find $N\geq 1$  such that $||Q_{K_nF_n}^{G\cup SG}(x_n)||_2\geq c,$ for all $n\geq N$. 
\vskip 0.05in
\noindent
{\it Proof of claim.} Let $y_n=(P_{K_n}\circ Q_{G})(y)$, then $||y_n-y||_2\leq \frac{2\varepsilon}{|S|}$. Also, let $x_n'=(P_{S}\circ Q_{F_n}^0)(x_n)$. Since $P_{S}$ and $Q_{F_n}^0$ commute, we get that $$||x_n-x_n'||_2=||x_n-Q_{F_n}^0(x_n)+Q_{F_n}^0(x_n-P_{S}(x_n))||_2\leq c_1+c_2\tag 5.a$$
Now,   $y_nx_n'$ belongs to the closed linear span of $\{B^{G}u_{k}(B^{F_n}\ominus\Bbb C1)u_{\gamma})\otimes N|k\in K_n, \gamma\in S\}$.
Since $u_k(B^{F_n}\ominus\Bbb C1)\subset (B^{K_nF_n}\ominus\Bbb C1)u_k$, for $k\in K_n$, we get that  $y_nx_n'\in\Cal H_{K_nF_n}^G$ (recall that $\Cal H_F^G=\Cal H_{F\cup G}\ominus\Cal H_G$).  Since $K_nF_n$ is disjoint from $G\cup SG$, we get that  $ y_nx_n'\in \Cal H_{K_nF_n}^{G\cup SG}.$ 

Next, we estimate $||yx_n-y_nx_n'||_2$. Lemma 5.1 gives that $||x_n'||\leq 2|S|$, hence $||y-y_n||_2||x_n'||\leq 4\varepsilon$. By using the triangle inequality and (5.a) we get that $$||yx_n-y_nx_n'||_2\leq ||y||\hskip 0.03in ||x_n-x_n'||_2+||y-y_n||_2||x_n'||\leq c_1+c_2+4\varepsilon\tag 5.b$$

Since $y_nx_n'\in\Cal H_{K_nF_n}^{G\cup SG}$, (5.b) implies that $||yx_n-Q_{K_nF_n}^{G\cup SG}(yx_n)||_2\leq c_1+c_2+4\varepsilon$. Since $||[x_n,y]||_2\rightarrow 0$,  we can find $N$ such that $||x_ny-Q_{K_nF_n}^{G\cup SG}(x_ny)||_2\leq c_1+c_2+5\varepsilon,$ for all $n\geq N$.
Since $x_ny\in\Cal U(p(M\overline{\otimes}N)p)$, we get that $$||Q_{K_nF_n}^{G\cup SG}(x_ny)||_2\geq ||p||_2-(c_1+c_2+5\varepsilon)\tag 5.c$$

\noindent
Using that $||P_{S}(x_n)||\leq |S|$ we get that $||x_ny-P_{S}(x_n)y_n||_2\leq ||x_n-P_{S}(x_n)||_2||y||+||P_{S}(x_n)||\hskip 0.01in ||y-y_n||_2\leq c_2+|S|\frac{2\varepsilon}{|S|}=c_2+2\varepsilon.$ Combining this inequality with (5.c) yields that $$||Q_{K_nF_n}^{G\cup SG}(P_{S}(x_n)y_n)||_2\geq ||p||_2-(c_1+2c_2+7\varepsilon)\tag 5.d$$

Since $y_n\in\Cal H_{G}$, Lemma 5.1 gives that $Q_{K_nF_n}^{G\cup SG}(P_{S}(x_n)y_n)=Q_{K_nF_n}^{G\cup SG}(P_{S}(x_n))y_n$ and that $||Q_{K_nF_n}^{G\cup SG}(P_{S}(x_n))||\leq 2|S|$. Thus,  since $||y||\leq 1$, for all $n\geq N$ we have that $$||Q_{K_nF_n}^{G\cup SG}(x_n)||_2\geq ||Q_{K_nF_n}^{G\cup SG}(P_S(x_n))||_2-c_2\geq $$ $$||Q_{K_nF_n}^{G\cup SG}(P_S(x_n))y||_2-c_2\geq$$ $$||Q_{K_nF_n}^{G\cup SG}(P_S(x_n))y_n||_2-||Q_{K_nF_n}^{G\cup SG}(P_S(x_n))(y_n-y)||_2-c_2.$$

Since $||y_n-y||_2\leq \frac{2\varepsilon}{|S|}$, (5.d) implies that the last term is greater of equal than
 $||p||_2-(c_1+2c_2+7\varepsilon)-2|S|\frac{2\varepsilon}{|S|}-c_2=c,$ as claimed.\hfill$\square$
\vskip 0.1in
\noindent
To prove the theorem, assume  by contradiction that $C\nprec_{M\overline{\otimes}N}B^{\Gamma}\overline{\otimes}N$. By using the claim for all $m\geq 1$  we  will construct finite subsets $G_m$ and $\{K_{m,n}\}_{n\geq 1}$ of $\Gamma$ such that

(i) $K_{m,n}$ is disjoint from $G_1,..,G_m$, $K_{1,n},..,K_{m-1,n}$, for all $n\geq 1$,

(ii) $\sup_{n\geq 1}|K_{m,n}|<\infty$ and

(iii) $\liminf_{n\rightarrow\infty}||Q_{K_{m,n}}^{G_{m}}(x_n)||_2\geq c$. 
\vskip 0.05in
Before proving this statement let us show how it leads to a contradiction. First, (i) implies that for every $m\geq 1$, the projections $Q_{K_{1,n}}^{G_1},Q_{K_{2,n}}^{G_{2}},..,Q_{K_{m,n}}^{G_m}$ are mutually orthogonal (since $Q_{F}^GQ_{F'}^{G'}=0$, whenever $F$ is disjoint from $G$, $F'$ and $G'$). Second, for every $m\geq 1$ we have that $$||p||_2^2=\liminf_{n\rightarrow\infty}||x_n||_2^2\geq\sum_{l=1}^m\liminf_{n\rightarrow\infty}||Q_{K_{l,n}}^{G_l}(x_n)||_2^2.$$ By (iii) this implies that $||p||_2^2\geq mc^2$, for all $m\geq 1$,  a contradiction.
\vskip 0.05in
So, we are left with proving the above statement. We proceed by using induction on $m$. For $m=1$, if we let $K_{1,n}=F_n$ and $G_1=\emptyset$, then the statement is true by the hypothesis. Next, assume that we have constructed $G_l,\{K_{l,n}\}_{n\geq 1}$,  for all $l\leq m-1$. Since $\sup_{n\geq 1}|K_{l,n}|<\infty$, for all $l$, and $\sup_{n\geq 1}|F_n|<\infty$ we get that $$L=\sup_{n\geq 1,\hskip 0.02in l\leq m-1}|(G_l\cup K_{l,n})F_n^{-1})|<\infty$$

Recall that  $M\overline{\otimes}N=(B^{\Gamma}\overline{\otimes}N)\rtimes\Gamma$. Since $C\nprec_{M\overline{\otimes}N}B^{\Gamma}\overline{\otimes}N$, by Theorem 1.3.2 we can find $y\in\Cal U(C)$ such that for every finite set $T\subset \Gamma$ of cardinality $|T|\leq L$, we have that $$||P_T(y)||_2=(\sum_{\gamma\in T}||E_{B^{\Gamma}\overline{\otimes}N}(yu_{\gamma}^*)||_2^2)^{\frac{1}{2}}\leq \frac{\varepsilon}{2|S|}\tag 5.e$$

Now, let $K,G\subset\Gamma$ be finite sets such that $||y-Q_G(y)||_2\leq\frac{\varepsilon}{|S|}$ and $||y-P_K(y)||_2\leq\frac{\varepsilon}{2|S|}$. For every $n\geq 1$, define $K_n=K\setminus (\cup_{l\leq m-1}(G_l\cup K_{l,n})F_n^{-1})$. Since $|K\setminus K_n|\leq L$, by (5.e) we deduce that $||P_K(y)-P_{K_n}(y)||_2\leq \frac{\varepsilon}{2|S|}$ and thus $||y-P_{K_n}(y)||_2\leq \frac{\varepsilon}{|S|}$, for all $n\geq 1$. Since $F_n\rightarrow\infty$, as $n\rightarrow\infty$, we can find $s\geq 1$ such that $KF_n\cap (G\cup SG)=\emptyset$, for all $n\geq s$. Thus, $K_nF_n$ is disjoint from $G\cup SG$, for all $n\geq s$. 
Altogether, the above claim yields that $$\liminf_{n\rightarrow\infty}||Q_{K_nF_n}^{G\cup SG}(x_n)||_2\geq c\tag 5.f$$

Finally, set $K_{m,n}=K_nF_n$, for all $n\geq s$, $K_{m,n}=\emptyset$, for all $n\leq s-1$,  and $G_m=G\cup SG$. Then (5.f) can be rewritten as $
\liminf_{n\rightarrow\infty}||Q_{K_{m,n}}^{G_{m}}(x_n)||_2\geq c$, hence (iii) is verified. Since $|K_{m,n}|\leq |K||F_n|$ and $\sup_{n\geq 1}|F_n|<\infty$, we get that (ii) also holds true.  Also, by definition it is clear that $K_{m,n}$ is disjoint from $G_l$ and $K_{l,n}$, for all $1\leq l\leq m-1$. This proves the above statement and the theorem.\hfill$\blacksquare$

\vskip 0.2in
\head \S 6. {A dichotomy result for subalgebras.}\endhead
\vskip 0.2in
In this section, we  combine the results of the last four sections to prove Theorem F (which we restate below as Theorem 6.2.): if $D$ is an abelian subalgebra of  $M=B^{\Gamma}\rtimes\Gamma$ whose normalizer has ``large intersection" with $L(\Gamma)$, then a corner of $D$ embedds into either $B^{\Gamma}$ or $L(\Gamma)$. In the more general context, when $M$ is a tensor product of two factors associated with Bernoulli actions, we obtain:
\proclaim {6.1 Theorem} 
Assume the notations from 3.1, i.e. $M_i=B_i^{\Gamma_i}\rtimes_{\sigma_i}\Gamma_i$, $M=M_1\overline{\otimes}M_2$, $A=B_1^{\Gamma_1}\overline{\otimes}B_2^{\Gamma_2}$ and $\Gamma=\Gamma_1\times\Gamma_2$. 
Let  $D\subset qMq$ be an abelian von Neumann subalgebra, for some projection $q\in L(\Gamma)$.  Denote  $\Lambda=\Cal N_{qMq}(D)\cap \Cal U(qL(\Gamma)q)$ and assume that $\Lambda''\nprec_{M}L(\Gamma_1)\otimes 1$ and $\Lambda''\nprec_{M}1\otimes L(\Gamma_2)$.
\vskip 0.05in
\noindent
Then one of the following holds true:

\noindent
(1) $D'\cap qMq$ is of type I and there exist a unitary $u\in M$ and a projection  $q_0\in A$ such that $uq_0u^*=q$ and $u(Aq_0)u^*\subset D'\cap qMq$.
 
\noindent
(2) $D\prec_{M}M_1\overline{\otimes}L(\Gamma_2)$ or $D\prec_{M}L(\Gamma_1)\overline{\otimes}M_2$.
\endproclaim
\noindent
Recall that $\Lambda''$  denotes the von Neumann algebra generated by $\Lambda$ inside $qL(\Gamma)q$.
\vskip 0.05in
\noindent
{\it Proof.} We first prove the theorem in the case $q=1$.
Assume that (2) is false. We will show that (1) holds true. By using the hypothesis and  Theorem 1.3.1 (see the comment following it) we can find a sequence $\{u_n\}_{n\geq 1}\subset\Lambda$ such that $$||E_{L(\Gamma_1)}(au_nb)||_2\rightarrow 0\hskip 0.06in\text{and}\hskip 0.06in||E_{L(\Gamma_2)}(au_nb)||_2\rightarrow 0,\hskip 0.04in\forall a,b\in L(\Gamma)\tag 6.a$$

Next, we claim that  $$(u_nxu_n^*)_n\in A^{\omega}\rtimes_{\sigma}\Gamma,\hskip 0.04in\forall x\in D\tag 6.b$$ To prove (6.b), let $\tilde M$ and $\{\theta_t\}_{t\in\Bbb R}\subset$  Aut$(\tilde M)$ be defined as in the statement of Theorem 3.2. Since ${\theta_t}_{|L(\Gamma)}=$ id$_{|L(\Gamma)}$ and $u_n\in L(\Gamma)$, we get that $\theta_t$ converge uniformly to id$_{\tilde M}$ on $\{u_nxu_n^*\}_{n\geq 1}$. As $D$ is abelian and $u_n$ normalizes $D$, we deduce that $(u_nxu_n^*)_n\in D'\cap M^{\omega}$. Since (2) is assumed false, (6.b) follows from Theorem 3.2. 

\vskip 0.05in
Towards proving (1), let $C=D'\cap M$. Since $D$ is abelian, we have that  $C'\cap M\subset C$ or, equivalently,  $C'\cap M=\Cal Z(C)$. The main part of this proof consists of showing that $$Cp\prec_{M}B_1^{\Gamma_1}\overline{\otimes}M_2,\forall p\in \Cal P(\Cal Z(C))\tag 6.c$$  By symmetry, we will also get that $Cp\prec_{M}M_1\overline{\otimes}B_2^{\Gamma_2}$. Finally, we will combine these two statements to first get that $Cp\prec_{M}A$ and then derive (1).
\vskip 0.05in
First, we reduce (6.c) to a weaker statement. By Lemma 1.3.3, it suffices to prove (6.c) for projections $p\in\Cal Z(C)$ which commute with the normalizer of $C$. Since $\Lambda$ normalizes $C$, we get that $[p,\Lambda]=0$. Since $\Lambda''\nprec_{L(\Gamma)} L(\Gamma_1)$ and $\Lambda''\nprec_{L(\Gamma)}L(\Gamma_2)$, [Va08, Lemma 4.2.] implies that $\Lambda'\cap M\subset L(\Gamma)$, so in particular $p\in L(\Gamma)$. 
Indeed, notice that if $B_1\simeq B_2\simeq L^{\infty}(X_0,\mu_0)$, then the action of $\Gamma$ on $A$ can be identified with the generalized Bernoulli action $\Gamma\curvearrowright (X_0,\mu_0)^I$, where $I=\Gamma_1\sqcup\Gamma_2$ and $(\gamma_1,\gamma_2)\cdot g_1=\gamma_1g_1$, $(\gamma_1,\gamma_2)\cdot g_2=\gamma_2g_2$, for all $\gamma_1, g_1\in\Gamma_1$ and $\gamma_2, g_2\in\Gamma_2$. By applying [Va08, Lemma 4.2] to $I_1=\emptyset$, we get our assertion. For not necessarily isomorphic $B_1$ and $B_2$, it is not hard to adapt the proof of [Va08, Lemma 4.2] to prove our assertion.

Altogether, we get that it suffices to prove (6.c) for projections $p\in \Lambda'\cap L(\Gamma)\cap \Cal Z(C)$. 
\vskip 0.05in

For $n\geq 1$, let $v_n=u_np\in L(\Gamma)$.
Since (2) is false, we can find $x\in\Cal U(D)$ such that $||E_{M_1\overline{\otimes}L(\Gamma_2)}(xp)||_2\leq\frac{1}{2}||p||_2$ and $||E_{L(\Gamma_1)\overline{\otimes}M_2}(xp)||_2\leq \frac{1}{2}||p||_2$. Since $p$ commutes with $u_n\in\Lambda$ and $x\in D$, we get that $v_n(xp)v_n^*=u_nxu_n^*p$.
Since $(u_nxu_n^*)_n\in A^{\omega}\rtimes\Gamma$ and $||u_nxu_n^*p||_2=||p||_2$, for all $n$, by applying Theorem 4.1 we deduce that $\lim_{n\rightarrow\omega}h(v_n)>0.$ By replacing $u_n$ with a subsequence, (6.a) and (6.b) are preserved while we may assume that there is $\delta>0$ such that
 $$h(v_n)\geq \delta,\hskip 0.04in\forall n\geq 1\tag 6.d$$

For  $n\geq 1$, let $\gamma_n=(\gamma_n^1,\gamma_n^2)\in\Gamma=\Gamma_1\times\Gamma_2$ such that $|\tau(v_nu_{\gamma_n}^*)|=h(v_n)$.  We claim that $\gamma_n^1\rightarrow\infty$.
If not, then $\gamma_{n_k}^1=\gamma_1$, for some increasing subsequence $\{n_k\}_{k\geq 1}$ of $\Bbb N$ and some $\gamma_1\in\Gamma_1$.
 But this would imply that 
$$||E_{ L(\Gamma_2)}(u_{n_k}p(u_{\gamma_1}^*\otimes 1))|_2=||E_{ L(\Gamma_2)}(v_{n_k}(u_{\gamma_1}^*\otimes 1))||_2\geq\tau (v_{n_k}u_{\gamma_{n_k}}^*)|\geq\delta,\hskip 0.04in\forall k\geq 1$$
in contradiction with (6.a). Similarly, it follows that $\gamma_n^2\rightarrow\infty$.
\vskip 0.05in
\noindent
Given finite subsets  $F$, $S$ of $\Gamma_1$ and $T$ of $\Gamma_2$, we denote by

\noindent
$\bullet$ $\Cal K_F$ the closed linear span of $\{(B_1^{F}\ominus\Bbb C1)u_{\gamma_1}\otimes M_2|\gamma_1\in\Gamma_1\}$.

\noindent
$\bullet$ $Q_F^0$ the orthogonal projection from $L^2(M)$ onto $\Cal K_F$.

\noindent
$\bullet$ $P_{S}$ the orthogonal projection onto the closed linear span of $\{B_1^{\Gamma_1}u_{\gamma_1}\otimes M_2|\gamma_1\in S\}$.

\noindent
$\bullet$ $R_{T}$ the orthogonal projection onto the closed linear span of $\{M_1\otimes B_2^{\Gamma_2}u_{\gamma_2}|\gamma_2\in T\}$.

\vskip 0.1in

\noindent{\it Claim 1.} For every $x\in (M)_1$, every  subset $F\subset\Gamma_1$ and all $n\geq 1$, we have that 

$||v_nxv_n^*-Q_{\gamma_n^1F}^0(v_nxv_n^*)||_2\leq ||x-Q_F^0(x)||_2+\sqrt{||p||_2^2-\delta^2}$.
\vskip 0.1in
\noindent
{\it Proof of claim 1.} Recall that $\gamma_n=(\gamma_n^1,\gamma_n^2)$. Since  $u_{\gamma_n}\Cal K_{F}=\Cal K_{\gamma_n^1F}$, $v_n\in L(\Gamma)$  and $\Cal K_{F}$ is a right  $L(\Gamma)$--module, we get that $u_{\gamma_n}Q_F^0(x)v_n^*\in\Cal K_{\gamma_n^1F}$.  Then, since  $||\tau(v_nu_{\gamma_n}^*)u_{\gamma_n}||=|\tau(v_nu_{\gamma_n}^*)|\leq 1$, the triangle inequality gives that

$$||(\tau(v_nu_{\gamma_n}^*)u_{\gamma_n})Q_F^0(x)v_n^*-v_nxv_n^*||_2\leq $$ $$||(\tau(v_nu_{\gamma_n}^*)u_{\gamma_n})(Q_F^0(x)-x)v_n^*||_2+||(\tau(v_nu_{\gamma_n}^*)u_{\gamma_n}-v_n)xv_n^*||_2\leq$$ $$||Q_F^0(x)-x||_2+||\tau(v_nu_{\gamma_n}^*)u_{\gamma_n}-v_n||_2=||Q_F^0(x)-x||_2+(||v_n||_2^2-|\tau(v_nu_{\gamma_n}^*)|^2)^{\frac{1}{2}}.$$
Together with (6.d) this implies the claim.\hfill$\square$
\vskip 0.05in

Since (2) is assumed false we can find $x\in\Cal U(D)$ such that $||E_{L(\Gamma_1)\overline{\otimes}M_2}(x)||_2<\frac{1}{2}(||p||_2-\sqrt{||p||_2^2-\delta^2})$. If we enumerate $\Gamma_1=\{g_i\}_{i\geq 1}$ and let $F_n=\{g_i\}_{i=1}^n$, then $||E_{L(\Gamma_1)\overline{\otimes}M_2}(x)-(x-Q_{F_n}^0(x))||_2\rightarrow 0$.
Thus, we can find a finite subset $F\subset\Gamma_1$ such that $||x-Q_F^0(x)||_2<\frac{1}{2}(||p||_2-\sqrt{||p||_2^2-\delta^2})$.
 If we let $x_n=v_nxv_n^*$, by claim 1 we get that $$||x_n-Q_{\gamma_n^1F}^0(x_n)||_2< \frac{1}{2}(||p||_2+\sqrt{||p||_2^2-\delta^2}),\hskip 0.04in\forall n\geq 1\tag 6.e$$

Since $p$ commutes with $u_n$, $D$ and $u_n$ normalizes $D$, we have that $x_n=v_nxv_n^*=(u_nxu_n^*)p\in\Cal U(Dp)$, for all $n\geq 1$. Thus, we deduce that $(x_n)_n\in \Cal U((Cp)'\cap (pMp)^{\omega})$. By  (6.a) we have that $(x_n)_n=(u_nxu_n^*p)_n\in A^{\omega}\rtimes_{\sigma}\Gamma$. Hence, we can find a finite subset $S\subset\Gamma_1$ such that $\limsup_{n\rightarrow\infty}||x_n-P_S(x_n)||_2 < \frac{1}{6}(||p||_2-\sqrt{||p||_2^2-\delta^2})$. By combining this inequality with (6.e) we get that 

$$\limsup_{n\rightarrow\infty} \hskip 0.04in (||x_n-Q_{\gamma_n^1F}^0(x_n)||_2+3||x_n-P_S(x_n)||_2)<||p||_2.$$
\noindent
Altogether, Theorem 5.2 yields that $Cp\prec_{M} B_1^{\Gamma_1}\overline{\otimes}M_2$. This finishes the proof of (6.c).

\vskip 0.1in
\noindent
{\it Claim 2.} $Cp\prec_{M}A=B_1^{\Gamma_1}\overline{\otimes}B_2^{\Gamma_2}$, for all $p\in\Cal P(\Cal Z(C))$.
\vskip 0.05in
\noindent
{\it Proof of claim 2.} By (6.c), we have that $Cp\prec_{M}B_1^{\Gamma_1}\overline{\otimes}M_2$. Since $C'\cap M=\Cal Z(C)$, by Theorem 1.3.1 we can find non--zero projections $p_0\in Cp$ and $q\in B_1^{\Gamma_1}\overline{\otimes}M_2$, a $*$--homomorphism $\psi:p_0Cp_0\rightarrow q(B_1^{\Gamma_1}\overline{\otimes}M_2)q$ and a non--zero partial isometry $v\in qMp_0$ such that $v^*v=p_0$ and $\psi(x)v=vx$, for all $x\in p_0Cp_0$. Thus, we have that $$x=v^*\psi(x)v,\forall x\in p_0Cp_0\tag 6.f$$ 

Now, let $z\leq p$ be the central support of $p_0$ in $C$. By symmetry, the proof of (6.c) implies that $Cz\prec_{M}M_1\overline{\otimes}B_2^{\Gamma_2}$. By reasoning as in the previous paragraph, we can find non--zero projections $p_1\in Cz$ and $r\in M_1\overline{\otimes}B_2^{\Gamma_2}$, a $*$--homomorphism $\rho:p_1Cp_1\rightarrow r(M_1\overline{\otimes}B_2^{\Gamma_2})r$ and a non--zero partial isometry $w\in rMp_1$ such that $$x=w^*\rho(x)w,\forall x\in p_1Cp_1\tag 6.g$$

Since $p_0$ and $p_1$ admit equivalent, non--zero subprojections, we can assume that (6.f) and (6.g) hold true for $p_0=p_1\leq p$.
 By Kaplansky's density theorem we can find finite subsets $S\subset\Gamma_1$, $T\subset\Gamma_2$ and $v',w'\in (M)_1$ satisfying $||v'-v||_2,||w'-w||_2\leq\frac{||p_0||_2}{6}$ and $v'=P_S(v')$, $w'=R_T(w')$. Using (6.f), (6.g) and the triangle inequality it follows that $$||x-P_{S^{-1}S}(x)||_2,||x-R_{T^{-1}T}(x)||_2\leq  \frac{||p_0||_2}{3},\hskip 0.04in \forall x\in \Cal U(p_0Cp_0)\tag 6.h$$

Finally, (6.h) implies that $||x-(P_{S^{-1}S}\circ R_{T^{-1}T})(x)||_2\leq \frac{2||p_0||_2}{3}$, or equivalently, that $||(P_{S^{-1}S}\circ R_{T^{-1}T})(x)||_2\geq \frac{\sqrt{5}}{3}||p_0||_2$, for all $x\in\Cal U(p_0Cp_0)$. 
Since $P_{S^{-1}S}\circ R_{T^{-1}T}$ is precisely the orthogonal projection on the closed linear span of $\{Au_{\gamma}|\gamma\in S^{-1}S\times T^{-1}T\}$, by Theorem 1.3.1 we get that $p_0Cp_0\prec_{M}A$. Thus, $Cp\prec_{M}A$.\hfill$\square$
 \vskip 0.1in
\noindent
{\it Claim 3.} $C$ is a type I algebra and there exists a unitary $u\in M$ such that $uAu^*\subset C$.
\vskip 0.05in
\noindent{\it Proof of claim 3.} The first assertion follows from claim 2. For the second assertion, let $C_0$ be a maximal abelian subalgebra of $C$. Then $C_0$ is maximal abelian in $M$. Indeed, $C_0$ contains the center of $C$, hence it contains $D$ and therefore $C_0'\cap M\subset D'\cap M=C$.

Next, let $p\in\Cal P(\Cal Z(C))$. By claim 2, $Cp\prec_{M}A$ and therefore $C_0p\prec_{M}A$. Since $C_0p\subset pMp$ and $A\subset M$ are maximal abelian, by [Po06c, Theorem A.1.] (see also [Va07, Lemma C.3.]) we can find a non--zero projection $p'\in C_0p$ and a unitary $v\in M$ such that $C_0p'\subset vAv^*$. Let $z$ be the central support of $p'$ in $C$. Since $C$ is of type I, we have that $C_0$ is regular in $C$, i.e. $\Cal N_{C}(C_0)''=C$.
It follows that we can find projections $p_1,p_2,..\in C_0p'$ and $u_1,u_2,..\in\Cal N_{C}(C_0)$ such that $z=\sum_{i\geq 1}u_ip_iu_i^*$.

Further, since $u_i\in\Cal N_{C}(C_0)$ we get that $C_0(u_ip_iu_i^*)=u_i(C_0p_i)u_i^*\subset(u_iv)A(u_iv)^*$, for all $i$. Since $A$ is regular in $M$ and $M$ is a II$_1$ factor, it follows that we can find a unitary $w\in M$ such that $C_0z\subset wAw^*$. 
We have altogether shown that for every non--zero projection $p\in \Cal Z(C)$, there is a non--zero projection $z\in\Cal Z(C)$ with $z\leq p$ such that $C_0z$ can be unitarily conjugated inside $A$.

Finally, let $\Cal S$ be the set of  families $\{p_i\}_{i\in I}\subset \Cal P(\Cal Z(C))$ of mutually orthogonal projections with the property that   $C_0p_i\subset u_iAu_i^*$, for some $u_i\in\Cal U(M)$, for all $i\in I$. By the above we get that if $\{p_i\}_{i\in I}$ is maximal element in $\Cal S$ with respect to inclusion, then $\sum_{i\in I}p_i=1$. Using again the fact that $A$ is regular in $M$, 
we deduce that $C_0$ can be unitarily conjugated inside $A$. Since $C_0$ is maximal abelian, we get that $C_0=uAu^*$, for some $u\in\Cal U(M)$.  This ends the proof of claim 3 and of the case $q=1$.\hfill$\square$

\vskip 0.1in
In general, let $D\subset qMq$ be a subalgebra as in the hypothesis and assume that (2) is false. 
Then (1) follows true by repeating the above proof once we show that (6.b) holds true in this context, i.e. $(u_nxu_n^*)_n\in A^{\omega}\rtimes_{\sigma}\Gamma$, for all $x\in D$ and any sequence $\{u_n\}_{n\geq 1}\subset\Lambda$. To see this, let $D_0\subset (1-q)M(1-q)$ be a unital, abelian von Neumann subalgebra such that $D_0\nprec_{M}M_1\overline{\otimes}L(\Gamma_2)$ and $D_0\nprec_{M}L(\Gamma_1)\overline{\otimes}M_2$. Then $D_1=D\oplus D_0$ has the same property. Since $D_1\subset M$ is unital  and $\theta_t$ converge uniformly to the identity on $(u_nxu_n^*)_n\in D_1'\cap M^{\omega}$, our claim follows from Theorem 3.2.
\hfill$\blacksquare$

\vskip 0.1in

Note that all of our results concerning tensor products of II$_1$ factors associated with Bernoulli actions (Theorems 3.3, 4.1,  6.1) have straightforward counterparts for single such II$_1$ factors. In this case, Theorem 6.1 reads as follows:

\proclaim {6.2 Theorem} Let $\Gamma$ be a  countable group, $B$ be an abelian von Neumann algebra and set $M=B^{\Gamma}\rtimes_{\sigma}\Gamma$, $A=B^{\Gamma}$. Let $q\in L(\Gamma)$ be a projection and $D\subset qMq$ be a unital abelian von Neumann subalgebra. Suppose that the group of unitary elements $u\in qL(\Gamma)q$ that normalize $D$ generates a diffuse von Neumann algebra. Then either: 
\vskip 0.03in
\noindent
(1) $D'\cap qMq$ is of type I and there exist a unitary $u\in M$ and a projection $q_0\in A$ such that $uq_0u^*=q$ and $u(Aq_0)u^*\subset D'\cap qMq$ or
\vskip 0.03in 
\noindent
(2) $D\prec_{M}L(\Gamma)$.
\endproclaim

We end this section by noticing that Theorem 6.2 can be used to give a proof of Popa's conjugacy criterion for  Bernoulli actions.

\proclaim {6.3 Theorem (Popa, [Po06b, Theorem 0.7.])} Let $\Gamma$ be a countable ICC group, $B$ be a non--trivial abelian von Neumann algebra. Denote  $M=B^{\Gamma}\rtimes\Gamma$, $A=B^{\Gamma}$ and let $q\in L(\Gamma)$ be a projection. 
Let $\rho:\Lambda\rightarrow$ Aut$(C)$ be a free ergodic action of a countable group $\Lambda$ on an abelian von Neumann algebra $C$. Denote $N=C\rtimes\Lambda$. 

\vskip 0.01in
\noindent
Let $\theta:N\rightarrow qMq$ be a $*$--isomorphism. 
Assume that $\theta(L(\Lambda))\subset qL(\Gamma)q$.
\vskip 0.05in
\noindent
Then $q=1$ and there exist a unitary $u\in M$, a character $\eta$ of $\Lambda$, a group isomorphism $\delta:\Lambda\rightarrow\Gamma$ 
such that $\theta(C)=uAu^*$ and $\theta(v_{\lambda})=\eta(\lambda)uu_{\delta(\lambda)}u^*$, for all $\lambda\in\Lambda$.
\endproclaim
\noindent
{\it Proof.} Denote $D=\theta(C).$ Since $\sigma$ is mixing and $D$ is regular in $qMq$, [Po06a, Theorem 3.1.] implies that $D\nprec_{M}L(\Gamma)$. Since $D$ is normalized by $\{\theta(v_{\lambda})\}_{\lambda\in\Lambda}\subset qL(\Gamma)q$ and $D$ is maximal abelian in $qMq$, Theorem 6.2 	implies that there exists a unitary $u\in M$ and a projection $q_0\in A$ such that $uq_0u^*=q$ and $D=u(Aq_0)u^*$. The conclusion now follows from  [Po06b, Theorem 5.2.].\hfill$\blacksquare$

\vskip 0.2in
\head \S 7. {Popa's conjugacy criterion for actions}.\endhead
\vskip 0.2in
In the proof of the strong rigidity results we will need the following version of Popa's conjugacy criterion for actions.  

\proclaim {7.1 Theorem [Po06b]} Let $\sigma:\Gamma\rightarrow$ Aut$(A)$ and $\beta:\Lambda\rightarrow$ Aut$(B)$ be two free, ergodic actions of two countable groups $\Gamma$ and $\Lambda$ on two abelian von Neumann algebras $A$ and $B$. Assume that $\beta$ is weakly mixing.  Denote $M=A\rtimes_{\sigma}\Gamma$ and $N=B\rtimes_{\beta}\Lambda$.  Let $\{u_{\gamma}\}_{\gamma\in\Gamma}\subset M$ and $\{v_{\lambda}\}_{\lambda\in\Lambda}\subset N$ be the canonical unitaries implementing $\sigma$ and $\beta$.

\vskip 0.01in
\noindent
Let $p\in L(\Gamma)$ be a non--zero projection and assume that $N$ is embedded (unitaly) inside $pMp$ such that $L(\Lambda)\subset pL(\Gamma)p$. Suppose that there exists a partial isometry $v\in M$ such that $v^*v=p$, $vv^*=q\in A$ and $vBv^*=Aq$.
\vskip 0.05in
\noindent
Then we can find maps $h:\Lambda\rightarrow S^1$, $\delta:\Lambda\rightarrow\Gamma$ and a unitary $z\in L(\Gamma)$ such that $v_{\lambda}=h(\lambda)(zu_{\delta(\lambda)}z^*)p$ and $[zu_{\delta(\lambda)}z^*,p]=0,$ for all $\lambda\in\Lambda$. Let $K=\{\gamma\in\Gamma|(zu_{\gamma}z^*)p\in\Bbb Cp\}$ and $\Gamma'=\{\delta(\lambda)k|\lambda\in\Lambda,k\in K\}$. Then $K,\Gamma'$ are subgroups of $\Gamma$, $|K|<\infty$ and $\Gamma'$ normalizes $K$.
Moreover, we have that:

\vskip 0.02in
\noindent
$\bullet$ If $\Gamma$ is torsion free, then  $h$ is a character and $\delta$ is a homomorphism. If
 in addition the restriction of $\sigma$ to $\delta(\Lambda)$ is ergodic, then $p=q=1$ and $B=zAz^*$.
\vskip 0.02in
\noindent
$\bullet$ In general, if the restriction of $\sigma$ to $\Gamma'$ is weakly mixing, then  $z^*pz\in L(K)\cap\Cal Z(L(\Gamma'))$, $z^*pz=|K|^{-1}\sum_{k\in K}\chi(k)u_k$, for some character $\chi$ of $K$, and $B=(zA^{K}z^*)p$, where $A^K=\{a\in A|\sigma_k(a)=a,\forall k\in K\}$. In particular, $\tau(p)^{-1}=|K|$ is an integer.

\endproclaim  
\vskip 0.05in

The statement of this result is very close to that of Theorem 5.2. in [Po06b].  However, rather than assuming that $\sigma$ is mixing, as in [Po06b], we only require that $\sigma_{|\delta(\Lambda)}$ is ergodic (in the torsion free case) and weakly mixing (in general). This generalization will be essential for us, as we will later apply Theorem 7.1 to the action $\sigma$ defined in 3.1, for which the ``global" condition $\sigma$ {\it is mixing} fails. On the other hand, since $\sigma$ is a concrete action (i.e. a product of two Bernoulli actions) it will be easy to verify the ``local" condition {\it $\sigma_{|\delta(\Lambda)}$ is weakly mixing}.

\vskip 0.1in
\noindent
{\it Proof.} In the beginning of the proof we repeat part of Vaes' proof of Popa's criterion (see [Va07, Proposition 9.3.]) which streamlines the original argument from [Po06b].
Identify $A$ with $L^{\infty}(X,\mu)$, where $(X,\mu)$ is a probability space. Let $Y\subset X$ be a Borel set such that $q=1_{Y}$. 
Note that $vv_{\lambda}v^*\in \Cal U(qMq)$ normalizes $Aq\cong L^{\infty}(Y)$, for all $\lambda\in\Lambda$, and denote by $\beta'$ the associated action of $\Lambda$ on $Y$. Since $\beta'$ is  isomorphic to $\beta$, we have that $\beta'$ is weakly mixing.
We consider the canonical embedding $\eta:L^{\infty}(X)\rtimes_{\sigma}\Gamma\subset L^{\infty}(X,\ell^2(\Gamma))$ given by $\eta(au_{\gamma})(x)=a(x)u_{\gamma}$, for all $a\in L^{\infty}(X),x\in X$ and $\gamma\in\Gamma$, as well as the  embedding $L(\Gamma)\subset\ell^2(\Gamma)$. 
Also, we view $S^1\times \Gamma$ as a subgroup of $\Cal U(L(\Gamma))$, in the natural way. 

Let $w=\tau(p)^{\frac{1}{2}}\eta(v)\in L^{\infty}(Y,\ell^2(\Gamma))$.
 Since $\beta$ is weakly mixing, [Po06b, Lemma 4.4.$(i)$] implies that $B=v^*Av$ is orthogonal to $L(\Gamma)$. The same argument as in the proofs of [Po06b, Lemma 6.1.] and [Va07, Proposition 9.3.] shows that  every essential value $w_0$ of the function $w:Y\rightarrow \ell^2(\Gamma)$ lies in $L(\Gamma)$ and satisfies $w_0^*w_0=p$. Set $v'=vw_0^*$, $p'={v'}^*v'=w_0pw_0^*=w_0w_0^*$, $w'=\tau(p)^{\frac{1}{2}}\eta(v')$ and notice that $p'$ is an essential value of $w'$.
By replacing $v$ with $v'$, $p$ with $p'$, $v_{\lambda}$ with $w_0v_{\lambda}w_0^*$ and $B$ with $w_0Bw_0^*$, we may assume that $p$ is an essential value of $w$, while the hypothesis is still satisfied. For the ``new" $v,p,v_{\lambda}$ we will prove  the conclusion for $z=1$, which clearly finishes the proof.

Let $\lambda\in\Lambda$.  Since $vv_{\lambda}v^*\in\Cal N_{qMq}(Aq)$, we can find $a_{\lambda}\in \Cal U(L^{\infty}(Y))$ such that $vv_{\lambda}v^*=a_{\lambda}\sum_{\gamma\in\Gamma}1_{\{y|\lambda^{-1}y=\gamma^{-1}y\}}u_{\gamma}$, where we denote  $\lambda y=\beta'_{\lambda}(y)$, for every $y\in Y$.
Let $\omega=\sum_{g\in\Gamma}b_{g}u_{g}\in M$, where $b_{g}\in A$. Then $(vv_{\lambda}v^*)\omega=\sum_{\gamma,g\in\Gamma}a_{\lambda}1_{\{y|\lambda^{-1}y=\gamma^{-1}y\}}(b_{g}\circ{\gamma}^{-1})u_{\gamma g}$. Thus, if we let $c(\lambda,y)=a_{\lambda}(y)u_{\gamma}\in S^1\times\Gamma\subset\Cal U(L(\Gamma))$, where $\gamma\in\Gamma$  is the unique element such that $\gamma^{-1}y=\lambda^{-1}y$, then
$$\eta((vv_{\lambda}v^*)\omega)(y)=\sum_{\gamma,g\in\Gamma}a_{\lambda}(y)1_{\{y|\lambda^{-1}y=\gamma^{-1}y\}}(y)b_{g}({\gamma}^{-1}y)u_{\gamma g}=\tag 7.a$$ $$(a_{\lambda}(y)u_{\gamma})(\sum_{g\in\Gamma}b_{g}({\lambda}^{-1}y)u_{g})=c(\lambda,y)\eta(\omega)(\lambda^{-1} y),\hskip 0.03in\text{for almost all}\hskip 0.03in y\in Y.$$

Since $(vv_{\lambda}v^*)v=vv_{\lambda}$ and $v_{\lambda}\in L(\Gamma)$, by applying $\eta$ and using (7.a), we get that $$c(\lambda,y)w(\lambda^{-1} y)=w(y)v_{\lambda},\hskip 0.1in\text{for almost all $y\in Y$} \tag 7.b$$

\vskip 0.05in

Next, we use an argument due to Popa (see the proofs of [Po07a, Proposition 3.5.] and [Va07, Lemma 4.8.]) to conclude the first assertion of the theorem. If $\varepsilon>0$, then, since $p$ is an essential value of $w$, the set $W=\{y\in Y|\hskip 0.03in||w(y)-p||_2\leq\varepsilon\}$ has positive measure. Fix $\lambda\in \Lambda$. Since $\beta'$ is weakly mixing, we can find $g\in\Lambda$ such that $\mu(gW\cap W)>0$ and $\mu((g\lambda^{-1})W\cap W)>0$.
Let $y\in gW\cap W$. Then $||w(y)-p||_2, ||w(g^{-1}y)-p||_2\leq\varepsilon$ and by using (7.b) we deduce that $$||v_g-c(g,y)p||_2=||pv_g-c(g,y)p||_2=$$ $$||(p-w(y))v_g-c(g,y)(p-w(g^{-1}y))||_2\leq 2\varepsilon.$$

Similarly,  if $z\in (g\lambda^{-1})W\cap W$, then $||v_{g {\lambda}^{-1}}-c(g\lambda^{-1},z)p||_2\leq 2\varepsilon$. 
Since $v_{\lambda}=(v_{g\lambda^{-1}})^*v_g$, the combination of the last two inequalities yields  $$||v_{\lambda}-p c(g\lambda^{-1},z)^*c(g,y)p||=||v_{g\lambda^{-1}}^*v_g-(c(g\lambda^{-1},z)p)^*c(g,y)p||_2\leq 4\varepsilon.$$
Therefore, since $c(g\lambda^{-1},z)^*c(g,y)\in S^1\times\Gamma$ and $\varepsilon>0$ is arbitrary, we can find two sequences $\{t_n\}_{n\geq 1}\subset S^1$ and $\{\gamma_n\}_{n\geq 1}\subset\Gamma$ such that $\lim_{n\rightarrow\infty}||v_{\lambda}-p(t_nu_{\gamma_n})p||_2=0$. By passing to a subsequence we can assume that $t_n\rightarrow t\in S^1$ and either $\gamma_n\rightarrow\gamma\in\Gamma$ or $\gamma_n\rightarrow\infty$, as $n\rightarrow\infty$. If $\gamma_n\rightarrow\infty$, then $t_nu_{\gamma_n}$ converges weakly to $0$, which would imply that $v_{\lambda}=0$, a contradiction. Thus, we get that $\gamma_n\rightarrow\gamma\in\Gamma$, hence $v_{\lambda}=tpu_{\gamma}p$. Since $v_{\lambda}\in\Cal U(pMp)$, we must have that $u_{\gamma}$ commutes with $p$. We altogether deduce that there exist maps $h:\Lambda\rightarrow S^1$  and $\delta:\Lambda\rightarrow\Gamma$ such that $$[u_{\delta(\lambda)},p]=0\hskip 0.1in\text{and}\hskip 0.1in v_{\lambda}=h(\lambda)u_{\delta(\lambda)}p,\forall \lambda\in\Lambda\tag 7.c$$

As in the hypothesis, denote $K=\{\gamma\in\Gamma|u_{\gamma}p\in\Bbb Cp\}$ and $\Gamma'=\{\delta(\lambda)k|\lambda\in\Lambda,k\in K\}$. Then  $K$ is  finite  and (7.c)  implies that $\Gamma'$ is a subgroup of $\Gamma$ which contains $K$ as a normal subgroup. Notice also that the map $\Lambda\ni\lambda\rightarrow\delta(\lambda)K\in\Gamma'/K$ is a group isomorphism.

\vskip 0.1in
We can now prove the moreover assertions.
\vskip 0.05in
\noindent
$\bullet$ If $\Gamma$ is {\it torsion free}, then $K=\{e\}$ and (7.c) implies that $\delta$ is a homomorphism and $h$ is character.
Further, assume that $\sigma_{|\delta(\Lambda)}$ is ergodic and let $\varepsilon>0$. Let $W_{\varepsilon}$ be the set of $y\in Y$ for which there exists $u\in S^1\times\Gamma\subset \Cal U(L(\Gamma))$ such that $||w(y)-up||_2\leq\varepsilon$. Since $p$ is an essential value of $w$, we have that $\mu(W_{\varepsilon})>0$. On the other hand, (7.b) and (7.c) imply that $W_{\varepsilon}$ is $\Lambda$--invariant. Indeed, if $y\in Y,\lambda\in\Lambda$ and $u\in S^1\times\Gamma$, then $u'=h(\lambda)c(\lambda,y)^*uu_{\delta(\lambda)}\in S^1\times\Gamma$ satisfies $$||w(\lambda^{-1} y)-u'p||_2=||c(\lambda,y)^*w(y)v_{\lambda}-h(\lambda)c(\lambda,y)^*uu_{\delta(\lambda)}p||=$$ $$||c(\lambda,y)^*w(y)v_{\lambda}-{c(\lambda,y)}^*upv_{\lambda}||_2=||w(y)-up||_2.$$ Since  $\beta'$ is ergodic, we get that $W_{\varepsilon}=Y$. By reasoning as above we conclude that there exists a measurable map $u:Y\rightarrow S^1\times \Gamma$ such that $w(y)=u(y)p$, for a.e. $y\in Y$.

 But then (7.b) rewrites as $c(\lambda,y)u(\lambda^{-1} y)p=h(\lambda)u(y)u_{\delta(\lambda)}p,$ for all $\lambda\in\Lambda$ and for almost all $y\in Y$. Since $c(\lambda,y)u(\lambda^{-1}y),h(\lambda)u(y)u_{\delta(\lambda)}\in S^1\times\Gamma$, by using the fact that $K=\{e\}$, we get that $$c(\lambda,y)u(\lambda^{-1} y)=h(\lambda)u(y)u_{\delta(\lambda)},\forall\lambda\in\Lambda \hskip 0.05in\text{and for almost all}\hskip 0.05in y\in Y\tag 7.d$$

Finally, let $\omega\in L^2(A\rtimes_{\sigma}\Gamma)\cong L^2(X,\ell^2(\Gamma))$ be defined by $\omega(x)=u(x)$, if $x\in Y$ and $\omega(x)=0$, otherwise. By using (7.a) and  (7.d) we get that $(vv_{\lambda}v^*)\omega=h(\lambda)\omega u_{\delta(\lambda)},$ for all $\lambda\in\Lambda$.
Since $\omega=q\omega$ and $vv_{\lambda}v^*\in\Cal U(qMq)$, we further derive that $$\omega^*\omega=u_{\delta(\lambda)}(\omega^*\omega)u_{\delta(\lambda)}^*,\forall\lambda\in\Lambda\tag 7.e$$
By the definition of $\omega$, we can find $a\in\Cal U(A)$ and a measurable partition $\{A_{\gamma}\}_{\gamma\in\Gamma}$ of $Y$ such that $\omega=a\sum_{\gamma\in\Gamma}1_{A_{\gamma}}u_{\gamma}$. Hence, we have that $\omega^*\omega=\sum_{\gamma\in\Gamma}u_{\gamma}^*1_{A_{\gamma}}u_{\gamma}=\sum_{\gamma\in\Gamma}1_{\gamma^{-1}A_{\gamma}}\in A.$ Since $\sigma_{|\delta(\Lambda)}$ is ergodic, by (7.e) we deduce that $\sum_{\gamma\in\Gamma}1_{\gamma^{-1}A_{\gamma}}\in\Bbb C1$. This clearly implies that $Y=X$, hence $p=q=1$. Also, we get that $w(y)=u(y)\in S^1\times\Gamma$, for almost all $y\in X$,  thus $v\in\Cal N_{M}(A)$. This shows that $B=A$, which concludes the proof of the torsion free case.

\vskip 0.1in
\noindent
$\bullet$ The  proof  of the {\it general case} is an adaptation of arguments from [Po06b, Section 6].  Assume that $\sigma_{|\Gamma'}$ is weakly mixing. Since  $p$ is an essential value for $w=\tau(p)^{\frac{1}{2}}\eta(v)$, we can find a decreasing sequence of non--zero projections $\{q_n\}_{n\geq 1}\subset Aq$ such that $$||q'v-\tau(p)^{-\frac{1}{2}} q'p||_2\leq 2^{-n-1}||q'||_2, \forall q'\in\Cal P(Aq_n)\tag 7.f$$ and every $n\geq 1$ (see [Po06b, Lemma 6.1.]).

Let us introduce some notation. Fix $n\geq 1$ and $\lambda\in\Lambda$. Since $p,v_{\lambda}\in L(\Gamma)$ we can decompose
 $p=\sum_{\gamma\in\Gamma}c_{\gamma}u_{\gamma}$ and $v_{\lambda}=\sum_{\gamma\in\Gamma}c_{\gamma}^{\lambda}u_{\gamma}$, with $c_{\gamma},c_{\gamma}^{\lambda}\in\Bbb C$, for all $\gamma\in\Gamma$. Thus, we have that $q_nv_{\lambda}q_n=\sum_{\gamma\in\Gamma}c_{\gamma}^{\lambda}\sigma_{\gamma}(q_n)q_nu_{\gamma}$. Also, we decompose $q_nvv_{\lambda}v^*q_n=\sum_{\gamma\in\Gamma}a_{\gamma}^{\lambda,n}u_{\gamma}$, with $a_{\gamma}^{\lambda,n}\in Aq$.

Let $F_n^{\lambda}=\{\gamma\in\Gamma|a_{\gamma}^{\lambda,n}\not=0\}$.  By using (7.f), [Po06b, Lemma 6.3] gives that $$|a_{\gamma}^{\lambda,n}|=\sigma_{\gamma}(q_n)q_n\hskip 0.05in\text{and}\hskip 0.05in |1-\tau(p)^{-1}|c_{\gamma}^{\lambda}||\leq 2^{-n},\forall \gamma\in F_n^{\lambda}\tag 7.g $$

 Note that since $v_{\lambda}=h(\lambda)u_{\delta(\lambda)}p$ and $|h(\lambda)|=1$, we get that $$ |c_{\gamma}^{\lambda}|=|c_{\delta(\lambda)^{-1}\gamma}|,\forall\gamma\in\Gamma,\lambda\in\Lambda\tag 7.h$$

Next, we prove the following:
\vskip 0.05in
\noindent
{\it Claim}. $\tau(p)^{-1}\in\Bbb N$  and $|c_{\gamma}|\in\{0,\tau(p)\}$, for all $\gamma\in\Gamma$.
\vskip 0.05in
\noindent
{\it Proof of claim.} Let $n\geq 1$ and $\varepsilon>0$. Set $F=\{\gamma\in\Gamma|\hskip 0.03in |c_{\gamma}|\geq\frac{\tau(p)}{2}\}$. By (7.g), if $\gamma\in F_n^{\lambda}$ then $|c_{\gamma}^{\lambda}|\geq\frac{\tau(p)}{2}$. Using (7.h) this shows that $\delta(\lambda)^{-1}\gamma\in F$. In other words,  $F_n^{\lambda}\subset\delta(\lambda)F,$ for all $\lambda\in\Lambda$.

Now, since $vv_{\lambda}v^*$ normalizes $Aq$ we get that $a_{\gamma}^{\lambda,n}$ are partial isometries in $Aq$. Thus, using the definitions of $F_n^{\lambda}$, $a_{\gamma}^{\lambda,n}$ and the fact that $v^*q_nv\in v^*(Aq)v=B$ we derive that $$\sum_{\gamma\in F_n^{\lambda}}\tau(|a_{\gamma}^{\lambda,n}|)=\sum_{\gamma\in\Gamma}\tau(|a_{\gamma}^{\lambda,n}|)=||q_nvv_{\lambda}v^*q_n||_2^2=\tau(q_nvv_{\lambda}v^*q_nvv_{\lambda}^*v^*q_n)=\tag 7.i$$ $$\tau(v_{\lambda}(v^*q_nv)v_{\lambda}^*(v^*q_nv))=\tau(\beta_{\lambda}(v^*q_nv)\hskip 0.03in(v^*q_nv)).$$
Combining (7.g) and (7.i) yields $$\sum_{\gamma\in F_n^{\lambda}}\tau(\sigma_{\gamma}(q_n)q_n)=\tau(\beta_{\lambda}(v^*q_nv)\hskip 0.03in(v^*q_nv)),\forall\lambda\in\Lambda\tag 7.j$$

Recall that $\Lambda\cong \Gamma'/K$ and let $\pi:\Gamma'\rightarrow\Lambda$ be the projection given by $\pi(\delta(\lambda)k)=\lambda$, for all $\lambda\in\Lambda$ and $k\in K$. Let $\rho:\Gamma'\rightarrow$ Aut$(A\overline{\otimes}B)$ be the action defined by $\rho_g=\sigma_g\otimes \beta_{\pi(g)}$, for all $g\in\Gamma'$. Since $\sigma_{|\Gamma'}$ and $\beta$ are both weakly mixing, we get that $\rho$ is weakly mixing.
Let $\tau_B$ denote the normalized trace on $B$, i.e. $\tau_B(p)=1$ and $\tau_B(b)=\tau(p)^{-1}\tau(b)$, for all $b\in B$.
Since $\rho$ is weakly mixing, we can find  a sequence $\{g_m\}_{\geq 1}\subset\Gamma'$ such that $\sigma_{g_m}(a)\rightarrow\tau(a)$ and $\beta_{\pi(g_m)}(b)\rightarrow\tau_B(b)$, weakly, for all $a\in A$ and $b\in B$, as $m\rightarrow\infty$. We may clearly assume that $g_m=\delta(\lambda_m)k$, for some sequence $\{\lambda_m\}_{m\geq 1}\subset\Lambda$ and $k\in K$. Thus, we have that $\sigma_{\delta(\lambda_m)}(a)\rightarrow\tau(a)$ and $\beta_{\lambda_m}(b)\rightarrow\tau_B(b)$, weakly, for all $a\in A$, $b\in B$.

In particular, it follows that we can find $\lambda\in\Lambda$ such that $$\tau(\sigma_{\delta(\lambda)g}(q_n)q_n)\leq (1+\varepsilon)\tau(q_n)^2, \forall g\in F\hskip 0.03in\text{and}\tag 7.k$$ 
$$\tau(\beta_{\lambda}(v^*q_nv)\hskip 0.03in(v^*q_nv))=\tau(p)\tau_B(\beta_{\lambda}(v^*q_nv)\hskip 0.03in(v^*q_nv))\geq$$ $$ (1-\varepsilon)\tau(p)\tau_B(v^*q_nv)^2=(1-\varepsilon)\tau(p)^{-1}\tau(q_n)^2.$$
Since $F_n^{\lambda}\subset\delta(\lambda)F$, (7.j) and (7.k) together imply that $(1+\varepsilon)|F_n^{\lambda}|\hskip 0.03in|\tau(q_n)|^2\geq (1-\varepsilon)\tau(p)^{-1}\tau(q_n)^2.$ Hence, $|F_n^{\lambda}|\geq (1-\varepsilon)^2\tau(p)^{-1}$.
From this, (7.g) and (7.h) we deduce that the cardinality of the set $S_n:=\{\gamma\in\Gamma|\hskip 0.03in |c_{\gamma}|\geq (1-2^{-n})\tau(p)\}$ is at least $(1-\varepsilon)^2\tau(p)^{-1}$. As $\varepsilon>0$ is arbitrary, we get that $|S_n|\geq\tau(p)^{-1}$, for all $n\geq 1$. Thus $S=\cap_{n\geq 1}S_n=\{\gamma|\in\Gamma|\hskip 0.03in |c_{\gamma}|\geq\tau(p)\}$ has cardinality at least $\tau(p)^{-1}$.
On the other hand, since $\sum_{\gamma\in\Gamma}|c_{\gamma}|^2=||p||_2^2=\tau(p)$, we get that $|S|\leq\tau(p)^{-1}$. The last two facts clearly give the claim.
\hfill$\square$

\vskip 0.1in
Going back to the proof of (2), notice that by combining the above claim and [Po06b, Lemma 5.3] we can find a finite group $K'\subset \Gamma$ and a  character $\chi$ of $K'$ such that $p=|K'|^{-1}\sum_{k'\in K'}\chi(k')u_{k'}$. Thus $K=\{\gamma\in\Gamma|u_{\gamma}p\in\Bbb Cp\}=K'$ and
the last part of the proof of [Po06b, Theorem 5.2.] (pages 439--441) gives the conclusion.		
\hfill$\blacksquare$

\vskip 0.2in
\head \S 8. {Strong rigidity for embeddings of II$_1$ factors.}\endhead
\vskip 0.2in
We are now ready to prove the main technical results of this paper. Let $M=A\rtimes_{\sigma}\Gamma$ and $N=D\rtimes_{\rho}\Lambda$ be the crossed product algebras associated with actions of countable groups $\Gamma$ and $\Lambda$ on abelian von Neumann algebras $A$ and $D$.
 We say that an embedding $\Delta:N\rightarrow M$ is {\it standard} if we can find a subgroup $\Gamma_0\subset\Gamma$ and a $\sigma(\Gamma_0)$--invariant subalgebra $A_0\subset A$ such that,  up to conjugation with a unitary element and modulo scalars, we have that $\Delta(\Lambda)=\Gamma_0$ and $\Delta(D)=A_0$. In other words, a standard embedding arises from a ``realization" of $\Lambda$ as a subgroup $\Gamma_0\subset\Gamma$ and of $\rho$ as a quotient of $\sigma_{|\Gamma_0}$.  

 Assume that $\Lambda$ is a ``rigid group" (e.g. $\Lambda$ has property (T) or is a product of two non--amenable groups) and that $\sigma$ is a Bernoulli action.
Our first result roughly says that any embedding $\Delta:N\rightarrow M$ is  standard unless it comes from an embedding of $N$ into $L(\Gamma)$.
As we will see in Section 10, in certain situations (e.g. if $\Gamma=\Bbb F_2\times\Bbb F_2)$, we can rule out the appearance of such {\it bad} embeddings and thus we can completely describe the embeddings of $N$ into $M$. This is why we use the phrase {\it strong rigidity for embeddings} to refer to the results of this section.
\vskip 0.05in

\proclaim {8.1 Theorem} Let $\Gamma$ be a torsion free, ICC, countable group, $B$ be an abelian von Neumann algebra and  denote $M=B^{\Gamma}\rtimes\Gamma$.
\vskip 0.02in
\noindent
Let $\rho:\Lambda\rightarrow$ Aut$(D)$ be a free action of a countable group $\Lambda$ on an abelian von Neumann algebra $D$ and denote $N=D\rtimes_{\rho}\Lambda$. 
Suppose that $\Lambda$ admits an infinite, almost normal subgroup $\Lambda_0$ such that either
\vskip 0.02in
\noindent
(1) the inclusion $(\Lambda_0\subset\Lambda)$ has the relative property (T) or

\noindent
(2) $\Lambda_0$ is generated by two commuting subgroups $\Lambda_1,\Lambda_2$, with $\Lambda_1$ non--amenable and $\Lambda_2$ infinite.
\vskip 0.02in
\noindent
 Let $\Delta:N\rightarrow M$ be a  $*$--homomorphism  and suppose that $\Delta(N)\nprec_{M}L(\Gamma)$.

\vskip 0.03in
\noindent
 Then $\Delta$ must be unital and there exist  a character $\eta$ of $\Lambda$, a group homomorphism $\delta:\Lambda\rightarrow\Gamma$ and a unitary $u\in M$ such that  $\Delta(D)\subset uAu^*$ and $\Delta(v_{\lambda})=\eta(\lambda)uu_{\delta(\lambda)}u^*$, for all $\lambda\in\Lambda$.
\endproclaim
Recall that a subgroup $\Lambda_0$ of a countable group $\Lambda$ is {\it almost normal} if $\lambda\Lambda_0\lambda^{-1}\cap\Lambda_0$ has finite index in both $\lambda\Lambda_0\lambda^{-1}$ and $\Lambda_0$, for any $\lambda\in\Lambda$.

\vskip 0.05in
We are also able to prove a strong rigidity result for the embeddings of $N$ into a tensor product $M=M_1\overline{\otimes}M_2$ of two II$_1$ factors coming from Bernoulli actions. If we represent $M$ as a crossed product $A\rtimes\Gamma$, as in 3.1, then  our second result gives a list of assumptions  which force an embedding $\Delta:N\rightarrow M$ to be standard, in the sense from above. 
This result will be crucial in deriving $W^*$--superrigidity (see Section 9).

\proclaim {8.2 Theorem} In the context from 3.1, i.e.  $M_i=B_i^{\Gamma_i}\rtimes\Gamma_i$, $M=M_1\overline{\otimes}M_2$, $A=B_1^{\Gamma_1}\overline{\otimes}B_2^{\Gamma_2}$, $\Gamma=\Gamma_1\times\Gamma_2$,  assume that $\Gamma$ is ICC (i.e. $\Gamma_1,\Gamma_2$ are ICC).
\vskip 0.02in
\noindent
Let $\rho:\Lambda\rightarrow$ Aut$(D)$ be a free action of a countable group $\Lambda$ on an abelian von Neumann algebra $D$ and denote $N=D\rtimes_{\rho}\Lambda$. Suppose that $\Lambda$ admits an infinite, almost normal subgroup $\Lambda_0$ such that the inclusion $(\Lambda_0\subset\Lambda)$ has the relative property (T).
\vskip 0.03in
\noindent
 Let $\Delta:N\rightarrow M$ be a  $*$--homomorphism and suppose that:

\noindent
(1) $\Delta(L(\Lambda_0))\nprec_{M}L(\Gamma_1)\otimes 1$ and $\Delta(L(\Lambda_0))\nprec_{M}1\otimes L(\Gamma_2)$.

\noindent
(2) $\Delta(N)\nprec_{M}L(\Gamma_1)\overline{\otimes}M_2$ and $\Delta(N)\nprec_{M}M_1\overline{\otimes}L(\Gamma_2)$.
\vskip 0.05in
\noindent
We have that
\vskip 0.02in
\noindent
$\bullet$ If $\Gamma$ is torsion free, then  $\Delta$ must be unital  and we can find  a character $\eta$ of $\Lambda$, a group homomorphism 
$\delta:\Lambda\rightarrow\Gamma$ and a unitary $u\in M$ such that  $\Delta(D)\subset uAu^*$ and $\Delta(v_{\lambda})=\eta(\lambda)uu_{\delta(\lambda)}u^*$, for all $\lambda\in\Lambda$.
\vskip 0.07in
\noindent
$\bullet$ In general,  we can find a finite index subgroup $\Lambda_1$ of $\Lambda$, a 1--1 map $\delta:\Lambda_1\rightarrow\Gamma$,  $u_1,..,u_n\in\Cal N_{M}(A)$ and $u\in\Cal U(M)$ such that if we let $\tilde\Delta=\text{Ad}(u)\circ\Delta:N\rightarrow M$ then

\vskip 0.03in
\noindent
(a)   
$\tilde\Delta(D)\subset A$, there exists a projection $0\not=p\in A$ and $a_i^{\lambda}\in (A)_1$, for all $i\in\{1,..,n\},\lambda\in\Lambda_1$, such that $p\leq\tilde\Delta(1)$ and $p\tilde\Delta(v_{\lambda})=\sum_{i=1}^na_i^{\lambda}u_{\delta(\lambda)}u_i$, for all $\lambda\in\Lambda_1$,
\vskip 0.03in
\noindent
(b) there exists a finite subgroup $K\subset\Gamma$ such that $\delta(\lambda)$ normalizes $K$ and we have that $\delta(\lambda\lambda')^{-1}\delta(\lambda)\delta(\lambda')$ belongs to $K$, for all $\lambda,\lambda'\in\Lambda_1$. 
\vskip 0.03in
\noindent
If in addition $\max\hskip 0.01in\{|K|\hskip 0.03in| K\hskip 0.03in\text{finite subgroup of}\hskip 0.03in\Gamma\}<\infty$, then  $\tau(\Delta(1))\in l^{-1}\Bbb N,$ where $l$ is the least common multiple of $|K|$ with $K$ being a finite subgroup of $\Gamma$.
\endproclaim

The conclusion of Theorem 8.2 (in the case when $\Gamma$ has torsion) is rather involved and technical and  might seem hard to work with.
However, it is  manageable enough to be useful in applications.  Thus, we will use it to prove $W^*$--superrigidity of Bernoulli actions of property (T) groups without assuming torsion freeness. This is of particular interest when applied to the linear groups SL$_n(\Bbb Z)$, for $n\geq 3$.

Theorem 8.1 can be also extended (just as Theorem 8.2) to cover groups with torsion.
 On the other hand, note that Theorem 8.1 applies to a larger class of groups $\Lambda$ than  Theorem 8.2. Let us explain why this difference appears.
The first step in the proofs of both 8.1 and 8.2 consists of showing that if $\Delta:N\rightarrow M$ is an embedding, then $L(\Lambda)$ is ``absorbed" by $L(\Gamma)$. When $\Lambda$ has property (T) this is granted by Theorem 2.1, in both the contexts from 8.1 and 8.2. If $\Lambda=\Lambda_1\times\Lambda_2$ is a product of two non--amenable groups, the absorption result still holds in the context of 8.1 (by [Po08]). 
 However, if $M=M_1\overline{\otimes}M_2$ is as in 8.2, then the absorption result fails in general. Indeed, one can easily construct a situation in which we have that  $L(\Lambda_1)\subset M_1$ and $L(\Lambda_2)\subset M_2$  without having that $L(\Lambda)$ is absorbed by $L(\Gamma)$.

We first  prove Theorem 8.2. Since the proof of Theorem 8.1 is analogous to the proof of Theorem 8.2, we will only sketch it, leaving most of the details to the reader.
\vskip 0.1in
\noindent
{\it Proof of Theorem 8.2}.
Let $\Delta:N\rightarrow M$ be a $*$--homomorphism satisfying (1)--(2) and set $q=\Delta(1)$. 
From now on, we identify $N$ with $\Delta(N)$ and forget about $\Delta$.   Since the inclusion $(\Lambda_0\subset\Lambda)$  has the  relative property (T), the inclusion $(L(\Lambda_0)\subset qMq)$ is rigid ([Po06c]).  By using (1) and the fact that $\Gamma$ is ICC, Theorem 2.1 implies that we can find a unitary $u\in M$ such that
$u(q\Cal N_{qMq}(L(\Lambda_0))'')u^*\subset L(\Gamma)$. Since $\Lambda_0$ is almost normal in $\Lambda$ we  deduce that $$uL(\Lambda)u^*\subset L(\Gamma)\tag 8.a$$ 

Next, we claim that $$D\nprec_{M}L(\Gamma_1)\overline{\otimes}M_2\hskip 0.05in\text{and}\hskip 0.05in D\nprec_{M}M_1\overline{\otimes}L(\Gamma_2) \tag 8.b$$ 

 Assume by contradiction that $D\prec_{M}L(\Gamma_1)\overline{\otimes}M_2$. Since $D$ is abelian, we can find a non--zero projection $r\in D'\cap qMq$ and $a,b\in M$ such that $(D)_1r\subset a(L(\Gamma_1)\overline{\otimes}M_2)_1b$. Recall that $M_1=B_1^{\Gamma_1}\rtimes\Gamma_1$ and let $\tilde M_1\supset M_1$, $\{\theta_t^1\}_{t\in \Bbb R}\subset$ Aut$(\tilde M_1)$ be the weakly malleable deformations defined in section 1.5. Set $\theta_t=\theta_t^1\otimes$ id$_{M_2}\in$ Aut$(\tilde M_1\overline{\otimes}M_2)$.

Let $\Cal V$ be a $\rho(\Lambda)$--invariant, $||.||_2$--dense subgroup of $\Cal U(D)$. Then $\Cal U=\{v_{\lambda}v|\lambda\in\Lambda,v\in\Cal V\}$ is a group of unitaries of $N$ which generates $N$ as a von Neumann algebra. By using (8.a) and the inclusion $\Cal Vr\subset a(M_2)_1b$, we get that $\Cal Ur\subset u^*(L(\Gamma))_1ua(L(\Gamma_1)\overline{\otimes} M_2)_1b.$ Since ${\theta_t}_{|L(\Gamma_1)\overline{\otimes}M_2}=$ id$_{L(\Gamma_1)\overline{\otimes}M_2}$, we can find $t>0$ such that $||\theta_t(x)r-xr||_2\leq \frac{||r||_2}{2},$ for all $x\in\Cal U$. Thus, for all $x\in\Cal U$ we have that  $\Re\hskip 0.02in\langle(\theta_t(x)r-xr),xr\rangle\geq -\frac{||r||_2^2}{2}$, hence $$\Re \hskip 0.02in\tau(x^*\theta_t(x)r)=\Re\hskip 0.02in\langle \theta_t(x)r,xr\rangle=\langle xr,xr\rangle + \Re\hskip 0.02in\langle(\theta_t(x)r-xr),xr\rangle\geq\frac{||r||_2^2}{2}.$$

Since $\Cal U$ is a group, by averaging (see e.g. the proof of 2.1) we deduce that there exists $0\not=v\in qM$ such that $v\theta_t(x)=xv$, for all $x\in\Cal U$. Since $\Cal U''=N$, we get that  $v\theta_t(x)=xv$, for all $x\in N.$
The proof of Theorem 2.1 (see Remark 2.2) implies that $N\prec_{M}L(\Gamma_1)\overline{\otimes}M_2.$ This contradicts assumption (2), hence, (8.b) is proven.

\vskip 0.05in

Further, by combining (8.a), (1) and (8.b),  Theorem 6.1 implies that $C=D'\cap qMq$ is a type I von Neumann algebra and that we may assume that $q=\Delta(1)\in A$ and  $Aq\subset C$.

Denote by $\Cal Z$  the center of $C$. Since $D$ is abelian,  $\Cal Z$ contains $D$. Since $C$ is of type I we can decompose $\Cal Z=\oplus_{i\geq 1}\Cal Z_i$, where $\Cal Z_i$ are abelian von Neumann algebras, such that $C=\oplus_{i\geq 1}(\Bbb M_{n_i}(\Bbb C)\otimes\Cal Z_i)$, for some (possibly finite) strictly increasing sequence $1\leq n_1< n_2<..$. Since $A\subset M$ is maximal abelian, we get that $Aq$ is maximal abelian in $C$. Since, by a classical result, any two maximal abelian subalgebras of a type I algebra are conjugate (see e.g. [Va07, Lemma C.2.]), we may assume that $Aq=\oplus_{i\geq 1}(\Bbb C^{n_i}\otimes\Cal Z_i)$, where $\Bbb C^{n_i}\subset\Bbb M_{n_i}(\Bbb C)$ is the subalgebra of diagonal matrices.    

Now, since $\Lambda=\{v_{\lambda}\}_{\lambda\in\Lambda}\subset\Cal U(qMq)$ normalizes $D$, it also normalizes  $C=D'\cap qMq$ and $\Cal Z=\Cal Z(C)$. Moreover, $\Lambda$ normalizes  $\Bbb M_{n_i}(\Bbb C)\otimes\Cal Z_i$ and $\Cal Z_i$ , for all $i\geq 1$. Denote by $\alpha$ the action of $\Lambda$ on $\Cal Z$ given by $\alpha_{\lambda}(z)=v_{\lambda}zv_{\lambda}^*$, for all $\lambda\in\Lambda$ and $z\in\Cal Z$. Since $\alpha$ leaves $\Cal Z_i$ invariant, for all $i\geq 1$, we can define an action $\beta$ of $\Lambda$ on $C$ by letting $$\beta_{\lambda}=\oplus_{i\geq 1}(\text{id}_{\Bbb M_{n_i}(\Bbb C)}\otimes{\alpha_{\lambda}}_{|\Cal Z_i})\in\hskip 0.02in\text{Aut}(C),\hskip 0.02in\forall\lambda\in\Lambda.$$

Let $\lambda\in\Lambda$.
Since the automorphisms $\beta_{\lambda}$ and Ad$(v_{\lambda})$ of $C$ are equal on its center, $\Cal Z$, by [KR97, Corollary 9.3.5.] we can find a unitary ${\omega}_{\lambda}\in C$ such that $\beta_{\lambda}=$ Ad$(v_{\lambda}{\omega}_{\lambda})$, as automorphisms of $C$. As $\beta$ leaves $Aq=\oplus_{i\geq 1}(\Bbb C^{n_i}\otimes\Cal Z_i)$ invariant, we get that $v_{\lambda}{\omega}_{\lambda}\in\Cal N_{qMq}(Aq)$.  Thus, we can find $z_{\lambda}\in\Cal U(Aq)$ and $V_{\lambda}\in\Cal N_{qMq}(Aq)$ of the form $V_{\lambda}=\sum_{\gamma\in\Gamma}p_{\gamma,\lambda}u_{\gamma}$, where $p_{\gamma,\lambda}\in\Cal P(Aq)$, such that $v_{\lambda}{\omega}_{\lambda}=V_{\lambda}z_{\lambda}$. 

We claim that $V_{\lambda}V_{\lambda'}=V_{\lambda\lambda'}$, for all $\lambda,\lambda'\in\Lambda$. 
To see this, let $v=V_{\lambda\lambda'}^*V_{\lambda}V_{\lambda'}\in\Cal U(qMq)$. Note that for all $x\in Aq$ we have that $V_{\lambda}xV_{\lambda}^*=V_{\lambda}z_{\lambda}xz_{\lambda}^*V_{\lambda}^*=\beta_{\lambda}(x)$. Using this fact,  we immediately get that $v\in (Aq)'\cap qMq=Aq$. On the other hand, $v$ is of the form $v=\sum_{\gamma\in\Gamma}p_{\gamma}u_{\gamma}$, for some $p_{\gamma}\in \Cal P(Aq)$ (this is because $V_{\lambda},V_{\lambda'}$ and $V_{\lambda\lambda'}^*$ are of this form). By combining the last two observations we get that $v=1$, as claimed.

We also claim that $E_A(V_{\lambda})=0$, for all $\lambda\in\Lambda\setminus\{e\}$. Otherwise, we can find $\lambda\not=e$ and a non--zero projection $p\in Aq$ such that $\beta_{\lambda}(x)=x$, for all $x\in Ap$. Thus, we can find a non--zero projection $p_0\in \Cal Z$ such that $\alpha_{\lambda}(x)=x$, for all $x\in \Cal Zp_0$. Since $D\subset\Cal Z$, we get that $\alpha_{\lambda}(x)p_0=xp_0$, for all $x\in D$. Let $0\not=p_1\in D$ be the support projection of $E_D(p_0)$. By projecting onto $D$ and noticing that ${\alpha_{\lambda}}_{|D}=\rho_{\lambda}$, as automorphisms of $D$, we get that $\rho_{\lambda}(x)p_1=xp_1$, for all $x\in D$. This, however, contradicts the freeness of $\rho$.

Altogether, we proved that there exist $w_{\lambda}={\omega}_{\lambda}z_{\lambda}^*\in\Cal U(C)$ such that $V_{\lambda}=v_{\lambda}w_{\lambda}$ normalizes $C$, $Aq$ and $\Cal Z$,  $V_{\lambda}V_{\lambda'}=V_{\lambda\lambda'}$ and $E_A(V_{\lambda})=\delta_{\lambda,e}q$, for all $\lambda,\lambda'\in\Lambda$. For $i\geq 1$, let $e_{i,1},..,e_{i,n_i}\in\Bbb C^{n_i}$ be the canonical projections. Then $\beta$ fixes $p_{i,j}:=e_{i,j}\otimes 1\in\Bbb C^{n_i}\otimes \Cal Z_i$.

Next, we prove the following:

\vskip 0.1in
\noindent
{\it Claim 1.} There exists a unitary $v\in\Cal U(M)$ such that $vV_{\lambda}v^*\in L(\Gamma),$ for all $\lambda\in\Lambda$.
\vskip 0.05in
\noindent
{\it Proof of claim 1.} Let $Q=\{V_{\lambda}|\lambda\in\Lambda_0\}''$ and $P=\{V_{\lambda}|\lambda\in\Lambda\}''$. Since the inclusion $(\Lambda_0\subset\Lambda)$ has the relative property (T), by [Po06c] the inclusion $(Q\subset M)$ is rigid. Since $\Lambda_0$ is almost normal in $\Lambda$, we get that  $P\subset q\Cal N_{qMq}(Q)''$. Since $\Gamma$ is ICC,  by Theorem 2.1, in order to prove the claim, it suffices to show that $Q\nprec_{M}L(\Gamma_1)$ and $Q\nprec_{M}L(\Gamma_2)$.

Let us show that $Q\nprec_{M}L(\Gamma_1)$. Since $L(\Lambda_0)\nprec_{M}L(\Gamma_1)$, Theorem 1.3.1 gives a sequence $\{\lambda_n\}_{n\geq 1}\subset\Lambda_0$ satisfying  $||E_{L(\Gamma_1)}(xv_{\lambda_n}y)||_2\rightarrow 0$, for all $x,y\in M$. We claim that $$\sup_{a\in (A)_1}||E_{L(\Gamma_1)}(xv_{\lambda_n}yaz)||_2\rightarrow 0,\hskip 0.04in\forall x,y,z\in M\tag 8.c$$ 

Let $a\in (A)_1$.
Since $uL(\Lambda_0)u^*\subset L(\Gamma)$, we can assume that $v_{\lambda_n}\in L(\Gamma)$, for the purpose of proving (8.c). Then we can decompose $v_{\lambda_n}=\sum_{\gamma\in\Gamma}c_{n,\gamma}u_{\gamma}$, where $c_{n,\gamma}\in\Bbb C$, for all $n\geq 1$ and $\gamma\in\Gamma$. It is clearly enough to prove (8.c) when $x=x_0u_g,y=y_0u_h, z=z_0u_k$, for some $x_0,y_0,z_0\in (A)_1$ and $g,h,k\in\Gamma$. In this case, we have that $$xv_{\lambda_n}yaz=\sum_{\gamma\in\Gamma}c_{n,\gamma}x_0\sigma_{g\gamma}(y_0)\sigma_{g\gamma h}(az_0)u_{g\gamma hk}.$$

\noindent
Since $a,x_0,y_0,z_0\in (A)_1$, we get that $|\tau(x_0\sigma_{g\gamma}(y_0)\sigma_{g\gamma h}(az_0))|\leq 1$, for all $\gamma\in\Gamma$.  Thus $$||E_{L(\Gamma_1)}(xv_{\lambda_n}yaz)||_2^2=\sum_{\gamma\in g^{-1}\Gamma_1k^{-1}h^{-1}}|c_{n,\gamma}|^2|\tau(x_0\sigma_{g\gamma}(y_0)\sigma_{g\gamma h}(az_0))|^2\leq$$ $$\sum_{\gamma\in g^{-1}\Gamma_1k^{-1}h^{-1}}|c_{n,\gamma}|^2=||E_{L(\Gamma_1)}(u_gv_{\lambda_n}u_{hk})||_2^2.$$ As we have that $||E_{L(\Gamma_1)}(u_gv_{\lambda_n}u_{hk})||_2\rightarrow 0$, (8.c) is therefore proven.

Next, fix $\varepsilon>0$ and $x,y\in M$. Let $q_0\in C$ be a projection such that $||q_0-q||_2\leq\varepsilon$ and $Cq_0=\oplus_{i=1}^k(\Bbb M_{n_i}(\Bbb C)\otimes\Cal Z_i)$, for some $k\geq 1$.  Remark that we can find $z_1,..,z_l\in C$ such that $(Cq_0)_1\subset\sum_{j=1}^l(A)_1z_j$. Moreover, for every $u\in\Cal U(C)$ we have that $||u(q-q_0)||_2=||q-q_0||_2\leq \varepsilon$. By using these facts and recalling that $w_{\lambda_n}\in\Cal U(C)$, we deduce that  $$||E_{L(\Gamma_1)}(xV_{\lambda_n}y)||_2=||E_{L(\Gamma_1)}(xv_{\lambda_n}w_{\lambda_n}y)||_2\leq $$ $$||E_{L(\Gamma_1)}(xv_{\lambda_n}w_{\lambda_n}(q-q_0)y)||_2+||E_{L(\Gamma_1)}(xv_{\lambda_n}(w_{\lambda_n}q_0)y)||_2\leq $$ $$\varepsilon||x||\hskip 0.02in||y||+\sum_{j=1}^l\sup_{a\in (A)_1}||E_{L(\Gamma_1)}(xv_{\lambda_n}az_jy)||_2\hskip 0.04in\forall n\geq 1.$$ This inequality in combination with (8.c) gives that $\limsup_{n\rightarrow\infty}||E_{L(\Gamma_1)}(xV_{\lambda_n}y)||_2\leq\varepsilon$. Since $\varepsilon>0$ is arbitrary, we deduce that $\lim_{n\rightarrow\infty}||E_{L(\Gamma_1)}(xV_{\lambda_n}y)||_2=0$, for all $x,y\in M$. Since $\{V_{\lambda_n}\}_{n\geq 1}\subset\Cal U(Q)$, we get  that $Q\nprec_{M}L(\Gamma_1)$.\hfill$\square$
\vskip 0.1in
Next, recall that $V_{\lambda}$ normalizes $Aq$ and $V_{\lambda}xV_{\lambda}^*=\beta_{\lambda}(x)$, for all $x\in Aq$, $\lambda\in\Lambda$.
\vskip 0.1in
\noindent
{\it Claim 2} ([Po06b]).  There exist a finite index subgroup  $\Lambda_1\subset\Lambda$  and a non--zero $\beta(\Lambda_1)$--invariant projection $p\in Aq\cap v^*L(\Gamma)v$ such that the restriction of the action  $\beta_{|\Lambda_1}$ to $Ap$  is weakly mixing. Moreover, we have that $p\leq p_{i,j}$, for some $i\geq 1$ and $j\in\{1,..,n_i\}$.
\vskip 0.05in
\noindent
{\it Proof of claim 2.} This claim follows from an argument of Popa (see the proofs of [Po06b, Lemma 4.5.] and [Va07, Theorem 9.1.]).  Let $A_0=Aq\cap v^*L(\Gamma)v$ and note that $q\in A_0$. It is easy to see that $A_0$ is a completely atomic, (abelian) von Neumann algebra. Since $V_{\lambda}\in v^*L(\Gamma)v$, for all $\lambda\in\Lambda$, we get that $V_{\lambda}$ normalizes $A_0$, hence $A_0$ is $\beta$--(globally) invariant. Let $p\in A_0$ be a minimal projection.  Then $p$ is invariant under some finite index subgroup $\Lambda_1$ of $\Lambda$. Denote $Q_1:=\{V_{\lambda}|\lambda\in\Lambda_1\}''\subset qMq$ and notice that $p\in Q_1'\cap qMq$. 

Towards  proving that the restriction of $\beta_{|\Lambda_1}$ to $Ap$  is weakly mixing, let $\Cal H\subset Ap$ be a  finite dimensional $\beta(\Lambda_1)$--invariant subspace. Since for every $\lambda\in\Lambda_1$ and $\xi\in\Cal H$, we have that $V_{\lambda}\xi=\beta_{\lambda}(\xi)V_{\lambda}\subset\Cal HV_{\lambda}$ and $\xi V_{\lambda}=V_{\lambda}\beta_{\lambda^{-1}}(\xi)\subset V_{\lambda}\Cal H$, we derive that $\Cal H$ is contained in the quasi--normalizer of $Q_1p$ in $pMp$.  
 Since $Q_1p\subset v^*L(\Gamma)v$, while $Q_1\nprec_{M}L(\Gamma_1)$ and $Q_1\nprec_{M}L(\Gamma_2)$ (by the proof of claim 1),  [Va08, Lemma 4.2.] gives that $\Cal H\subset v^*L(\Gamma)v$. Thus, we get that $\Cal H\subset Ap\cap v^*L(\Gamma)v=A_0p$. Since $p$ is a minimal projection in $A_0$, we must have that $\Cal H=\Bbb Cp$, which proves the weak mixingness assertion.
The moreover assertion is clear since  $p_{i,j}$ is $\beta$--invariant,  for all $i\geq 1,j\in\{1,..,n_i\}$. 
\hfill$\square$
\vskip 0.05in

 To summarize, we have that 

$\bullet$ $\{V_{\lambda}\}_{\lambda\in\Lambda}\subset v^*L(\Gamma)v\cap \Cal U(qMq)$ normalizes $Aq$ 

$\bullet$ $p\in Aq\cap v^*L(\Gamma)v$ commutes with $V_{\lambda}$, for all $\lambda\in\Lambda_1$, and

$\bullet$ the action $\Lambda_1\ni\lambda\rightarrow {\beta_{\lambda}}_{|Ap}=$ Ad$(V_{\lambda}p)\in$ Aut$(Ap)$ is weakly mixing. 
\vskip 0.05in
Let $p'=vpv^*\in L(\Gamma)$ and $\Lambda_1=\{v(V_{\lambda}p)v^*\}_{\lambda\in\Lambda_1}\subset  \Cal U(p'L(\Gamma)p')$. Then $\Lambda_1$ acts on  $B=v(Ap)v^*\subset p'(A\rtimes_{\sigma}\Gamma)p'$ by conjugation and the resulting action (denoted $\beta'$) is weakly mixing.
Moreover, we have that $\{B,\Lambda_1\}''\cong B\rtimes_{\beta'}\Lambda_1$. Indeed, this is because $E_B(v(V_{\lambda}p)v^*)=E_{Ap}(V_{\lambda}p)=E_A(V_{\lambda})p=0$, for all $\lambda\in\Lambda_1\setminus\{e\}$.

By applying Theorem 7.1 we can find maps $h:\Lambda_1\rightarrow S^1$, $\delta:\Lambda_1\rightarrow\Gamma$ and a unitary $w\in L(\Gamma)$ such that $[wu_{\delta(\lambda)}w^*,p']=0$ and
 $$v(V_{\lambda}p)v^*=h(\lambda)(wu_{\delta(\lambda)}w^*)p',\forall \lambda\in\Lambda_1\tag 8.d$$

Let $K=\{\gamma\in\Gamma|(wu_{\gamma}w^*)p'\in\Bbb Cp'\}$ and $\Gamma'=\{\delta(\lambda)k|\lambda\in\Lambda_1,k\in K\}$. Then $K,\Gamma'$ are subgroups of $\Gamma$, $|K|<\infty$ and $\Gamma'$ normalizes $K$. It is also clear that $\delta(\lambda\lambda')^{-1}\delta(\lambda)\delta(\lambda')\in K$, for all $\lambda,\lambda'\in\Lambda_1$, thus proving condition (b) in the general case.
To show that $\sigma_{|\Gamma'}$ is weakly mixing, let $\pi_i:\Gamma\rightarrow \Gamma_i$ be the projection $\pi_i(\gamma_1,\gamma_2)=\gamma_i$. Then $\sigma$ is the diagonal product of the generalized Bernoulli actions $\sigma_i\circ\pi_i:\Gamma\rightarrow$ Aut$(B_i^{\Gamma_i})$. Thus, if $\sigma_{|\Gamma'}$ is not weakly mixing then either $(\sigma_1\circ\pi_1)_{|\Gamma'}$ or $(\sigma_2\circ\pi_2)_{|\Gamma'}$ is not weakly mixing. By [Po07a, Lemma 4.5.] we would get that some finite index subgroup $\Gamma''$ of $\Gamma'$ is contained into either $\Gamma_1$ or $\Gamma_2$.
 By using (8.d), this would imply that a corner of $Q_1=\{V_{\lambda}|\lambda\in\Lambda_1\}''$ embeds into either $L(\Gamma_1)$ or $L(\Gamma_2)$, a contradiction, by the proof of claim 1.

\vskip 0.1in
\noindent
For the rest of the proof we analyze separately the case when $\Gamma$ is torsion free and the general case:
\vskip 0.05in
\noindent
$\bullet$ If $\Gamma$ is {\it torsion free} then $K=\{e\}$, $h$ is a character, $\delta$ is a 1--1 homomorphism and $\Gamma'=\delta(\Lambda_1)$. Since $\sigma_{|\delta(\Lambda_1)}$ is ergodic, by Theorem 7.1 we	get that $p=1$ and $B=wAw^*$.

Since $p\leq q=\Delta(1)$, $\Delta$ must be unital.
Also, since by claim 2 we have that $p\leq p_{i,j}$, for some $i\geq 1$ and $j\in\{1,..,n_i\}$, we deduce that $C$ is abelian and $\Cal Z=A=C$. Moreover, we get that $\beta_{|\Lambda_1}$ is weakly mixing and hence $\beta$ itself is weakly mixing. Thus, in claim 2 and after we may take  $\Lambda_1=\Lambda$. If $z=v^*w$, then $z$ normalizes $A$ and (8.d) can be rewritten as $$V_{\lambda}=h(\lambda)zu_{\delta(\lambda)}z^*,\forall\lambda\in\Lambda\tag 8.e$$

Let $\theta=$ Ad$(z)\in$ Aut$(A)$. Then (8.e) implies that $\beta_{\lambda}=\theta\circ\sigma_{\delta(\lambda)}\circ\theta^{-1},$ for all $\lambda\in\Lambda$. In other words, the actions $\beta$ and $\sigma_{|\delta(\Lambda)}$ are conjugate.

Next, recall that $v_{\lambda}=V_{\lambda}w_{\lambda}^*,$ for all $\lambda\in\Lambda$, for some unitary $w_{\lambda}\in C=A$. Since $v_{\lambda}v_{\lambda'}=v_{\lambda\lambda'}$ and $v_{\lambda}=\beta_{\lambda}(w_{\lambda}^*)V_{\lambda}$, for all $\lambda,\lambda'\in\Lambda$, we get that the map $\Lambda\ni\lambda\rightarrow \beta_{\lambda}(w_{\lambda}^*)\in\Cal U(A)$ is a 1--cocycle for $\beta$.  By the proof of claim 2, $Q=\{V_{\lambda}|\lambda\in\Lambda_0\}''$ satisfies $Q\nprec_{M}L(\Gamma_1)$ and $Q\nprec_{M}L(\Gamma_2)$. By reasoning as above we deduce that $\sigma_{|\delta(\Lambda_0)}$ is weakly mixing. Note that the inclusion $(\delta(\Lambda_0)\subset\delta(\Lambda))$ has the relative property (T) and is almost normal. Finally, note that $\sigma$ is {\it s--malleable} (see [Po07a, section 4] for definition and proof).

Altogether, we can apply Popa' cocycle superrigidity theorem ([Po07a, Corollary 5.2.]) and deduce that the action $\sigma_{|\delta(\Lambda)}$ is $\Cal U_{fin}$--cocycle superrigid (in the sense of [Po07a, 5.6.0.]). Thus, $\beta_{|\Lambda}$ is cocycle superrigid, so we can find $U\in\Cal U(A)$ and a character $h'$ of $\Lambda$ such that $\beta(\lambda)(w_{\lambda}^*)V_{\lambda}=h'(\lambda)UV_{\lambda}U^*,$ for all $\lambda\in\Lambda$. Together with (8.e) this provides a character $\eta=hh'$ of $\Lambda$ and  $W=Uz\in\Cal N_{M}(A)$  such that $v_{\lambda}=\eta(\lambda)Wu_{\delta(\lambda)}W^*,$ for all $\lambda\in\Lambda$. As $A=C\supset D$, this gives the conclusion in the case when $\Gamma$ is torsion free.
\vskip 0.05in

\noindent
$\bullet$ For {\it general} $\Gamma$, since $\sigma_{|\Gamma'}$ is weakly mixing, by Theorem 7.1 we deduce that $p_0=w^*p'w$ is of the form $p_0=|K|^{-1}\sum_{k\in K}\chi(k)u_k\in L(K)$, for some character $\chi$ of $K$, and that $B=w(A^{K}p_0)w^*$, where $A^K=\{a\in A|\sigma(k)(a)=a,\forall k\in K\}$.

 In particular, this shows that $\tau(p)=|K|^{-1}\in l^{-1}\Bbb N$, for every minimal projection $p$ of  $Aq\cap vL(\Gamma)v^*$  and so we get that $\tau(\Delta(1))=\tau(q)\in l^{-1}\Bbb N$ (if we assume that $l<\infty$).

Now, let $z=v^*w$. Then (8.d) implies that $V_{\lambda}p=h(\lambda)z(u_{\delta(\lambda)}p_0)z^*$, for all $\lambda\in\Lambda_1$. Also, since $B=v(Ap)v^*=w(A^{K}p_0)w^*$, we get that $$Ap=z(A^{K}p_0)z^*\tag 8.f$$
Next, let $q_0\in A$ be a projection such that the projections $\{\sigma(k)(q_0)\}_{k\in K}$ form a partition of unity in $A$ (here, we use the fact that $K$ acts freely). Then it is easy to verify that $\xi=|K|^{\frac{1}{2}}p_0q_0$ is a partial isometry such that $\xi^*\xi=q_0,\xi\xi^*=p_0$ and $A^Kp_0=\xi(Aq_0)\xi^*$. Thus,  (8.f) gives that $z\xi(Aq_0)(z\xi)^*=Ap$ and so we can find $\eta\in\Cal N_{M}(A)$ extending the partial isometry $z\xi$. Since $\eta\xi^*=zp_0$ we get that $$ V_{\lambda}p=h(\lambda)\eta(\xi^*u_{\delta(\lambda)}\xi)\eta^*,\forall\lambda\in\Lambda_1\tag 8.g$$

Recall that $v_{\lambda}=V_{\lambda}w_{\lambda}^*$, where $w_{\lambda}\in\Cal U(C)$, and that  $C=\oplus_{i\geq 1}(\Bbb M_{n_i}(\Bbb C)\otimes\Cal Z_i)\supset Aq=\oplus_{i\geq 1}(\Bbb C^{n_i}\otimes\Cal Z_i)$. Let $i$ such that $p\leq p_{i,j}$, for some $j\in\{1,..,n_i\}$.  Since $\Bbb C^{n_i}\otimes\Cal Z_i$ is regular in $\Bbb M_{n_i}(\Bbb C)\otimes\Cal Z_i$ and $\Bbb C^{n_i}\otimes\Cal Z_i=Aq_i$, for some $q_i\in\Cal P(Aq)$,    we can find
 $u_1,..,u_m\in\Cal N_{M}(A)$ such that $(\Bbb M_{n_i}(\Bbb C)\otimes\Cal Z_i)_1\subset\sum_{i=1}^{m}(A)_1u_i$.  Since $pC\subset\Bbb M_{n_i}(\Bbb C)\otimes\Cal Z_i$, we derive that $pw_{\lambda}^*\in \sum_{i=1}^{m}(A)_1u_i$, for all $\lambda\in\Lambda_1$. Combining this with (8.g), the definition of $\xi$ and the fact that $\eta$ normalizes $A$ gives that $$pv_{\lambda}=(V_{\lambda}p)(pw_{\lambda}^*)\in
\eta(\xi^*u_{\delta(\lambda)}\xi)\eta^*\sum_{i=1}^m(A)_1u_i\subset \tag 8.h$$ $$|K|\sum_{k,k'\in K}\sum_{i=1}^m (A)_1 \eta u_ku_{\delta(\lambda)}u_{k'}\eta^*u_i,\forall\lambda\in\Lambda_1. $$ 
Recall that $N\subset qMq$ (unitally), $D\subset Aq$ and $p\in Aq$.  Finally, let $\tilde\Delta=$ Ad$(\eta^*)_{|N}:N\rightarrow M$.
Since $\delta(\lambda)$ normalizes $K$ and  $\eta\in\Cal N_{M}(A)$, (8.h) implies that $$(\eta^*p\eta)\tilde\Delta(v_{\lambda})\in |K|\sum_{k,k'\in K}\sum_{i=1}^m(A)_1u_{\delta(\lambda)}u_{kk'}\eta^*u_i\eta,\forall\lambda\in\Lambda_1.$$
Since $\tilde\Delta(D)\subset A$, $\eta^*p\eta\leq \tilde\Delta(1_N)=\eta^*q\eta$ and $u_k$ ($k\in K$), $\eta$, $u_i$ all normalize $A$, this completes the proof. 
\hfill$\blacksquare$
\vskip 0.1in

\noindent {\it Proof of Theorem 8.1}.  Let $\Delta:N\rightarrow qMq$ be a unital $*$--homomorphism, for some $q\in\Cal P(M)$, and identify $N$ with $\Delta(N)$. Assume that $N\nprec_{M}L(\Gamma)$. Since $D$ is regular in $N$ and $\sigma$ is mixing, by [Po06a, Theorem 3.1]  we get that $D\nprec_{M}L(\Gamma).$ 

Further, we claim that we can find  $u\in\Cal U(M)$ such that $uL(\Lambda)u^*\subset L(\Gamma).$
In case (a), this follows from Theorem 2.1, while in case (b) it follows from [Po08, Lemma 5.2] and [CI08, Theorem 3]. By combining there two facts and Theorem 6.2 and by reasoning as in the proof of 8.2 (the torsion free case), the conclusion follows.
\hfill$\blacksquare$

\vskip 0.2in
\head \S 9.{ W$^*$--superrigidity}.\endhead
\vskip 0.2in
In this section we prove Theorem A. Let us first restate it using the language of von Neumann algebras.

\proclaim {9.1 Theorem} Let $\Gamma$ be an ICC countable group which admits an infinite, almost normal subgroup $\Gamma_0$ such that the inclusion $(\Gamma_0\subset\Gamma)$ has the  relative property (T). Let $B$ be a non--trivial abelian von Neumann algebra.  Denote $M=B^{\Gamma}\rtimes\Gamma$ and $A=B^{\Gamma}$. 
\vskip 0.01in
\noindent
Let $\rho:\Lambda\rightarrow$ Aut$(C)$ be a free ergodic action of a countable group $\Lambda$ on an abelian von Neumann algebra $C$. Denote $N=C\rtimes\Lambda$. 

\vskip 0.01in
\noindent
Let $\theta:N\rightarrow qMq$ be a $*$--isomorphism, for some projection $q\in M$. 

\vskip 0.05in
\noindent
Then $q=1$ and there exist a unitary $u\in M$, a character $\eta$ of $\Lambda$, a group isomorphism $\delta:\Lambda\rightarrow\Gamma$ 
such that $\theta(C)=uAu^*$ and $\theta(v_{\lambda})=\eta(\lambda)uu_{\delta(\lambda)}u^*$, for all $\lambda\in\Lambda$.

\endproclaim

Define the $*$--homomorphism $\Delta:N\rightarrow N\overline{\otimes}N$ by letting $\Delta(cv_{\lambda})=cv_{\lambda}\otimes v_{\lambda},$ for all $c\in C$ and $\lambda\in\Lambda$.
Note that $\Delta$ has been introduced in the proof of [PV09, Lemma 3.2.].

Since $N\cong qMq$, we can view $\Delta$ as a non--unital embedding of $M$ into $M\overline{\otimes}M$. In the proof of Theorem 9.1 we will apply Theorem 8.2 to $\Delta$. We will then need certain properties of $\Delta$, that we  record in the following lemma.

\proclaim {9.2 Lemma}  Let $Q\subset N$ be a von Neumann subalgebra. We have that:
\vskip 0.01in
\noindent
(1) If $\Delta(Q)\prec_{N\overline{\otimes}N}N\otimes 1$, then $Q\prec_{N}C$.
\vskip 0.01in
\noindent
(2) If $\Delta(Q)\prec_{N\overline{\otimes}N}1\otimes N$, then $Q$ is not diffuse.
\vskip 0.01in
\noindent
(3) If $\Delta(N)\prec_{N\overline{\otimes}N}Q\overline{\otimes}N$, then $C\prec_{N}Q$.
\vskip 0.01in
\noindent
(4) If $\Delta(N)\prec_{N\overline{\otimes}N}N\overline{\otimes}Q$, then $L(\Lambda)\prec_{N}Q$. Moreover, if $\Delta(N)s\prec_{N\overline{\otimes}N}N\overline{\otimes}Q$, for every non--zero projection $s\in\Delta(N)'\cap N\overline{\otimes}N$, then $L(\Lambda)r\prec_{N}Q$, for every non--zero projection $r\in L(\Lambda)'\cap N$.
\endproclaim
\noindent
{\it Proof of Lemma 9.2}. (1) Arguing by contradiction, if $Q\nprec_{N}C$, we can find a sequence $\{u_n\}_{n\geq 1}\subset\Cal U(Q)$ such that $||E_C(au_nb)||_2\rightarrow 0$, for all $a,b\in N$. We claim that $$||E_{N\otimes 1}(x\Delta(u_n)y)||_2\rightarrow 0,\hskip 0.04in\text{for all}\hskip 0.04in x,y\in N\overline{\otimes}N\tag 9.a$$ 

Since $E_{N\otimes 1}$ is $N\otimes 1$--bimodular, we may assume that $x=1\otimes au_g$, $y=1\otimes bu_h$, for some $a,b\in C$ and $g,h\in\Lambda$.
Then for  all $n\geq 1$, we have that $$x\Delta(u_n)y=\sum_{\lambda\in\Lambda}E_{C}(u_nv_{\lambda}^*)v_{\lambda}\otimes av_gv_{\lambda}bv_h$$
Thus we get that $||E_{N\otimes 1}(x\Delta(u_n)y)||_2=|\tau(a\rho_{h^{-1}}(b)|\hskip 0.02in ||E_C(u_nv_{\lambda}^*)||_2\rightarrow 0$, which proves (9.a). Since $u_n\in\Cal U(Q)$, (9.a) implies that $\Delta(Q)\nprec_{N\overline{\otimes}N} N\otimes 1$, and we are done.

\vskip 0.05in
\noindent
(2) The proof of this part is similar to the proof of (1), so we leave it as an exercise.
\vskip 0.05in
\noindent
(3) If $\Delta(N)\prec_{N\overline{\otimes}N}Q\overline{\otimes}N$, then  $\Delta(C)\prec_{N\overline{\otimes}N}Q\overline{\otimes}N$.  The conclusion follows easily after noticing that $\Delta(C)=C\otimes 1$.
\vskip 0.05in
\noindent
(4) Arguing by contradiction, assume that $L(\Lambda)r\nprec_{N}Q$, for some non--zero projection $r\in L(\Lambda)'\cap N$. Since the group $\{v_{\lambda}r\}_{\lambda\in\Lambda}$ generates $L(\Lambda)r$ as a von Neumann algebra, we can find a sequence $\{\lambda_n\}_{n\geq 1}\subset\Lambda$ such that $||E_Q(av_{\lambda_n}rb)||_2\rightarrow 0$, for all $a,b\in N$. 

Define $s:=1\otimes r\in\Delta(N)'\cap N\overline{\otimes}N$. We claim that $$||E_{N\overline{\otimes}Q}(x\Delta(v_{\lambda_n})sy)||_2\rightarrow 0,\hskip 0.04in\forall x,y\in N\overline{\otimes}N\tag 9.b$$

Since $E_{N\overline{\otimes}Q}$ is $N\otimes 1$--bimodular, we may take $x=1\otimes a$, $y=1\otimes b$, for some $a,b\in N$. Then $E_{N\overline{\otimes}Q}(x\Delta(v_{\lambda_n}s)y)=v_{\lambda_n}\otimes E_Q(av_{\lambda_n}rb).$
Since $||E_Q(av_{\lambda_n}rb)||_2\rightarrow 0$, this proves (9.b).  Now, since $\Delta(v_{\lambda_n})s\in\Cal U(\Delta(N)p)$, (9.b) implies that $\Delta(N)s\nprec_{N\overline{\otimes}N}N\overline{\otimes}Q$ and our proof by contradiction is over. 
\hfill$\blacksquare$

\vskip 0.1in
\noindent {\it Remark.}  
The moreover part of (4) was first noticed by Stefaan Vaes who pointed out to us that it can be used to simplify our initial proof of Theorem 9.1. Specifically, the moreover part of (4) enables us to avoid repeating the proof of Theorem 6.1 in the proof of {\bf Case (5)} below.

\vskip 0.1in
\noindent
{\it Proof of Theorem 9.1}. Let $\theta:N\rightarrow qMq$ be an isomorphism, for some $q\in\Cal P(M)$. 
Let $t=\tau(q)^{-1}$ and $k\geq t$ be an integer. View $\Bbb C^{k}$ as the algebra of diagonal matrices in $\Bbb M_k(\Bbb C)$ and let $r\in \Bbb C^k\otimes C$ be a projection of normalized trace $\frac{t}{k}$ in $\Bbb M_k(\Bbb C)\otimes C$. Consider the $*$--homomorphism id$_{\Bbb M_k(\Bbb C)}\otimes\Delta:\Bbb M_k(\Bbb C)\otimes N\rightarrow \Bbb M_k(\Bbb C)\otimes (N\overline{\otimes}N)\cong (\Bbb M_k(\Bbb C)\otimes N)\overline{\otimes}N$. Since (id$_{\Bbb M_k(\Bbb C)}\otimes \Delta)(r)=r\otimes 1$ and $r(\Bbb M_k(\Bbb C)\otimes N)r\cong N^t$, we get a unital $*$--homomorphism $\Delta_t:N^t\rightarrow N^t\overline{\otimes}N$.

Fix an embedding $N\subset N^t$ and view $\Delta_t$ as $*$--homomorphism $\Delta_t:N^t\rightarrow N^t\overline{\otimes}N^t$. It is easy to see that Lemma 9.2 holds true if we replace $N$ and $\Delta$ by $N^t$ and $\Delta_t$, throughout.
We denote by $C^t=r(\Bbb C^k\otimes C)r\subset N^t$. Remark that $C^t\subset N^t$ is a Cartan subalgebra and that $\Delta_t(x)=x\otimes q$, for all $x\in C^t$. For simplicity, we denote $1\otimes q$ by $q$. 
From now on, we identify   $N^t\cong M$ (via $\theta^t$, the $t$--amplification of $\theta$)
 and view $\Delta_t$ as $*$--homomorphism $\Delta_t:M\rightarrow M\overline{\otimes}M$ with $\Delta_t(1)=q$.

\vskip 0.05in
To continue, assume for the moment that $(\diamond)$ there is a non--zero projection $s\in\Delta_t(M)'\cap q(M\overline{\otimes}M)q$ such that $\Delta_t(M)s\nprec_{M\overline{\otimes}M}M\overline{\otimes}L(\Lambda)$. Let $\tilde\Delta_t:M\rightarrow M\overline{\otimes}M$ be given by $\tilde\Delta_t(x)=\Delta_t(x)s$.
Since $\tilde\Delta_t(M)\nprec_{M\overline{\otimes}M}M\overline{\otimes}L(\Lambda)$, 
by applying Theorem 8.2 to $\tilde\Delta_t$ we conclude that we are in one of the following cases:

\vskip 0.05in
\noindent
 {\bf (1)} There exist a finite index subgroup $\Gamma_1$ of $\Gamma$, a 1--1 map $\delta:\Gamma_1\rightarrow\Gamma\times\Gamma$, $u_1,..,u_n\in\Cal N_{M\overline{\otimes}M}(A\overline{\otimes}A)$ and  $u\in\Cal U(M\overline{\otimes}M)$ such that  $\tilde\Delta:=$ Ad$(u)\circ\tilde\Delta_t:M\rightarrow M\overline{\otimes}M$ satisfies

\noindent
 ({\it a})  $\tilde\Delta(A)\subset A\overline{\otimes}A$, there exists  $0\not=p\in \Cal P(A\overline{\otimes}A)$ and $a_i^{\gamma}\in (A\overline{\otimes}A)_1$, for all $i\in\{1,..,n\}$, $\gamma\in\Gamma_1$  such that $p\leq \tilde\Delta(1)$ and $p\tilde\Delta(u_{\gamma})=\sum_{i=1}^na_i^{\gamma}u_{\delta(\gamma)}u_i$, for all $\gamma\in\Gamma_1$, and

\noindent
({\it b}) there is a finite subgroup $K$ of $\Gamma\times\Gamma$ such that $\delta(\gamma\gamma')^{-1}\delta(\gamma)\delta(\gamma')\in K$ and $\delta(\gamma)$ normalizes $K$, for all $\gamma,\gamma'\in\Gamma_1$.
\vskip 0.05in
\noindent
{\bf (2)} $\tilde\Delta_t(L(\Gamma_0))\prec_{M\overline{\otimes}M}L(\Gamma)\otimes 1$.
\vskip 0.05in
\noindent
{\bf (3)} $\tilde\Delta_t(L(\Gamma_0))\prec_{M\overline{\otimes}M}1\otimes L(\Gamma)$.
\vskip 0.05in
\noindent
{\bf (4)} $\tilde\Delta_t(M)\prec_{M\overline{\otimes}M}L(\Gamma)\overline{\otimes}M$.
\vskip 0.05in

On the other hand, if $(\diamond)$ fails then we are in the following case:
\vskip 0.05in
\noindent
{\bf (5)} $\Delta_t(M)s\prec_{M\overline{\otimes}M}M\overline{\otimes}L(\Gamma),$ for any non--zero projection $s\in\Delta_t(M)'\cap q(M\overline{\otimes}M)q$.

\vskip 0.1in

\noindent
For the rest of the proof, we analyze each one of these cases and show that they either lead to a contradiction or imply that ($\star$) $C^t\prec_{M}A$ or that ($\star\star$) $uL(\Lambda)u^*\subset L(\Gamma)$, for some $u\in \Cal U(M)$. If ($\star$) holds, since $C^t$ and $A$ are Cartan subalgebras of $N^t\cong M$,  [Po06c, Theorem A.1.] provides a unitary $v\in M$ such that $A=vC^tv^*$.
 Popa's cocycle/OE superrigidity theorems ([Po07a, 5.2. and 5.6.]) then imply that $t=1$ (thus $q=1$) and that the isomorphism $\theta:N\rightarrow 
M$ is of the desired form. If ($\star\star$) holds, then the conclusion follows from [Po06b, Theorem 0.7.] (see Theorem 6.3).

\vskip 0.05in
\noindent
{\bf Case (1)}. Let $G$ be the set of all $g\in\Gamma_1$ for which we can find $k\in K$ such that $\delta(g)k\in \Gamma\times\{e\}$. Since $\delta$ satisfies ({\it b}),   $G$ is a subgroup of $\Gamma_1$. 

Let us prove that $G$ is finite. If $g\in G$, then $\delta(g)\in (\Gamma\times\{e\})K$. Condition ({\it a}) implies that $p\tilde\Delta(u_{g})\in\sum_{i=1}^n\sum_{k\in K}(M\overline{\otimes}A)_1u_ku_i$.
This shows in particular that $\tilde\Delta(L(G))\prec_{M\overline{\otimes}M}M\overline{\otimes}A$. On the other hand, we may assume that {\bf (2)} and {\bf (3)} are false since these cases will be dealt with later. Since the inclusion $\tilde\Delta(L(\Gamma_0))\subset\tilde\Delta(L(\Gamma))$ is  rigid and quasi--regular, by  Theorem 2.1  we can find a unitary $v\in M\overline{\otimes}M$ such that $\tilde\Delta(L(G))\subset\tilde\Delta(L(\Gamma))\subset v(L(\Gamma)\overline{\otimes}L(\Gamma))v^*$. The last two facts readily imply that $\tilde\Delta(L(G))\prec_{M\overline{\otimes}M}L(\Gamma)\otimes 1$. 

By the construction of $\tilde\Delta$ it follows that $\Delta_t(L(G))\prec_{M\overline{\otimes}M}L(\Gamma)\otimes 1$.
 But then Lemma 9.2 (1) (which holds for $\Delta_t$) yields that $L(G)\prec_{M}C^t$. 

 Further, by [Va08, Lemma 3.5.] we deduce that $C^t\prec_{M} L(G)'\cap M$. If $G$ is infinite, then $L(G)$ is diffuse. Hence, since $\sigma$ is mixing, [Po06a, Theorem 3.1.] implies that $L(G)'\cap M\subset L(\Gamma)$. Thus, we would get that $C^t\prec_{M}L(\Gamma)$. Since $C^t$ is regular in $M$, by applying [Po06a, Theorem 3.1.] (see also Lemma 3.4)  we derive a contradiction. This shows that $G$ is finite.

\vskip 0.05in
\noindent
Using the fact that $G$ is finite we next deduce the following claim:
\vskip 0.1in
\noindent{\it Claim.} Let $\{x_m\}_{m\geq 1}\subset M$ be a sequence with $||x_m||\leq 1$. Denote $x_m^{\gamma}=E_A(x_mu_{\gamma}^*)$ and assume that $||x_m^{\gamma}||_2\rightarrow 0$, for all $\gamma\in\Gamma$. Then we have that $$||E_{M\otimes 1}(yp\tilde\Delta(x_m)z)||_2\rightarrow 0,\hskip 0.06in\text{for all}\hskip 0.04in y,z\in M\overline{\otimes}M\tag 9.c$$
\vskip 0.1in
\noindent{\it Proof of the claim.} Let $M_1=A\rtimes\Gamma_1$. Since $\Gamma_1$ has finite index in $\Gamma$, we may assume that $x_m\in M_1$, for all $m\geq 1$.

Since $E_{M\otimes 1}$ is $M\otimes 1$--bimodular,  we may assume that $y=1\otimes au_g$ and $z=1\otimes bu_h$, for $a,b\in A$ and $g,h\in\Gamma$. Then $E_{M\otimes 1}(yxz)=E_{M\otimes 1}((1\otimes a\sigma_g(b))x(1\otimes u_{gh})$ for all  $x\in M\overline{\otimes}M$.  Thus, to prove (9.c), we may assume that $y\in A\overline{\otimes}A$ and $z=u_g$, for some $g\in\Gamma\times\Gamma$.

\noindent
Since $||E_{M\otimes 1}(yp\tilde\Delta(x_m)u_g)||_2\leq ||E_{M\overline{\otimes} A}(yp\tilde\Delta(x_m)u_g)||_2\leq ||y||\hskip 0.02in ||E_{M\overline{\otimes}A}(p\tilde\Delta(x_m)u_g)||_2$ (here we use that $y\in M\overline{\otimes}A$), it is sufficient to prove that 
$$||E_{M\overline{\otimes}A}(p\tilde\Delta(x_m)u_g)||_2\rightarrow 0,\hskip 0.06in\text{for all}\hskip 0.04in g\in\Gamma\times\Gamma\tag 9.d$$

\noindent
Now, ({\it a}) gives that $p\tilde\Delta(x_m)=\sum_{\gamma\in\Gamma_1}\tilde\Delta(x_m^{\gamma})p\tilde\Delta(u_{\gamma})=\sum_{i=1}^n\sum_{\gamma\in\Gamma_1}\tilde\Delta(x_m^{\gamma})a_i^{\gamma}u_{\delta(\gamma)}u_i$. As $\tilde\Delta(x_m^{\gamma})a_i^{\gamma}\in A\overline{\otimes}A$, this identity implies that in order to get (9.d) it suffices to show that $$||\sum_{\gamma\in\Gamma_1}\tilde\Delta(x_m^{\gamma})a_i^{\gamma}E_{M\overline{\otimes}A}(u_{\delta(\gamma)}u_iu_g)||_2\rightarrow 0,\hskip 0.06in\text{for all}\hskip 0.04in i\in\{1,..,n\}\tag 9.e$$

To prove (9.e), fix $i\in\{1,..,n\}$ and denote $v=u_iu_g$. Since $u_{\delta(\gamma)}v$ normalizes $A\overline{\otimes}A$ we can find $p_{\gamma}\in\Cal P(A\overline{\otimes}A)$ such that $E_{M\overline{\otimes}A}(u_{\delta(\gamma)}v)=p_{\gamma}u_{\delta(\gamma)}v$. Since $\delta$ is 1-1, the $A\overline{\otimes}A$--bimodules $\{(A\overline{\otimes}A)u_{\delta(\gamma)}v\}_{\gamma\in\Gamma_1}$ are mutually orthogonal. 
Altogether, (9.e) rewrites as $$\sum_{\gamma\in\Gamma_1}\tau(|\tilde\Delta(x_m^{\gamma})|^2|a_i^{\gamma}|^2p_{\gamma})\rightarrow 0,\hskip 0.06in \text{as}\hskip 0.04in m\rightarrow\infty \tag 9.f$$

Next, we claim that $\sum_{\gamma\in\Gamma_1}\tau(p_{\gamma})<\infty.$ Identify $A=L^{\infty}(X)$, where $(X,\mu)$ is a probability space. Then $\sigma$ induces a free ergodic p.m.p. action $\Gamma\curvearrowright X$. On $(X^2,\mu^2):=(X\times X,\mu\times\mu)$ we consider the direct product action of $\Gamma\times\Gamma$: $(\gamma_1,\gamma_2)\circ (x_1,x_2)=(\gamma_1x_1,\gamma_2x_2)$. Let $\phi$ be the automorphism of $X^2$  given by $vav^*=a\circ\phi^{-1}$, for all $a\in L^{\infty}(X^2)$. Then $p_{\gamma}$ is precisely the characteristic function of the set of $x\in X^2$ satisfying
$(\delta(\gamma)\circ\phi)^{-1}(x)\in (\Gamma\times\{e\})x$.

Hence, if we let $Y_{\gamma}=\{x\in X^2|\phi^{-1}(x)\in (\Gamma\times\{e\})\delta(\gamma)x\}$, then $\tau(p_{\gamma})=\mu^2(Y_{\gamma})$. If $\gamma,\gamma'\in\Gamma$ satisfy $\mu^2(Y_{\gamma}\cap Y_{\gamma'})>0$, then the freeness of the action $\Gamma\times\Gamma\curvearrowright X^2$ implies that $(\Gamma\times\{e\})\delta(\gamma)\cap (\Gamma\times\{e\})\delta(\gamma')\not=\emptyset$. Thus $\delta(\gamma)\delta(\gamma')^{-1}\in \Gamma\times\{e\}$, hence $\gamma{\gamma'}^{-1}\in G$. This shows that if $G\delta(\gamma)\cap G\delta(\gamma')=\emptyset$, then $\mu^2(Y_{\gamma}\cap Y_{\gamma'})=0$. Since $\delta$ is 1--1, it follows that $\sum_{\gamma\in\Gamma}\mu^2(Y_{\gamma})\leq |G|$, hence $\sum_{\gamma\in\Gamma}\tau(p_{\gamma})<\infty$, as claimed.

To prove (9.f) and the claim, fix $\varepsilon>0$ and let $F\subset\Gamma_1$ be finite such that $\sum_{\gamma\in \Gamma_1\setminus F}\tau(p_{\gamma})<\varepsilon$.
Since $||x_m^{\gamma}||,||a_i^{\gamma}||,||p_{\gamma}||\leq 1$ and $\tau(|\tilde\Delta(x_m^{\gamma})|^2)=||\tilde\Delta(x_m^{\gamma})||_2^2\leq ||x_m^{\gamma}||_2^2$, we get that $$\sum_{\gamma\in\Gamma_1}\tau(|\tilde\Delta(x_m^{\gamma})|^2|a_i^{\gamma}|^2p_{\gamma})\leq \varepsilon+\sum_{\gamma\in F}||x_m^{\gamma}||_2^2.$$

Since $||x_m^{\gamma}||_2\rightarrow 0$ and $\varepsilon>0$ is arbitrary, this indeed proves (9.f). \hfill$\square$
\vskip 0.05in
Recall that $\tilde\Delta(x)=u\Delta_t(x)su^*$, for all $x\in M$, and that $\Delta_t(x)=(x\otimes 1)q$, for all $x\in C^t$.
Since $p\leq \tilde\Delta(1)=usu^*$, it follows that $u^*p\tilde\Delta(x)u=u^*pu(x\otimes 1)$ and hence that $||u^*p\tilde\Delta(x)u||_2=||p||_2$, for all $x\in\Cal U(C^t)$. Then the above claim implies that we cannot find a sequence $\{x_m\}_{m\geq 1}\subset\Cal U(C^t)$ such that $||E_A(x_mu_{\gamma}^*)||_2\rightarrow 0$, for all $\gamma\in\Gamma$.
Thus $C^t\prec_{M}A$, which settles Case {\bf (1)}.
\vskip 0.05in
\noindent
{\bf Case (2)}. If $\tilde\Delta_t(L(\Gamma_0))\prec_{M\overline{\otimes}M}L(\Gamma)\otimes 1$, then $\Delta_t(L(\Gamma_0))\prec_{M\overline{\otimes}M}L(\Gamma)\otimes 1$ and Lemma 9.2 (1) implies that $L(\Gamma_0)\prec_{M}C^t$. Since $\Gamma_0$ is infinite, we derive a contradiction by arguing as the proof of {\bf (1)}.
\vskip 0.05in
\noindent
{\bf Case (3)}. If $\tilde\Delta_t(L(\Gamma_0))\prec_{M\overline{\otimes}M}1\otimes L(\Gamma)$, then  $\Delta_t(L(\Gamma_0))\prec_{M\overline{\otimes}M}1\otimes L(\Gamma)$ and Lemma 9.2 (2) implies that $L(\Gamma_0)$ cannot be diffuse. This contradicts the fact that $\Gamma_0$ is infinite.
\vskip 0.05in
\noindent
{\bf Case (4)}. If  $\tilde\Delta_t(M)\prec_{M\overline{\otimes}M}L(\Gamma)\overline{\otimes}M$, then $\Delta_t(M)\prec_{M\overline{\otimes}M}L(\Gamma)\overline{\otimes}M$ and Lemma 9.1 (3) entails $C^t\prec_{M}L(\Gamma)$. Since $\sigma$ is mixing and $C^t$ is regular in $M$, [Po06a, Theorem 3.1.] leads to a contradiction.
\vskip 0.05in
\noindent
{\bf Case (5)}. Assume that $\Delta_t(M)s\prec_{M\overline{\otimes}M}M\overline{\otimes}L(\Gamma),$ for any non--zero projection $s\in\Delta_t(M)'\cap q(M\overline{\otimes}M)q$. Then the moreover part of Lemma 9.2 (4) gives that $L(\Lambda)r\prec_{M}L(\Gamma)$, for any non--zero projection $r\in L(\Lambda)'\cap N$. Since $\Gamma$
is ICC, by reasoning as in the end of the proof of Theorem 2.1 we can find a unitary $u\in M$ such that $uL(\Lambda)u^*\subset\Gamma$.
 This finishes the proof of the last case of Theorem 9.1.\hfill$\blacksquare$

\vskip 0.2in

\head \S 10. {Further applications}\endhead

\vskip 0.2in
In this section, we derive several applications of the results of Section 8. 
\vskip 0.1in
\noindent
(I) {\bf Group von Neumann algebra decomposition}. First, we use Theorem 8.2 to provide a large class of II$_1$ factors which cannot be decomposed as the group von Neumann algebra, $L(G)$,  of some countable group $G$. 

\proclaim  {10.1 Corollary} Let $\Gamma$ be a countable ICC group which admits an infinite, almost normal subgroup $\Gamma_0$ such that the inclusion $(\Gamma_0\subset\Gamma)$ has the relative property (T). 
\vskip 0.01in
\noindent
Assume that $\sup\hskip 0.02in\{|K|\hskip 0.02in |K$ is a finite subgroup of $\Gamma\}$ is finite. Let $l=l(\Gamma)$ be the least common multiple of $|K|$, with $K$ ranging over all finite subgroups of $\Gamma$. 
\vskip 0.01in
\noindent
 Let $B$ be a non--trivial abelian von Neumann algebra. Set $M=B^{\Gamma}\rtimes_{\sigma}\Gamma$ and let $p\in\Cal P(M)$.  

\vskip 0.05in
\noindent
 If $\tau(p)\not\in \{\frac{1}{l},\frac{2}{l},..,\frac{l}{l}\}$, then the {\rm II}$_1$ factor $pMp$ is not isomorphic to a group von Neumann algebra. 
In particular, if $\Gamma$ is torsion free (i.e. $l=1$), then $pMp$ is not isomorphic to a group von Neumann algebra, for any $p\not=1$.
\endproclaim

\noindent
{\it Proof.} Assume that $pMp\cong N:=L(G)$, for some projection $p\in M$ and a countable (necessarily)  ICC group  $G$. In other words, we have that $M\cong N^t$, where $t=\tau(p)^{-1}$.

Next, let $\Delta:N\rightarrow N\overline{\otimes} N$ be the $*$--homomorphism given by $\Delta(u_g)=u_g\otimes u_g$, for all $g\in G$. By amplifying $\Delta$ (as in the proof of Theorem 9.1), we get a  $*$--homomorphism $\Delta_t:N^t\rightarrow N^t\overline{\otimes}N^t$ such  that $\tau(\Delta_t(1))=\tau(p)$. We continue with:

\proclaim {10.2 Lemma} Let $Q\subset N^t$ be a not necessarily unital von Neumann subalgebra. Then we have the following:
\vskip 0.03in
\noindent
(1) If $\Delta_t(Q)\prec_{N^t\overline{\otimes}N^t}N^t\otimes 1$ or $\Delta_t(Q)\prec_{N^t\overline{\otimes}N^t}1\otimes N^t$, then $Q$ cannot be diffuse.

\noindent
(2) If $\Delta_t(N^t)\prec_{N^t\overline{\otimes}N^t}Q\overline{\otimes}N^t$ or $\Delta_t(N^t)\prec_{N^t\overline{\otimes}N^t}N^t\overline{\otimes}Q$, then $N^t\prec_{N^t}Q$.
\endproclaim
\noindent
 The proof of this lemma is analogous to the proof of Lemma 9.1 and we therefore leave it as an exercise.
\vskip 0.05in

Now, since $M\cong N^t$, we can view $\Delta_t$ as a $*$--homomorphism $\Delta_t:M\rightarrow M\overline{\otimes}M$. Thus, we can apply Theorem 8.2 to $\Delta_t$.  Lemma 10.2 guarantees that conditions (1)--(2) from 8.2 are satisfied. Finally,  by the last assertion of Theorem 8.2 we deduce that $\tau(p)=\tau(\Delta_t(1))\in l^{-1}\Bbb N$, as claimed.
\hfill$\blacksquare$
\vskip 0.1in

\vskip 0.05in
\noindent
{\bf 10.3 Remarks.} (1) 
Recall that $A\wr\Gamma=(\oplus_{\Gamma}A)\rtimes\Gamma$, where $\Gamma$ acts on $\oplus_{\Gamma}A$ by shifting the copies of $A$, is the {\it wreath product} of two groups $A$ and $\Gamma$.
Corollary 10.1 implies that if $n\geq 3$ and  $t\in (0,1)\setminus\Bbb Q$, then $L(\Bbb Z\wr\text{SL}_n(\Bbb Z))^t$ is not a group von Neumann algebra. Indeed, SL$_n(\Bbb Z)$ is ICC and has property (T) by Kazhdan's classical result ([Ka67]). Additionally, if $m\geq 3$ then the kernel of the natural quotient SL$_n(\Bbb Z)\rightarrow$ SL$_n(\Bbb Z/m\Bbb Z)$ is a torsion free group. This shows that $|K|\leq |$SL$_n(\Bbb Z/m\Bbb Z)|$, for every finite subgroup $K$ of SL$_n(\Bbb Z)$.
\vskip 0.03in

\noindent
(2) The first examples of II$_1$ factors which are not group von Neumann algebras were exhibited by Connes by means of his $\chi(M)$ invariant ([C75]). Recently,  Popa's deformation/rigidity theory has been used to give new examples of such factors ([IPP08], [PV08a]). In all of these cases, one moreover proves that the II$_1$ factors involved do not have anti--automorphisms, i.e. $M\not\cong M^{op}$. 

\noindent
(3) If $M$ is as in Corollary 10.1, then $pMp$ admits an involutary anti--automorphism, for every $p\in\Cal P(M)$. Indeed, the formula $\Phi(au_{\gamma})=u_{\gamma^{-1}}a$, for $a\in B^{\Gamma}$ and $\gamma\in\Gamma$, defines an anti--automorphism of $M$. If we take $p\in B^{\Gamma}$, then $\Phi(p)=p$ and  $\Phi_p=\Phi_{|pMp}$ is an involutary anti--automorphism of $pMp$.

The first examples of II$_1$ factors which are not group von Neumann algebras and yet have anti--automorphisms have been constructed in [Jo80]. We point out that, as opposed to our examples, the examples constructed in [Jo08] have no involutary anti--automorphism.

\vskip 0.1in
\noindent
(II) {\bf Symmetry groups of II$_1$ factors.} Given a II$_1$ factor $M$, there are three natural objects capturing the symmetries of $M$:
 $\Cal F(M)$, the fundamental group of $M$,  Out$(M)$, the outer automorphism group of $M$ and Bimod$(M)$, the set of bi--finite Hilbert $M$--bimodules $\Cal H$, i.e. such that  dim$(_{M}\Cal H),$ dim$(\Cal H_{M})<\infty$.

The calculation of these invariants is an extremely challenging problem and it is only recently that explicit calculations were obtained for large families of II$_1$ factors by using Popa's deformation/rigidity theory. While there is, by now, an extensive literature on this problem (see e.g. the introduction of [Va08]), we only mention here that the existence of II$_1$ factors $M$ for which $\Cal F(M)$, Out$(M)$ and Bimod$(M)$ are trivial has been proven in [Po06], [IPP08] and [Va09], respectively.

In this subsection, we consider three natural invariants for II$_1$ factors which generalize the ones from above.
\vskip 0.03in
\noindent
{\bf 10.4 Definition.} Let $M$ be a II$_1$ factor. Then we define 
\vskip 0.03in
\noindent
$\bullet$ $\Cal F_s(M)$, {\it the fundamental semigroup of $M$}, to be the set of $t>0$ for which there exists a unital 
$*$--homomorphism $\theta:M\rightarrow M^t$.
\vskip 0.03in
\noindent
$\bullet$ End$(M)$, {\it the endomorphism semigroup of $M$}, to be the set of unital 

\noindent
$*$--homomorphisms 
$\theta:M\rightarrow M$, and
\vskip 0.03in
\noindent
$\bullet$ LFBimod$(M)$, to be the set of {\it left--finite} Hilbert $M$--bimodules $\Cal H$, i.e. such that dim$(_{M}\Cal H)<\infty$.
\vskip 0.1in

It is clear that $\Cal F_s(M)$ is a semigroup (with respect to multiplication)  which contains $\Cal F(M)$.  Thus, whenever $\Cal F(M)=(0,+\infty)$ (e.g. if $M$ is the hyperfinite II$_1$ factor), we get that $\Cal F_s(M)=(0,+\infty)$. Also, we have that $\Bbb N\subset \Cal F_s(M)$.  Furthermore, we have a dichotomy: either $\Cal F_s(M)\cap (0,1)=\emptyset$ or $\Cal F_s(M)=(0,+\infty)$. This is a consequence of the following two facts: $(\star)$ if $\{t_n\}_{n\geq 1}\in\Cal F_s(M)$ is a sequence such that $t:=\sum_{n\geq 1}t_n<\infty$, then $t\in\Cal F_s(M)$, and $(\star\star$) if $s\in (0,1)$, then every $t>0$ can be written as $t=\sum_{n\geq 1}t_n$, with $t_n$ being a power of $s$, for all $n$.

Let us elaborate on the definition of LFBimod$(M)$. Any left--finite Hilbert $M$--bimodule $\Cal H$ comes from a unital embedding of $M$ into one of its amplifications ([C80a],  see also [Po86]). Indeed, if  $t=$ dim$_{M}(\Cal H)$, then $M'\cap \Bbb B(\Cal H)\cong (M^{op})^t$ and thus we obtain a unital $*$--homomorphism $\theta:M\rightarrow M^t$. Conversely, any unital $*$--homomorphism $\theta:M\rightarrow M^t$ induces a left--finite $M$--bimodule, which we denote $\Cal H_{\theta}$. For this, first represent $M^t=p(\Bbb M_n(C)\otimes M)p$, where $n\geq t$ is an integer and $p\in\Bbb M_n(\Bbb C)\otimes M$ has normalized trace equal to $\frac{t}{n}$. Then set $\Cal H_{\theta}=(\Bbb M_{1,n}(\Bbb C)\otimes L^2(M))p$ and define the left and right module actions by $x\cdot \xi=x\xi$ and $\xi\cdot x=\xi\theta(x)$.
Note also that if $\Cal H$ and $\Cal K$ are two left--finite $M$--bimodules then so is their {\it Connes tensor product}, $\Cal H\otimes_{M} \Cal K$ (see [Po86, 1.3.4.]). 
\vskip 0.05in

\vskip 0.05in
As a consequence of Theorem 8.1 and of results from [OP07] we obtain the first (partial) calculations of these invariants. We refer the reader to [OP07, 2.4.] for the definition of the {\it complete metric approximation property} (abbreviated CMAP).

\proclaim {10.5 Corollary}   Let $\Gamma_1$, $\Gamma_2$ be two ICC, torsion free groups, with  $\Gamma_2$  non--amenable.  Assume that $L(\Gamma_1)$, $L(\Gamma_2)$ have the CMAP (e.g. $\Gamma_1=\Bbb F_m,\Gamma_2=\Bbb F_n$, for $2\leq m,n\leq\infty$) and set $\Gamma=\Bbb\Gamma_1\times\Bbb \Gamma_2$. 
\vskip 0.01in
\noindent
Let $B$ be a non--trivial abelian von Neumann algebra. Let $\sigma:\Gamma\rightarrow$ Aut$(B^{\Gamma})$ be the Bernoulli action. Denote $M=B^{\Gamma}\rtimes\Gamma$ and $A=B^{\Gamma}$.

\vskip 0.05in
\noindent
If $\theta:M\rightarrow M^t$ is a unital $*$--homomorphism for some $t\leq 1$, then $t=1$. Moreover, there exist a character $\eta$ of $\Gamma$, an injective group morphism $\delta:\Gamma\rightarrow\Gamma$ and a unitary $u\in M$ such that $\theta(A)\subset uAu^*$ and $\theta(u_{\gamma})=\eta(\gamma)uu_{\delta(\gamma)}u^*$, for all $\gamma\in\Gamma$.
In particular,  every $\theta\in$ {\rm End}$(M)$ is irreducible, i.e. $\theta(M)'\cap M=\Bbb C1$.
\endproclaim

In other words, if $\Cal H$ is a Hilbert $M$--bimodule with dim$_{M}(\Cal H)\leq 1$, then we must have that dim$_{M}(\Cal H)=1$. Equivalently, $\Cal F_s(M)\cap (0,1)=\emptyset$, i.e. $M$ does not embed into $M^t$, for $t<1$. Moreover, Corollary 10.5 gives a description of End$(M)$. Roughly speaking, every endomorphism $\theta:M\rightarrow M$ is determined  by the following data:

\vskip 0.01in
\noindent
$\bullet$ a character $\eta$ of $\Gamma$, a group embedding $\delta:\Gamma\rightarrow\Gamma$ and

\noindent
$\bullet$ a ``realization" of $\sigma$ as a quotient of $\sigma_{|\delta(\Gamma)}$.
\vskip 0.03in

In view of Corollary 10.5, it seems reasonable to conjecture that there are II$_1$ factors $M$ for which LFBimod$(M)$ is trivial, i.e. such that every left--finite Hilbert $M$--bimodule  is isomorphic to $L^2(M)\otimes\ell^2_n$, for some $n\geq 1$, where $\ell^2_n$ is the  Hilbert space of dimension $n$. In this case, we would necessarily have that $\Cal F_s(M)=\Bbb N$ and End$(M)=$ Int$(M)$. To support our conjecture, we notice below that, for certain II$_1$ factors $M$,  LFBimod$(M)$ is ``countably generated" (see Remark 10.8).
\vskip 0.1in
\noindent
{\it Proof of Corollary 10.5}. By Theorem 8.1 (applied to the situation $N=M$) we only need to argue that $\theta(M)\nprec_{M}L(\Gamma).$ If this were not the case, then, since $L(\Gamma)$ has the CMAP and $M$ is a factor, we would get that $M\cong \theta(M)$ also has the CMAP. Notice that the restriction of $\sigma$ to $\Gamma_2$ can be seen as the Bernoulli action of $\Gamma_2$ on $\tilde B^{\Gamma_2}$, where $\tilde B=B^{\Gamma_1}$. Since $\Gamma_1$ is infinite we have that $\tilde B\cong L(\Bbb Z)$ and thus $M\supset B^{\Gamma}\rtimes_{\sigma_{|\Gamma_2}}\Gamma_2\cong L(\Bbb Z\wr\Gamma_2)$.  Altogether, we would derive that $L(\Bbb Z\wr\Gamma_2)$ has the CMAP. Since $\Gamma_2$ is non--amenable this contradicts  [OP07, Corollary 2.12.].

The irreducibility assertion follows easily from the first part. \hfill$\blacksquare$
\vskip 0.1in
\noindent 
(III) {\bf  Group measure space decompositions of rigid factors.}  In this section, we consider II$_1$ factors $M$ which have property (T), in the sense of Connes--Jones [CJ85]. Our main goal is to prove that $M$  admits at most countably many group measure space decompositions (see Theorem 10.7). We start with the following lemma:

\proclaim {10.6 Lemma} Let $M$ and $\Cal M$ be two separable II$_1$ factors, with $M$ having property (T). Let $\{\Delta_i:M\rightarrow\Cal M|i\in I\}$ be an uncountable family of not necessarily unital $*$--homomorphisms. 
\vskip 0.03in
\noindent
$\bullet$ Then there exists an uncountable subset $J$ of $I$ such that for every $i,j\in I$ we can find $0\not=v\in\Cal M$ satisfying $v=\Delta_i(1)v\Delta_j(1)$ and $\Delta_i(x)v=v\Delta_j(x)$, for all $x\in M$.
\vskip 0.03in
\noindent
$\bullet$  Assume moreover that  $\Delta_i$ is irreducible (i.e. $\Delta_i(M)'\cap \Delta_i(1)\Cal M\Delta_i(1)=\Bbb C\Delta_i(1)$), for every $i\in I$. Then  for all $i,j\in J$ we have that $\tau(\Delta_i(1))=\tau(\Delta_j(1))$ and that $\Delta_i=$ Ad$(u)\circ\Delta_j$, for some $u\in\Cal U(\Cal M)$. In particular, the set $\{\tau(\Delta_i(1))|i\in I\}$ is countable.
\endproclaim

This lemma is proven by a separability argument going back to [C80b] (see also  [Po86]). For completeness, let us sketch its proof.

\vskip 0.1in
\noindent
{\it Proof.} 
 Let $\varepsilon\in (0,\frac{1}{4})$. After replacing $I$ with an uncountable subset, we can assume that  $\tau(\Delta_i(1))\in ((1-\varepsilon)t,t)$, for all $i\in I$, for some $t\in (0,1)$. Let $p\in\Cal M$ be a projection of trace $t$. Since $\tau(\Delta_i(1))\leq \tau(p)$, we may assume that $\Delta_i(1)\leq p$, for all $i\in I$.

Next, for every $i,j\in I$, let $\Cal H_{i,j}=L^2(\Delta_i(1)\Cal M\Delta_j(1))$ and $\xi_{i,j}=\Delta_i(1)\Delta_j(1)\in\Cal H_{i,j}$. Endow $\Cal H_{i,j}$ with the Hilbert $M$--bimodule structure given by $x\cdot\xi\cdot y=\Delta_i(x)\xi\Delta_j(y)$. Then for every $x\in \Cal U(M),$ we have that $$||x\cdot\xi_{i,j}-\xi_{i,j}\cdot x||=||\Delta_i(x)\Delta_j(1)-\Delta_i(1)\Delta_j(x)||_2\leq $$ $$||\Delta_i(x)-\Delta_j(x)||_2+||\Delta_i(1)-p||_2+||\Delta_j(1)-p||_2\leq 2\sqrt{\varepsilon t}+||\Delta_i(x)-\Delta_j(x)||_2.$$ 

Also, it is easy to see that $||\xi_{i,j}||\geq \sqrt{t}(1-2\sqrt{\varepsilon})>0$. Thus, we have that $$\frac{ ||x\cdot\xi_{i,j}-\xi_{i,j}\cdot x||} {||\xi_{i,j}||}\leq \frac{2\sqrt{\varepsilon}}{1-2\sqrt{\varepsilon}}+\frac{||\Delta_i(x)-\Delta_j(x)||_2} {\sqrt{t}(1-2\sqrt{\varepsilon})},\forall x\in \Cal U(M).$$

Since $L^2(\Cal M)$ is a separable Hilbert space, this estimate implies that for every finite set $F\subset \Cal U(M)$  we can find an uncountable set $J\subset I$ such that $\frac{||x\cdot\xi_{i,j}-\xi_{i,j}\cdot x||}{||\xi_{i,j}||}\leq \frac{3\sqrt{\varepsilon}}{1-2\sqrt{\varepsilon}}$, for all $x\in F$ an every $i,j\in J$. Thus, if we choose $\varepsilon$  small enough, then property (T) quarantees that $\Cal H_{i,j}$ has a central vector, for every $i\not=j\in J$ ([CJ85],[Po06c]). In other words, we can find $v\in \Delta_i(1)\Cal M\Delta_j(1)$ such that $\Delta_i(x)v=v\Delta_j(x)$, for all $x\in M$. This proves the first assertion. 
The second assertion is  immediate once we notice that $vv^*\in \Delta_i(M)'\cap \Delta_i(1)\Cal M\Delta_i(1)$ and $v^*v\in\Delta_j(M)'\cap \Delta_j(1)\Cal M\Delta_j(1)$.
\hfill$\blacksquare$
\vskip 0.1in

We can now prove:

\proclaim {10.7 Theorem} Let $M$ be a property (T) II$_1$ factor. Then
\vskip 0.01in
\noindent
(1) There exist at most countably many non--conjugate, free ergodic p.m.p. actions $\Gamma\curvearrowright X$ such that $M\cong L^{\infty}(X)\rtimes\Gamma$. 

\noindent
(2) The set of $t>0$ such that $M^t$ is isomorphic to the von Neumann algebra of a free ergodic p.m.p. action is countable. 

\noindent
(3) The set of $t>0$ such that $M^t$ is isomorphic to a group von Neumann algebra is countable.

\endproclaim

Part (3) has been first proven in [Po07b, Section 4]. Below, we give a proof relying on a completely different argument. 
\vskip 0.05in
\noindent
{\it Proof.} (1) Assume by contradiction that there exists an uncountable set $I$  of mutually non--conjugate, free ergodic actions $\Gamma_{i}\curvearrowright X_{i}$ such that $M=L^{\infty}(X_i)\rtimes\Gamma_{i}$. For $i\in I$, define $\Delta_{i}:M\rightarrow M\overline{\otimes}M$ by $\Delta_{i}(au_{\gamma})=au_{\gamma}\otimes u_{\gamma}$, for all $a\in L^{\infty}(X_i)$ and every $\gamma\in\Gamma_i$. 

By Lemma 10.6 we can find an uncountable set $J\subset I$ such that for every $i,j\in J$, there exist $0\not= v\in M\overline{\otimes}M$ satisfying $\Delta_i(x)v=v\Delta_j(x)$, for all $x\in M$. Since $\Delta_i(a)=a$, for all $a\in L^{\infty}(X_i)$, we get that $\Delta_j(L^{\infty}(X_i))\prec_{M\overline{\otimes}M}M\otimes 1$. 

Lemma 9.2 then gives that $L^{\infty}(X_i)\prec_{M}L^{\infty}(X_j)$. Since $L^{\infty}(X_i)$ and $L^{\infty}(X_j)$ are Cartan subalgebras of $M$, [Po06c, Theorem A.1.] implies that they are unitarily conjugate. Therefore, the actions $\Gamma_i\curvearrowright X_i$ and $\Gamma_j\curvearrowright X_j$  are orbit equivalent ([Si55], [FM77]). Since $M=L^{\infty}(X_i)\rtimes\Gamma_i$ has property (T), $\Gamma_i$ also has property (T), for every $i\in I$. Altogether, we have found an uncountable family of mutually non--conjugate, orbit equivalent, free ergodic actions of property (T) groups. This contradicts [PV08b, Corollary 6.3.].

\vskip 0.05in
\noindent
(2) Since property (T) for II$_1$ factors is closed under amplifications, it suffices to show that the set $I=\{\tau(p)|$ $pMp$ is the von Neumann algebra of a free, ergodic action$\}$ is countable. Assume by contradiction that $I$ is uncountable. For every $t\in I$, let $p_t\in M$ be a projection of trace $t$ and let $\Gamma_t\curvearrowright X_t$ be a free, ergodic action such that $p_tMp_t=L^{\infty}(X_t)\rtimes\Gamma_t$.
Define $\theta_t:p_tMp_t\rightarrow p_tMp_t\overline{\otimes}p_tMp_t$ by $\theta_t(au_{\gamma})=au_{\gamma}\otimes u_{\gamma}$, for every $a\in L^{\infty}(X_t)$ and $\gamma\in\Gamma_t$. 

By taking  amplifications we get a $*$--homomorphism $\Delta_t:M\rightarrow M\overline{\otimes}M$ satisfying $\tau(\Delta_t(1))=t$. Since $I$ is uncountable, Lemma 10.6 gives an uncountable set $J\subset I$ such that for every $s,t\in J$ there is $0\not=v\in\Delta_t(1)(M\overline{\otimes}M)\Delta_s(1)$ satisfying $\Delta_t(x)v=v\Delta_s(x)$, for all $x\in M$.

By reasoning similarly to part (1), it follows that the actions $\Gamma_t\curvearrowright X_t$ and $\Gamma_s\curvearrowright X_s$ are stably orbit equivalent, for all $s,t\in J$.
More precisely, if  $\Cal R_t$ denotes the equivalence relation induced by  $\Gamma_t\curvearrowright X_t$, then  $\Cal R_t^{\frac{1}{t}}\cong \Cal R_s^{\frac{1}{s}},$ for all $s,t\in J$.  Here, $\Cal R^t$ denotes the $t$--amplification of an equivalence relation $\Cal R$ (for the definition, see e.g. [Io08, Section 4]). Thus, if we fix $t\in J$, then $\Cal R_t^{\frac{s}{t}}\cong\Cal R_s$, for every $s\in J$. On the other hand,  [Io08, Theorem 5.9.] gives that the set $\{q>0|\Cal R_t^q$ is induced by a free action of a countable group$\}$ is countable. Thus, $J$ must be countable, a contradiction.

\vskip 0.05in
\noindent
(3) It suffices to prove that $I=\{\tau(p)| pMp$ is a group von Neumann algebra$\}$ is countable. For $t\in I$, let $p_t\in M$ be a projection of trace $t$ and let $\Gamma_t$ be a group such that $p_tMp_t=L(\Gamma_t)$. Let $\theta_t:L(\Gamma_t)\rightarrow L(\Gamma_t)\overline{\otimes}L(\Gamma_t)$ be given by $\theta_t(u_{\gamma})=u_{\gamma}\otimes u_{\gamma}$, for each $\gamma\in\Gamma_t$. 

Since $M$ is a factor, $\Gamma_t$ is ICC. Thus, the set $\{(\gamma g\gamma^{-1},\gamma h\gamma^{-1})|\gamma\in\Gamma_t\}$ is infinite, for all $(g,h)\in(\Gamma_t\times\Gamma_t)\setminus\{(e,e)\}$. This fact implies that $\theta_t(L(\Gamma_t))'\cap L(\Gamma_t)\overline{\otimes}L(\Gamma_t)=\Bbb C1$, for all $t$. Since $p_tMp_t=L(\Gamma_t)$, by amplifying $\theta_t$ we get a $*$--homomorphism $\Delta_t:M\rightarrow M\overline{\otimes}M$ verifying $\tau(\Delta_t(1))=t$ and $\Delta_t(M)'\cap \Delta_t(1)(M\overline{\otimes}M)\Delta_t(1)=\Bbb C\Delta_t(1)$.

Finally, Lemma 10.6 implies that  $I=\{\tau(\Delta_t(1))|t\in I\}$ is countable. 
\hfill$\blacksquare$

\vskip 0.1in
\noindent{\bf 10.8 Remarks}. (1) Let $M$ be a II$_1$ factor. Recall that every left--finite $M$--bimodule is of the form $\Cal H_{\theta}$, for some unital $*$--homomorphism $\theta:M\rightarrow M^t$. Notice that $\Cal H_{\theta}$ is irreducible iff $\theta(M)'\cap M^t=\Bbb C1$. Lemma 10.6 thus implies that if $M$ has property (T), then there are only countably many  irreducible left--finite Hilbert $M$--bimodules (up to isomorphism).

\noindent
(2) A  II$_1$ factor $M$ is  {\it solid} if $A'\cap M$ is atomic, for any completely non--amenable von Neumann subalgebra $A\subset M$. N. Ozawa proved that $L(\Gamma)$ is solid, for any hyperbolic group $\Gamma$ ([Oz04]). 
If $M$ is a non--amenable solid factor, then  $\theta(M)'\cap M^t$ is atomic, for every unital $*$--homomorphism $\theta:M\rightarrow M^t$. This fact implies that  any left--finite  $M$--bimodule is the direct sum of countably many irreducible left--finite $M$--bimodules.

\noindent
(3) Let us explain how (1) and (2) imply that there are II$_1$ factors $M$ with ``few" left--finite bimodules.
Indeed, take $M=L(\Gamma)$, for some hyperbolic, property (T) group $\Gamma$. 
Combining (1) and (2) yields that there exist a countable family of Hilbert $M$--bimodules $\{\Cal H_n\}_{n\geq 1}$ such that every left--finite Hilbert $M$--bimodule $\Cal H$ is the direct sum of some of the $\Cal H_n$'s. 

\noindent
(4) Finally, note that N. Ozawa has very recently shown that if $\Gamma$ is a hyperbolic, property (T), ICC group (e.g. any ICC lattice $\Gamma<$ Sp$(1,n)$), then $L(\Gamma)$ is a property (T) II$_1$ factor which does not admit a group measure space decomposition ([Oz10]).

\vskip 0.1in
\noindent 
(IV) {\bf II$_1$ factors not isomorphic to twisted group von Neumann algebras.} We end the paper by noticing that an extension of the results from Section 8 can be used to give examples of II$_1$ factors that are not isomorphic to twisted group von Neumann algebras.
Let us begin by recalling the construction of the von Neumann algebra $L_{\alpha}(G)$ arising from a  countable group $G$ and a 2--cocycle $\alpha\in$ H$^2(G,\Bbb T)$, i.e. a map $\alpha:G\times G\rightarrow \Bbb T$ satisfying $(\diamond)$\hskip 0.05in $\alpha(g,h)\alpha(gh,k)=\alpha(g,hk)\alpha(h,k)$, for all $g,h,k\in G$. 

First,  the formula

$$u_g^{\alpha}(\delta_h)=\alpha(g,h)\delta_{gh},\hskip 0.05in\forall g,h\in G$$
defines a projective unitary representation $u^{\alpha}:G\rightarrow \Cal U(\ell^2G)$, where $\{\delta_{h}\}_{h\in G}$ is the usual orthonormal basis of $\ell^2G$. More precisely, $u_g^{\alpha}u_h^{\alpha}=\alpha(g,h)u_{gh}^{\alpha}$, for all $g,h\in G$.

 Then $L_{\alpha}(G)$ is defined as the von Neumann algebra generated by $\{u_g^{\alpha}\}_{g\in G}$. Note that  $\tau:L_{\alpha}(G)\rightarrow\Bbb C$ given by $\tau(u_g^{\alpha})=\delta_{g,e}\alpha(e,e)$ is a faithful normal trace. 

We continue with two useful facts about twisted group von Neumann algebras.

\proclaim {10.9 Lemma} Let $G$ be a countable group and $\alpha\in$ H$^2(G,\Bbb T)$ be a 2--cocycle. Let $\overline{\alpha}\in$ H$^2(G,\Bbb T)$ be given by $\overline{\alpha}(g,h)=\overline{\alpha(g,h)}$, for all $g,h\in G$. Then 

\vskip 0.05in
\noindent
(1) The opposite von Neumann algebra $L_{\alpha}(G)^{op}$ is isomorphic to $L_{\overline{\alpha}}(G)$.

\noindent
(2) The map $\Delta:G\rightarrow L_{\alpha}(G)\overline{\otimes}L_{\alpha}(G)\overline{\otimes}L_{\overline{\alpha}}(G)$ given by $\Delta(u_g^{\alpha})=u_g^{\alpha}\otimes u_g^{\alpha}\otimes u_g^{\overline{\alpha}}$, for every $g\in G$, extends to a unital $*$--homomorphism $\Delta:L_{\alpha}(G)\rightarrow L_{\alpha}(G)\overline{\otimes}L_{\alpha}(G)\overline{\otimes}L_{\overline{\alpha}}(G)$.
\endproclaim
\noindent
{\it Proof.} (1) Denote $\beta=\alpha(e,e)$ and notice that the cocycle identity $(\diamond)$ implies that $\alpha(e,g)=\beta$, for all $g\in G$.

Now, for every $g\in G$,  define $\theta(u_g^{\alpha})=\beta\alpha(g,g^{-1})u_{g^{-1}}^{\overline{\alpha}}\in L_{\overline{\alpha}}(G)$. 
We claim that $\theta$ extends to a $*$--isomorphism $\theta:L_{\alpha}(G)^{op}\rightarrow L_{\overline{\alpha}}(G)$. To see this, it suffices to prove that  $\theta$ is trace preserving and that $\theta(u_h^{\alpha}u_g^{\alpha})=\theta(u_g^{\alpha})\theta(u_h^{\alpha})$, for all $g,h\in G$. Fix $g,h\in G$.

The first assertion holds because $\tau(\theta(u_g^{\alpha}))=\beta\alpha(g,g^{-1})\delta_{g,e}\overline{\alpha(e,e)}=\delta_{g,e}\beta=\tau(u_g^{\alpha}).$

For the second assertion, using the cocycle identity twice yields that $$\alpha(h,g)\alpha(hg,g^{-1}h^{-1})\alpha(g^{-1},h^{-1})=\alpha(h,h^{-1})\alpha(g,g^{-1}h^{-1})\alpha(g^{-1},h^{-1})=\tag 10.a$$ $$=\alpha(e,h^{-1})\alpha(g,g^{-1})\alpha(h,h^{-1}) $$
Since $\theta(u_h^{\alpha}u_g^{\alpha})=\alpha(h,g)\theta(u_{hg}^{\alpha})=\beta\alpha(h,g)\alpha(hg,g^{-1}h^{-1})u_{g^{-1}h^{-1}}^{\overline{\alpha}}$, while $\theta(u_g^{\alpha})\theta(u_h^{\alpha})=\beta^2\alpha(g,g^{-1})\alpha(h,h^{-1})u_{g^{-1}}^{\overline{\alpha}}u_{h^{-1}}^{\overline{\alpha}}=\beta^2\alpha(g,g^{-1})\alpha(h,h^{-1})\overline{\alpha(g^{-1},h^{-1})}u_{g^{-1}h^{-1}}^{\overline{\alpha}},$ the second assertion follows by using (10.a).

\noindent
(2) Since $\Delta_{|G}$ is clearly multiplicative and trace preserving, the conclusion follows.
\hfill$\blacksquare$
\vskip 0.05in

Next, let us point out a generalization of the results from Section 8.
Let $\Gamma$ be a torsion free, ICC group which admits an infinite almost normal subgroup such that the inclusion $(\Gamma_0\subset\Gamma)$ has the relative property (T). Let $B$ be a non--trivial abelian von Neumann algebra and define $M=B^{\Gamma}\rtimes\Gamma$.

 Theorems 8.1 and 8.2 classify  embeddings of $M$ into itself and into $M\overline{\otimes}M$, repsectively. A straightforward modification of the proof of Theorem 8.2 shows that the embeddings of $M$ into $\overline{\otimes}_{i=1}^nM$ can be classified in a similar way, for every $n\geqslant 1$. 

\proclaim {10.10 Theorem} Let $M=B^{\Gamma}\rtimes\Gamma$ be as above and denote $A=B^{\Gamma}$. 

\noindent
Let  $\Delta:M\rightarrow M^S$ be a  $*$--homomorphism, for some finite set $S$, and
suppose that:

\noindent
(1) $\Delta(L(\Gamma_0))\nprec_{M^S}L(\Gamma)^{S\setminus\{s\}}$, for any $s\in S$.

\noindent
(2) $\Delta(M)\nprec_{M^S}L(\Gamma)^{\{s\}}\overline{\otimes}M^{S\setminus\{s\}}$, for any $s\in S$.
\vskip 0.05in
\noindent
Then $\Delta$ must be unital and we can find a character $\eta$ of $\Gamma$, group homomorphisms $\delta_s:\Gamma\rightarrow\Gamma$, for every $s\in S$, and a unitary $u\in M^{S}$ such that $u\Delta(A)u^*\subset A^{S}$ and $u\Delta(u_{\gamma})u^*=\eta(\gamma)(\otimes_{s\in S}u_{\delta_s(\gamma)})$, for every $\gamma\in\Gamma$. 
\endproclaim

Here we denote by $M^S$ the tensor product algebra $\overline{\otimes}_{s\in S}(M)_s$. For a subset $S'\subset S$ and a von Neumann subalgebra $Q\subset M$, we view $Q^{S'}$ as a von Neumann subalgebra of $M^S$, in the obvious way.

\vskip 0.05in
Finally, let us combine the last two results to give examples of II$_1$ factors which are not isomorphic to twisted group von Neumann algebras.

\proclaim {10.11 Corollary} Let $M$ be as above and $p\in M\setminus\{1\}$ be a projection. Then  $pMp$ is not isomorphic to $L_{\alpha}(G)$, for any countable group $G$ and any $\alpha\in$ H$^2(G,\Bbb T)$.

\endproclaim
\noindent
{\it Proof.} Assume that $pMp$ can be written as $L_{\alpha}(G)$, for some projection $p\in M$. Our goal is to prove that $p=1$. First, by Remark 10.3 (3), we have that $pMp$ is anti--isomorphic to itself, i.e. $(pMp)^{op}\cong pMp$. On the other hand, part (1) of Lemma 10.9 gives that $L_{\alpha}(G)^{op}\cong L_{\overline{\alpha}}(G)$. Altogether, we deduce the existence of a $*$--isomorphism $\theta:L_{\overline{\alpha}}(G)\rightarrow L_{\alpha}(G)$.

Next, let $S=\{1,2,3\}$ and consider the $*$--homomorphism $$\Delta_1:=(id\otimes id\otimes\theta)\circ\Delta:L_{\alpha}(G)\rightarrow L_{\alpha}(G)^S=\overline{\otimes}_{i=1}^3L_{\alpha}(G),$$ where $\Delta$ is as  in Lemma 10.9 (2).
Concretely, we have that $\Delta_1(u_g^{\alpha})=u_g^{\alpha}\otimes u_g^{\alpha}\otimes\theta(u_g^{\overline{\alpha}})$, for all $g\in G$.
To simplify notation, denote $N=L_{\alpha}(G)$. Since $\Delta_1$ comes from a ``diagonal embedding" of $G$, we have the following:

\vskip 0.05in
\noindent
{\it Claim.} Let $Q\subset N$ be a von Neumann subalgebra. 

\noindent
(1) If $\Delta_1(Q)\prec_{N^S}N^{S\setminus\{s\}}$,  for some $s\in S$, then $Q$ is not diffuse.

\noindent
(2) If $\Delta_1(N)\prec_{N^S}Q^{\{s\}}\overline{\otimes}N^{S\setminus\{s\}}$, then $N\prec_{N}Q$.

\vskip 0.05in
The proof of this Claim is  analogous to the proof of Lemma 9.2 and so we leave the details to the reader. 
\vskip 0.05in

 Now, since $N=pMp$, by amplifying $\Delta_1:N\rightarrow N^S$, we get a $*$--homomorphism $\Delta_2:M\rightarrow M^S$, such that $\tau(\Delta_2(1))=\tau(p)^2.$ By using the Claim it is immediate that $\Delta_2$ verifies conditions (1) and (2) in Theorem 10.10. Thus, by applying Theorem 10.10 it follows that $\Delta_2$ is unital, hence $p=1$.\hfill$\blacksquare$

\head References\endhead
\item {[BHV08]} B. Bekka, P. de la Harpe, A. Valette: {\it Kazhdan's property (T),} New Mathematical
Monographs, 11. Cambridge University Press, Cambridge, 2008.
\item {[BrOz08]} N.P. Brown, N. Ozawa: {\it C$^*$--algebras and finite-dimensional approximations,} Graduate
Studies in Mathematics, 88. American Mathematical Society, Providence, RI,
2008. xvi+509 pp.
\item {[Bu91]} M. Burger: {\it Kazhdan constants for} SL$(3,\Bbb Z)$, J. Reine Angew. Math. {\bf 413} (1991), 36--67.
\item {[C75]} A. Connes: {\it Sur la classification des facteurs de type II}, C. R. Acad. Sci. Paris 281
(1975), 13-15.
\item {[C80a]} A. Connes: {\it Correspondences}, handwritten notes, 1980.
\item {[C80b]} A. Connes: {\it A factor of type II$_1$ with countable fundamental group},  J. Operator Theory  {\bf 4} (1980), no. 1, 151--153. 
\item {[CFW81]} A. Connes, J. Feldman, B. Weiss: {\it An amenable equivalence relation is generated
by a single transformation},  Ergodic. Th. and Dynam. Sys {\bf 1} (1981), no. 4, 431--450.
\item {[CI08]} I. Chifan, A. Ioana: {\it Ergodic Subequivalence Relations Induced by a Bernoulli Action}, 
Geom. Funct. Anal. Vol. {\bf 20} (2010) 53–-67.
\item {[CJ82]} A. Connes, V.F.R. Jones: {\it  A II$_1$ factor with two non--conjugate Cartan subalgebras,} Bull.
Amer. Math. Soc. {\bf 6} (1982), 211--212.
\item {[CJ85]} A. Connes, V.F.R. Jones: {\it Property (T) for von Neumann algebras}, Bull. London Math.
Soc. {\bf 17} (1985), 57--62.
\item {[Fu99]} A. Furman: {\it Orbit equivalence rigidity,} Ann. of Math. {\bf 150} (1999), 1083--1108.
\item {[Fu09]} A. Furman: {\it A survey of Measured Group Theory,} Geometry, Rigidity, and Group Actions, 296--374, The University of Chicago Press, Chicago and London, 2011, available at arXiv:0901.0678.
\item {[FM77]} J. Feldman, C.C. Moore: {\it Ergodic equivalence relations, cohomology, and von Neumann algebras, II}, Trans. Amer. Math. Soc. {\bf 234} (1977), 325--359.
\item {[Io07]} A. Ioana: {\it Rigidity results for wreath product II$_1$ factors}, 
J. Funct. Anal. {\bf 25} (2007) 763--791.
\item {[Io08]} A. Ioana: {\it Cocycle superrigidity for profinite actions of property (T) groups}, 
Duke Math. J. Volume {\bf 157}, Number 2 (2011), 337--367.
\item {[Io09]} A. Ioana: {\it Non-orbit equivalent actions of $\Bbb F_n$}, 
Ann. Sci. \'Ec. Norm. Sup\'er., {\bf 42}, fascicule 4 (2009), 675--696. 
\item {[IPP08]} A. Ioana, J. Peterson, S. Popa: {\it Amalgamated free products of weakly rigid factors
and calculation of their symmetry groups}, Acta Math. {\bf 200} (2008), no. 1, 85-–153.
\item {[J080]} V.F.R. Jones: {\it A II$_1$ factor anti-isomorphic to itself but without involutory 
antiautomorphisms}, Math. Scand. {\bf 46} (1980), no. 1, 103--117.
\item {[KR97] } R.V. Kadison, J.R. Ringrose: {\it Fundamentals of the theory of operator algebras.
Vol. II,} Academic Press, Orlando, 1986.
\item {[Ka67]} D. Kazhdan: {\it Connection of the dual space of a group with the structure of its
closed subgroups, Funct. Anal. and its Appl.}, {\bf 1} (1967), 63--65.
\item {[Ki06]} Y. Kida: {\it Measure equivalence rigidity of the mapping class group,}   Ann. of Math.  {\bf 171} (2010), No. 3, 1851–-1901.                         
\item {[Ki09]} Y. Kida: {\it Rigidity of amalgamated free products in measure equivalence theory}, to appear in J. Topology, preprint arXiv:0902.2888.
\item {[Ma82]} G. Margulis: {\it Finitely-additive invariant measures on Euclidian spaces}, Ergodic.
Th. and Dynam. Sys. {\bf 2} (1982), 383--396.
\item {[MvN36]} F. Murray, J. von Neumann: {\it On rings of operators,} Ann. of Math. {\bf 37} (1936), 116--229.
\item {[MvN43]} F. Murray, J. von Neumann: {\it Rings of operators IV,} Ann. Math. {\bf 44} (1943), 716--808.
\item {[OW80]} D. Ornstein, B. Weiss: {\it Ergodic theory of amenable groups. I. The Rokhlin lemma.},
Bull. Amer. Math. Soc. (N.S.) {\bf 1} (1980), 161-–164.
\item {[Oz04]} N. Ozawa: {\it Solid von Neumann algebras},  Acta Math.  {\bf 192}  (2004),  no. 1, 111--117.
\item {[OP07]} N. Ozawa, S. Popa: {\it On a class of II$_1$ factors with at most one Cartan subalgebra},
Ann. Math. Vol. {\bf 172} (2010), No. 1, 713–749. 
\item {[Oz10]} N. Ozawa: {\it Examples of groups which are not weakly amenable}, 

 preprint arXiv:1012.0613.
\item {[Pe09]} J. Peterson: {\it $L^2$-rigidity in von Neumann algebras},  Invent. Math.  {\bf 175}  (2009),  no. 2, 417--433.
\item {[Pe10]} J. Peterson: {\it Examples of group actions which are virtually W*-superrigid}, 
 preprint arXiv:1002.1745.
\item {[Po86]} S. Popa: {\it Correspondences}, INCREST preprint 1986, unpublished.
\item {[Po06a]} S. Popa: {\it Strong Rigidity of II$_1$ Factors Arising from Malleable Actions of w-Rigid Groups. I.}, Invent. Math. {\bf 165} (2006), 369--408.
\item  {[Po06b]} S. Popa: {\it Strong Rigidity of II$_1$ Factors Arising from Malleable Actions of w-Rigid Groups. II.}, Invent. Math. {\bf 165} (2006), 409--451.
\item {[Po06c]} S. Popa: {\it On a class of type II$_1$ factors with Betti numbers invariants}, Ann. of Math. {\bf 163} (2006), 809--889. 
\item {[Po06d]} S. Popa: {\it Some rigidity results for non-commutative Bernoulli shifts}, J. Funct. Anal.
 {\bf 230} (2006), 273--328.
\item {[Po07a]} S. Popa: {\it Cocycle and orbit equivalence superrigidity for malleable actions of $w$-rigid groups},  Invent. Math.  {\bf 170}  (2007),  no. 2, 243--295.
\item {[Po07b]} S. Popa: {\it Deformation and rigidity for group 	
actions and von Neumann algebras},  International Congress of Mathematicians. 
Vol. I,  445--477, 
Eur. Math. Soc., Z$\ddot{\text{u}}$rich, 2007.
\item {[Po08]} S. Popa: {\it On the superrigidity of malleable actions with spectral gap}, 
J. Amer. Math. Soc. {\bf 21} (2008), 981--1000. 
\item {[PV08a]} S. Popa, S. Vaes: {\it  Strong rigidity of generalized Bernoulli actions and computations of their 
symmetry groups.} Adv. Math. {\bf 217} (2008),  no. 2, 833--872.
\item {[PV08b]} S. Popa, S. Vaes: {\it On the fundamental group of II$_1$ factors and equivalence relations
 arising from group actions}, Quanta of maths,  519–-541, Clay Math. Proc., 11, Amer. Math. Soc., Providence, RI, 2010.
\item {[PV09]} S. Popa, S. Vaes: {\it Group measure space decomposition of II$_1$ factors and} 

{\it W*-superrigidity}, Invent. Math. {\bf 182} (2010), no. 2, 371--417.
\item {[Si55]} I.M. Singer: {\it Automorphisms of finite factors}, Amer. J. Math. {\bf 77} (1955), 117--133.
\item {[Va07]} S. Vaes: {\it Rigidity results for Bernoulli actions and their von Neumann algebras (after Sorin Popa)}, S\'eminaire Bourbaki. Vol. 2005/2006.  Ast\'erisque  No. {\bf 311}  (2007), Exp. No. 961, viii, 237--294.
\item {[Va08]} S. Vaes: {\it Explicit computations of all finite index bimodules for a family of II$_1$ factors}, Ann. Sci. \'Ec. Norm. Sup\'er. (4) {\bf 41} (2008), no. 5, 743--788.
\item {[Va09]} S. Vaes: {\it Factors of type II$_1$ without non-trivial finite index subfactors,} Trans. Amer. Math. Soc. {\bf 361} (2009), no. 5, 2587--2606.

\enddocument